\documentclass[a4paper,dutch,reqno]{article}
\usepackage[margin=1.2in]{geometry}
\usepackage[english]{babel}
\usepackage{graphicx}
\usepackage{amsmath}
\usepackage{amsfonts}
\usepackage{amssymb}
\usepackage{amsthm}
\usepackage{amsfonts}
\usepackage{mathtools}
\usepackage{tabu}
\usepackage{verbatim}
\usepackage{float}
\usepackage{tikz-cd}
\usepackage[utf8x]{inputenc}
\usepackage{enumitem}
\usetikzlibrary{positioning}
\usetikzlibrary{matrix}
\usepackage{ wasysym }
\usepackage{wrapfig}
\usepackage{blkarray}
\usepackage[toc,page]{appendix}
\usepackage{hyperref}

\DeclareUnicodeCharacter{120779}{\gammagamma}

\usepackage[resetcount,linesnumbered,titlenotnumbered,noend]{algorithm2e}

\usepackage{bbm}

\theoremstyle{plain}
\newtheorem{stelling}{Stelling}[section]

\newtheorem{lemma}[stelling]{Lemma}

\newtheorem{proposition}[stelling]{Proposition}
\newtheorem{corollary}[stelling]{Corollary}
\newtheorem{theorem}[stelling]{Theorem}

\theoremstyle{definition}

\newtheorem{definition}[stelling]{Definition}
\newtheorem{remark}[stelling]{Remark}
\newtheorem{example}[stelling]{Example}
\newtheorem{Algorithm}[stelling]{Algorithm}

\newcommand{\Z}		{\ensuremath{\mathbb{Z}}}
\newcommand{\Q}		{\ensuremath{\mathbb{Q}}}
\newcommand{\R}		{\ensuremath{\mathbb{R}}}
\newcommand{\C}		{\ensuremath{\mathbb{C}}}
\newcommand{\F}		{\ensuremath{\mathbb{F}}}
\newcommand{\m}		{\ensuremath{\mathfrak{m}}}
\newcommand{\n}		{\ensuremath{\mathfrak{n}}}
\newcommand{\p}		{\ensuremath{\mathfrak{p}}}
\newcommand{\q}		{\ensuremath{\mathfrak{q}}}

\newcommand{\mapsfrom}{\mathrel{\reflectbox{\ensuremath{\mapsto}}}}

\newcommand{\tensor} {\otimes}

\DeclareMathOperator{\Hom}{Hom}
\DeclareMathOperator{\Aut}{Aut}

\DeclareMathOperator{\divs}{|}

\DeclareMathOperator{\supp}{supp}
\DeclareMathOperator{\id}{id}
\DeclareMathOperator{\im}{im}
\DeclareMathOperator{\ord}{ord}
\DeclareMathOperator{\lcm}{lcm}
\DeclareMathOperator{\spec}{spec}

\DeclareMathOperator{\wreath}{wr}
\DeclareMathOperator{\bauble}{ba}
\DeclareMathOperator{\obj}{obj}
\DeclareMathOperator{\Iso}{Iso}

\DeclareMathOperator{\minspec}{minspec}
\DeclareMathOperator{\maxspec}{maxspec}

\DeclareMathOperator{\aut}{aut}
\DeclareMathOperator{\iso}{iso}
\DeclareMathOperator{\End}{End}
\DeclareMathOperator{\prid}{prid}
\DeclareMathOperator{\sqrfree}{rad}

\DeclareMathOperator{\Tr}{Tr}

\newcommand{\Xe} {X_{e}}
\newcommand{\Xew} {X_{e}^{!}}
\newcommand{\Xet} {X_{e}^{!!}}

\newcommand{\catt}[1] {\texttt{\textup{#1}}}
\renewcommand{\Upsilon} {\textup{Y}}

\newcommand{\us} {\rule[.2pt]{.70em}{.5pt}}
\newlength{\dhatheight}
\newcommand{\widehhat}[1] {\settoheight{\dhatheight}{\ensuremath{\widehat{#1}}}%
    \addtolength{\dhatheight}{-0.35ex}%
    \widehat{\vphantom{\rule{1pt}{\dhatheight}}%
    \smash{\widehat{#1\vphantom{f}}}}}
\newcommand{\hatf} { \widehat{\,\cdot\,} }
\newcommand{\hhatf} { \widehhat{\,\cdot\,} }
\DeclareMathOperator{\ev}{ev}
\DeclareMathOperator{\characteristic}{char}

\allowdisplaybreaks[3]

\begin{document}

\begin{titlepage}
\begin{center}

{\Large\bf 
D.M.H. van Gent
} 

\vspace{1em} 

{\LARGE\bf 
Algorithms for finding \\ the gradings of reduced rings
} 

\vspace{10em} 

{\large\bf 
Master thesis
} 

\vspace{1em}

{\large\bf 19 June 2019
}

\vspace{10em} 

{\large\bf
\begin{tabular}{ll}
Thesis supervisor: & Prof. Dr. H.W. Lenstra Jr
\end{tabular}
}

\vfill

\includegraphics{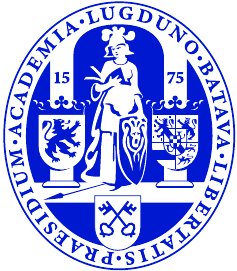}\\

\vspace{2em}

{\large\bf 
Leiden University\\
Mathematical Institute\\
}

\end{center}
\end{titlepage}

\tableofcontents
\clearpage

\section{Introduction}\label{sec:introduction}

In this thesis we will study gradings of \(\Q\)-algebras and {\em orders}, which are commutative rings for which the additive group is isomorphic to \(\Z^r\) for some \(r\in\Z_{\geq 0}\).
An element \(x\) of a ring \(R\) is called {\em nilpotent} if \(x^n=0\) for some integer \(n\geq 1\), and \(R\) is called {\em reduced} if \(0\) is the only nilpotent of \(R\).
In \cite{Lenstra2018}, H.W. Lenstra and A. Silverberg showed that each reduced order has a universal grading, for which we also give an alternative proof (Theorem~\ref{thm:pgood_theorem}).
Our goal is to present an algorithm to compute this universal grading in polynomial time once we bound the size of the minimal spectrum of the order.

\begin{definition}\label{def:decomposition}
Let \(k\) be a commutative ring.
A {\em decomposition} of a \(k\)-module \(M\) is a family of \(k\)-submodules \(\{M_s\}_{s\in S}\) indexed by some set \(S\) such that the natural map \(\bigoplus_{s\in S} M_s\to M\) is an isomorphism.
For any decomposition \(\{M_s\}_{s\in S}\) and map \(f:S\to T\) we define \(f_*(\{M_s\}_{s\in S}) = \{ \sum_{s\in f^{-1} t} M_s \}_{t\in T}\), which is again a decomposition of \(M\).
We say a decomposition \(\mathcal{M}\) is a {\em refinement} of a decomposition \(\mathcal{N}\), written \(\mathcal{N}\preceq\mathcal{M}\), if there exists a map \(f\) such that \(f_*\mathcal{M}=\mathcal{N}\).
A {\em pre-grading} of a \(k\)-algebra \(R\) is a decomposition \(\{R_s\}_{s\in S}\) of \(R\) as \(k\)-module such that for all \(s,t\in S\) there exists some \(u\in S\) such that \(R_s\cdot R_t \subseteq R_u\).
\end{definition}

With \(R=k[X]\), an example of a pre-grading would be the decomposition \(\{k X^i\}_{i\in \Z_{\geq 0}}\), since \(kX^i \cdot k X^j \subseteq k X^{i+j}\).
For the order \(\Z[\sqrt{2}]\) one can take \(\{\Z\sqrt{2^i}\}_{i\in \Z/2\Z}\), where \(\Z\sqrt{2^i}\) is well defined since \(\sqrt{2^2}\in\Z\).
In both examples, the index set \(S\) comes equipped with a binary operator \(*:S^2\to S\) satisfying \(R_s\cdot R_t \subseteq R_{s*t}\).
This is not a surprise, since if \(R_s\cdot R_t\neq 0\), then any \(u\in S\) such that \(R_s\cdot R_t\subseteq R_u\) must be unique if it exists. 

\begin{definition}
Let \(k\) be a commutative ring and \(R\) be a (not necessarily commutative) \(k\)-algebra.
A {\em group-grading of \(R\)} is a triple \(\overline{R}=(R,G,\mathcal{R})\) with \(G\) a group and \(\mathcal{R}=\{R_g\}_{g\in G}\) a \(G\)-indexed pre-grading of \(R\), such that \(1\in R_1\) and for all \(g,h\in G\) we have \(R_g\cdot R_h \subseteq R_{gh}\).
For brevity we call \(\overline{R}\) a {\em\(G\)-grading} when \(\overline{R}\) is a group-grading with the group \(G\) as index set.
We call a \(G\)-grading an {\em abelian group-grading} if \(G\) is an abelian group.
For any morphism of groups \(\varphi:G\to H\) we define \(\varphi_*(\overline{R})=(R,H,\varphi_*(\mathcal{R}))\), which is again a grading.
The class of gradings of \(R\) forms a category (see Appendix A) where the {\em morphisms} \(\varphi:\overline{R}\to\overline{S}\) are simply the morphisms of groups such that \(\varphi_*(\overline{R})=\overline{S}\).
\end{definition}

The pre-grading \(\{\Z\sqrt{2^i}\}_{i\in \Z/2\Z}\) of \(\Z[\sqrt{2}]\) trivially becomes a \(\Z/2\Z\)-grading.
For \(k[X]\) we may take \(R_i=k X^i\) when \(0\leq i\) and \(R_i=0\) when \(i<0\) to turn \(\{R_i\}_{i\in\Z}\) into a \(\Z\)-grading.
We say a group-grading \(\overline{U}\) of \(R\) is {\em universal} if it is an initial object in the category of group-gradings of \(R\).
Concretely, this means that \(\overline{U}\) is the finest grading of \(R\) possible.
By Yoneda's lemma, having a universal grading is equivalent to having a representation of the functor \(F:\catt{Grp}\to\catt{Set}\) that sends each group \(G\) to the set of \(G\)-gradings of \(R\).
We analogously define a universal abelian group-grading by considering the category of abelian group gradings of \(R\) instead.
We cannot expect in general that a universal grading exists (Example~\ref{intersection_counter}).
However, for reduced orders a universal abelian group-grading always exists \cite{Lenstra2018}.
Since the grading of \(\Z[\sqrt{2}]\) we gave earlier cannot be further refined, we must conclude that it is a universal abelian group-grading. 

For our algorithms, \(k\) will either be \(\Z\) or \(\Q\), elements of which have an obvious binary encoding.
Additionally, \(R\) will either be a reduced order or reduced commutative \(\Q\)-algebra, which both are free \(k\)-modules.
Elements of \(k^n\) are simply represented by a list of \(n\) elements of \(k\) and a morphism of \(k\)-modules \(f:k^m\to k^n\) is represented as a matrix in \(k^{n\times m}\).
Submodules \(M\) of \(k^n\) are represented by any injection \(f:k^m\to k^n\) with image \(M\).
For \(k\)-algebras free of rank \(n\) as \(k\)-module we additionally encode the multiplicative structure by structure constants \((a_{hij})_{hij}\in k^{n\times n\times n}\) such that \(e_h\cdot e_i := \sum_{j=1}^n a_{hij} e_j\) for all \(1\leq h,i\leq n\), with \(e_i\) the \(i\)-th standard basis vector of \(k^n\).
We encode a finite abelian group by its entire multiplication table and we then get an obvious encoding for our gradings.
Now that we have specified our encoding we can formally talk about algorithmic complexity.

We call a prime ideal \(\p\) of a commutative ring \(R\) {\em minimal} if for all prime ideals \(\q\subseteq \p\) we have \(\p=\q\), and write \(\minspec R\) for the set of minimal prime ideals of \(R\).

\begin{theorem}\label{AlgAbGrpThm}
There is a deterministic algorithm that takes a pair \((E,q)\) as input, where \(E\) is a reduced commutative \(\Q\)-algebra and \(q\) is a prime power, and produces all \(\Z/q\Z\)-gradings of \(E\) in time \(n^{O(m)}\), where \(n\) is the length of the input and \(m=\#\minspec E\).
\end{theorem}

Our main result is the following.

\begin{theorem}\label{AlgGridThm}
There is a deterministic algorithm that takes a reduced order \(R\) as input and produces a universal abelian group grading of \(R\) in time \(n^{O(m)}\), where \(n\) is the length of the input and \(m=\#\minspec R\).
\end{theorem}

If we bound the size of the minimal spectrum by a constant, the above algorithms run in polynomial time. 
For example, when \(R\) is a domain, or equivalently an order in a number field, then \(\#\minspec R = 1\) and thus is the algorithm polynomial.
For the proof of Theorem~\ref{AlgAbGrpThm} and Theorem~\ref{AlgGridThm}, see Section~\ref{sec:AlgAbGrpThm} respectively Section~\ref{sec:AlgGridThm}.

\begin{example}\label{ex:why_grid}
One might be tempted to think that all group-gradings of a commutative ring are abelian group-gradings.
Indeed, for a \(G\)-grading \(\overline{R}\) we have \(R_g \cdot R_h \subseteq R_{gh} \cap R_{hg}\) for all \(g,h\in G\), so if \(R_g\cdot R_h \neq 0\) we have \(R_{gh}=R_{hg}\neq 0\) and thus \(gh=hg\).
If \(R\) is not necessarily a domain, say \(R=\Z[\sqrt{2}]^2\), we may provide a counterexample.
Let \(G=\Aut_\catt{Set}(\{1,2,3\})\) be the symmetric group of degree 3 and define \(R_1=\Z\times \Z\), \(R_{(1\,2)}=(\Z\sqrt{2})\times 0\), \(R_{(2\, 3)}=0\times(\Z\sqrt{2})\) and \(R_g=0\) for all other \(g\in G\).
Then \((R,G,\{R_g\}_{g\in G})\) is a non-abelian group-grading of a commutative ring.  
We can give a \((\Z/2\Z)^2\)-grading \(\overline{S}\) of \(R\) with \(S_{(0,0)}=R_1\), \(S_{(1,0)}=R_{(1\, 2)}\), \(S_{(0,1)}=R_{(2\, 3)}\) and \(S_{(1,1)}=0\) that has the same non-zero components as \(\overline{R}\) but has an abelian group.
Moreover, there exists no morphism \(f\) between the groups \(G\) and \((\Z/2\Z)^2\) in either direction such that \(f_*\overline{R}=\overline{S}\) or \(f_*\overline{S}=\overline{R}\).
\end{example}

We consider Example~\ref{ex:why_grid} as evidence that there perhaps is a structure that is better suited to grade a ring than a group.
Another piece of evidence is Example~\ref{ex:grid_not_group}.

\begin{definition}
A {\em grid} is a quadruple \(G=(S,1,D,*)\), where \(S\) is a set containing \(1\), such that \(D\subseteq S^2\) satisfies \((1,s),(s,1)\in D\) for all \(s\in S\) and \(*:D\to S\) is a binary operator such that for all \(s\in S\) we have \(s*1=1*s=s\), and \(s*s=s\) implies \(s=1\). 
Through abuse of notation we will often write \(g\in G\) where we mean \(g\in S\).
We say \(G\) is {\em abelian} when for all \((s,t)\in D\) we have \((t,s)\in D\) and \(s*t=t*s\).
With \(H=(T,1,E,\star)\) another grid, a {\em morphism} from \(G\) to \(H\) is a map \(f:S\to T\) such that for all \((s,t)\in D\) we have \((f(s),f(t))\in E\) and \(f(s*t)=f(s)\star f(t)\).
\end{definition}

It follows immediately that \(f(1)=1\) for all grid morphisms \(f\) from the fact that \(1\) is the only idempotent.
The class of grids forms a category \catt{Grd}, which contains the category of groups as a full subcategory.

\addtocounter{stelling}{-2}
\begin{example}[continued]
An example of a grid would be the grid \(V\) on the elements \(\{1,a,b\}\) such that \(a*a=b*b=1\) and \(a*b\) and \(b*a\) are undefined.
One can informally view this grid as \((\Z/2\Z)^2=\{1,a,b,ab\}\) with \(ab\) and all products resulting in \(ab\) removed.
We have grid morphisms \(f:V\to(\Z/2\Z)^2\) given by \(a\mapsto(1,0)\) and \(b\mapsto (0,1)\), and \(g:V\to G=\Aut_\catt{Set}(\{1,2,3\})\) given by \(a\mapsto (1\,2)\) and \(b\mapsto (2\,3)\).
The ring \(\Z[\sqrt{2}]^2\) similarly has a \(V\)-grading \(\overline{T}\) under the proper definition of a grid-grading, which extends the definition of a group-grading, and the \((\Z/2\Z)^2\)-grading \(\overline{S}\) and \(G\)-grading \(\overline{R}\) of \(\Z[\sqrt{2}]^2\) can be obtained as \(f_*\overline{T}\) and \(g_*\overline{T}\).
\end{example}
\addtocounter{stelling}{1}

\begin{definition}
Let \(k\) be a commutative ring and \(R\) a \(k\)-algebra.
A {\em grid-grading} of \(R\) is a triple \(\overline{R}=(R,G,\mathcal{R})\) with \(G\) a grid and \(\mathcal{R}=\{R_g\}_{g\in G}\) a \(G\)-indexed pre-grading of \(R\), which satisfies \(1\in R_1\) and additionally for all \(g,h\in G\) such that \(R_g\cdot R_h \neq 0\) we have that \(g*h\) is defined and \(R_g\cdot R_h \subseteq R_{g*h}\).
We call \(\overline{R}\) an {\em abelian grid-grading} if \(G\) is an abelian grid.
For any morphism of grids \(f:G\to H\) we define \(f_*(\overline{R})=(R,H,f_*(\mathcal{R}))\), which is again a grading.
\end{definition}

The proof of \cite{Lenstra2018} generalizes effortlessly, which we present in Section~\ref{sec:UniGridGradThm}.

\begin{theorem}\label{UniGridGradThm}
Let \(R\) be a reduced order. 
Then \(R\) has a universal abelian group-grading {\em\cite{Lenstra2018}}, a universal group-grading and a universal grid-grading.
\end{theorem}

Using this theorem it is now not difficult to see that \(\overline{T}\) from Example~\ref{ex:why_grid} is a universal grid-grading.
For Theorem~\ref{UniGridGradThm} to work we require the somewhat artificial property of grids that \(1\) is the only idempotent.
If we drop this property, then the pre-grading \(\{R_s\}_{s\in \{0,1\}}\) of \(\Z^2\) with \(R_1=(1,1)\Z\) and \(R_0=(1,0)\Z\) can be turned into a `grid'-grading \(\overline{R}\). 
If a universal grading should still exist, then \(\overline{R}\) must be it since it cannot be refined, but reversing the coordinates of \(\Z^2\) gives a different `grid'-grading \(\overline{R}'\) of \(\Z^2\) between the two of which no mapping exists.

With some more effort, we also generalize Theorem~\ref{AlgGridThm}.

\begin{theorem}\label{AlgGridThm2}
There is a deterministic algorithm that takes a reduced order \(R\) as input and produces a universal grid grading \(\overline{U}=(R,\Upsilon,\{U_\upsilon\}_{\upsilon\in\Upsilon})\) of \(R\) in time \(n^{O(m)}\), where \(n\) is the length of the input and \(m=\#\minspec R\).
There is a deterministic algorithm that takes this \(\overline{U}\) as input and produces a finite presentation of a group \(\Gamma\) and a morphism of grids \(f:\Upsilon\to\Gamma\) such that \(f_*\overline{U}\) is a universal group grading of \(R\) in time \(n^{O(1)}\).
\end{theorem}

Note that the algorithm of Theorem~\ref{AlgGridThm2} does not compute \(f_*\overline{U}\).
The finitely presented group \(\Gamma\) will generally not be finite, as for \(\Z[\sqrt{2}]^2\) as in Example~\ref{ex:why_grid} it is \(\langle a,b \,|\, a^2=b^2=1\rangle\) by Proposition~\ref{prop:hierarchy}.
Hence the algorithm only returns a presentation.
More importantly, we are thus far unable to explicitly compute the decomposition of \(R\) corresponding to the grading \(f_*\overline{U}\) as it requires us to solve the word problem in \(\Gamma\), which for general groups is undecidable \cite{10.2307/1970103}.
We know of no algorithm that given \(f\), \(\overline{U}\) and \(x\in R\) decides whether \(x\) is homogeneous in \(f_*\overline{U}\), i.e. there is some \(\gamma\in\Gamma\) such that \(x\in\sum_{\upsilon\in f^{-1}\gamma} U_\upsilon\).

\clearpage
\section{Graded rings}

In this section we derive some properties of grids and gradings and introduce some tools to construct gradings from smaller gradings.
This is followed by some examples highlighting the difference between gradings graded with grids, groups and abelian groups.

\subsection{Basic properties}

In this section \(k\) will be a commutative ring and \(R\) will be a (not necessarily commutative) \(k\)-algebra.
For each grid-grading \(\overline{R}=(R,G,\{R_g\}_{g\in G})\) the \(k\)-module \(R_1\) is closed under multiplication and contains \(1\), and is thus a \(k\)-subalgebra of \(R\). Similarly \(R_g\) is an \(R_1\)-\(R_1\)-bimodule for all \(g\in G\).
We will use this fact without reference.

\begin{definition}\label{def:subgrid}
Let \(G=(S,1,D,*)\) be a grid. 
We call a \((T,1,E,\star)\) a {\em subgrid} of \(G\) if it is a grid such that \(T\subseteq S\), \(E=D\cap T^2\) and \(s \star t = s * t\) for all \((s,t)\in E\).
For \(X\subseteq S\) we write \(\langle X\rangle_\catt{Grd}\) or simply \(\langle X\rangle\) for the smallest subgrid of \(G\) containing \(X\).
\end{definition}

We similarly write \(\langle X \rangle_\catt{Grp}\) and \(\langle X\rangle_\catt{Ab}\) for the smallest subgroup respectively abelian subgroup containing \(X\). Note that if \(G\) is a group, the notation \(\langle X \rangle\) has become ambiguous. For \(G=\Z\), we have \(\langle 1 \rangle_\catt{Grd} = \Z_{\geq 0}\)  while \(\langle 1\rangle_\catt{Grp} =\Z\). 

\begin{definition}
Let \(\mathcal{C}\in\{\catt{Grd},\catt{Grp},\catt{Ab}\}\) and let \(\overline{R}=(R,G,\mathcal{R})\) be a \(\mathcal{C}\)-grading of \(R\). 
We say \(\overline{R}\) is {\em efficient} if \(G=\langle g\in G\,|\, R_g\neq 0\rangle_\mathcal{C}\) and say \(\overline{R}\) is {\em loose} if \(R_g \cdot R_h\neq 0\) for all \(g,h\in G\) such that \(\{g,h\}\neq \{1\}\) and \(g*h\) is defined.
\end{definition}

Note that whether a grading is efficient depends on the underlying category. By the same argument as before a group-grading can be efficient as group-grading but be not efficient as grid-grading. Note that a loose grading is necessarily efficient.

\begin{lemma}\label{lem:sub_initial}
Let \(\mathcal{C}\in\{\catt{Grd},\catt{Grp},\catt{Ab}\}\) and let \(\overline{R}=(R,G,\{R_g\}_{g\in G})\) be an efficient \(\mathcal{C}\)-grading.
If \(H\in\obj(\mathcal{C})\) and morphisms \(\varphi,\psi:G\to H\) satisfy \(\varphi_*\overline{R}=\psi_*\overline{R}\), then \(\varphi=\psi\).
\begin{proof}
Since \(\varphi_*\overline{R}=\psi_*\overline{R}\), we must have that \(\varphi\) and \(\psi\) agree on \(X=\{g\in G\,|\, R_g\neq 0\}\).
Let \(K=\{ x\in G \,|\, \varphi(x)=\psi(x) \}\supseteq X\). For \(x,y\in K\) such that \(x*y\) is defined we have \(\varphi(x*y)=\varphi(x)\varphi(y)=\psi(x)\psi(y)=\psi(x*y)\), so \(x*y\in K\) and \(K\) is a subgrid of \(G\). Additionally, if \(G\) is a group \(K\) is also closed under taking inverses, making \(K\) a subgroup of \(G\).
It follows that \(K\subseteq G=\langle X\rangle_\mathcal{C} \subseteq K\) since \(X\subseteq K\), so \(\varphi=\psi\).
\end{proof}
\end{lemma}

\begin{proposition}\label{prop:universal_implies_efficient}
Let \(\mathcal{C}\in\{\catt{Grd},\catt{Grp},\catt{Ab}\}\) and let \(\overline{R}=(R,G,\{R_g\}_{g\in G})\) be a \(\mathcal{C}\)-grading.
Let \(H=\langle g\in G\,|\, R_g\neq 0\rangle_\mathcal{C}\) and let \(\varphi:H\to G\) be the inclusion.
Then \(\overline{S}=(R,H,\{R_h\}_{h\in H})\) is an efficient \(\mathcal{C}\)-grading such that \(\varphi_*\overline{S}=\overline{R}\).
If for each \(\mathcal{C}\)-grading \(\overline{T}\) there exists a morphism \(\psi_*:\overline{R}\to\overline{T}\), then \(\overline{S}\) is a universal grading.
Moreover, universal \(\mathcal{C}\)-gradings are efficient.
\begin{proof}
It is an easy verification that \(\overline{S}\) is a grading and that \(\varphi_*\overline{S}=\overline{R}\).
If morphisms \(\psi_*:\overline{R}\to\overline{T}\) exist for all \(\mathcal{C}\)-gradings \(\overline{T}\), then the same holds for \(\overline{S}\) by composition with \(\varphi_*\). Then by Lemma~\ref{lem:sub_initial} such a map is unique, hence \(\overline{S}\) is universal.
If \(\overline{R}\) was already universal, then there exists a morphism \(\psi:G\to H\) such that \((\varphi\circ \psi)_*\overline{R}=\overline{R}\), hence \(\varphi\circ \psi=\id_G\) by uniqueness of morphisms. Similarly \(\psi\circ \varphi=\id_H\), hence \(\overline{R}=\overline{S}\) and \(\overline{R}\) is efficient.
\end{proof}
\end{proposition}

Proposition~\ref{prop:universal_implies_efficient} allows us to restrict ourselves to the category of efficient gradings when looking for a universal grading, since each grading is the image of an efficient grading under some morphism.

\begin{remark}\label{rem:reduced_ring_monoid}
Note that by the definition of a grid grading, if \(R_{\gamma_1}\dotsm R_{\gamma_m} \neq 0\) for \(\gamma_1,\dotsc,\gamma_m\in\Gamma\), then \(\gamma_1*\dotsm*\gamma_m\) is uniquely defined regardless of grouping of factors.
In particular, if \(x\in R_\gamma\) is not nilpotent then \(0\neq R_\gamma\dotsm R_\gamma\) and thus \(\gamma^n\) is uniquely defined for all \(n\geq 0\).
\end{remark}

We call an element \(x\) of a ring \(R\) {\em regular} if for all \(y\in R\) such that \(xy=0\) or \(yx=0\) we have \(y=0\).
Any \(x\in R\) that is not regular is called a {\em zero-divisor}.

\begin{proposition}\label{prop:basic_grad_props}
Let \(\overline{R}=(R,\Gamma,\{R_\gamma\}_{\gamma\in\Gamma})\) be a grid-grading of a reduced (not necessarily commutative) algebra \(R\) over a commutative ring \(k\) such that \(X=\{\gamma\in\Gamma\,|\,R_\gamma\neq0\}\) is finite. Then the following holds.
\begin{enumerate}[topsep=0pt,itemsep=-1ex,partopsep=1ex,parsep=1ex,label={\em(\arabic*)}]
\item For each \(\gamma\in X\) we have that \(\langle\gamma\rangle_\catt{Grd}\subseteq X\) is a finite cyclic group.
\item If \(R_1\) is a domain, then for each \(\gamma\in\Gamma\) all non-zero \(x\in R_\gamma\) are regular, and \(X=\langle X\rangle_\catt{Grd}\) is a group. 
If additionally \(R\) is commutative, then \(X\) is abelian.
\item If \(R_1\) is a field, then for each \(\gamma\in\Gamma\) all non-zero \(x\in R_\gamma\) are invertible with \(x^{-1}\in R_{\gamma^{-1}}\). 
Additionally, \(\dim_{R_1} R_\gamma=1\) for all \(\gamma\in X\).
\end{enumerate}
\begin{proof}
(1) Let \(\gamma\in\Gamma\) be such that \(R_\gamma\neq 0\).
Then \(\gamma^n\) is uniquely defined for all \(n\in\Z_{\geq 0}\) by Remark~\ref{rem:reduced_ring_monoid}, and since \(\{1,\gamma,\gamma^2,\dotsc\}\subseteq X\) is finite we must have \(\gamma^n=\gamma^m\) for some \(0\leq n < m\). 
Now \(\gamma^{2m(m-n)}=\gamma^{(2m-1)(m-n)}=\dotsm=\gamma^{m(m-n)}\) since \((2m-a)(m-n)\geq m\) for all \(a\leq m\).
Hence \(\gamma^{m(m-n)}\) is idempotent and \(\gamma^{m(m-n)}=1\).
It follows that \(\langle\gamma\rangle_\catt{Grd}\) is a finite group.

(2) Assume \(R_1\) is a domain.
Let \(\gamma_1,\dotsc,\gamma_n\in X\) and for each \(i\) let \(x_i\in R_{\gamma_i}\) be non-zero.
Applying induction to \(n\) we show that \(x_1\dotsm x_n\neq 0\).
Assume \(0\neq x:=x_1 \dotsm x_{n}\), then \(x \in R_\delta\) for some \(\delta\in X\). 
Using (1) we get \(0\neq x_{n+1}^{\ord(\gamma_{n+1})}, x^{\ord(\delta)}\in R_1\) and thus \(0\neq x^{\ord(\delta)} x_{n+1}^{\ord(\gamma_{n+1})}\) since \(R_1\) is a domain.
In particular, \(x_1\dotsm x_{n+1}=x\cdot x_{n+1}\neq 0\). 
It follows from Remark~\ref{rem:reduced_ring_monoid} that \(\gamma=\gamma_1*\dotsm*\gamma_n\) is uniquely defined regardless of grouping of factors and that \(0\neq R_\gamma\), so \(\gamma\in X\).
We conclude that \(X=\langle X\rangle_\catt{Grd}\) and that \(X\) is a monoid.
It then follows from (1) that \(X\) is a group.
If \(R\) is commutative, \( 0\neq R_\gamma R_\delta \cap R_\delta R_\gamma \subseteq R_{\gamma*\delta}\cap R_{\delta*\gamma}\) implies \(\gamma*\delta=\delta*\gamma\) for all \(\gamma,\delta\in X\), making \(X\) abelian.

Let \(x\in R_\delta\) be non-zero and \(y\in R\) such that \(xy=0\).
Write \(y=\sum_{\gamma} y_\gamma\) with \(y_\gamma\in R_\gamma\) for all \(\gamma\in\Gamma\).
Then \(xy=\sum_{\gamma} xy_\gamma = \sum_\gamma xy_{\delta^{-1}*\gamma}\) with \(xy_\gamma\in R_{\delta*\gamma}\).
By uniqueness of decomposition \(xy_\gamma=0\) and by the above \(y_\gamma=0\) for all \(\gamma\in\Gamma\).
It follows that \(y=0\), hence \(x\) is not a left zero-divisor.
Analogously \(x\) is not a right zero-divisor, so \(x\) is regular.

(3) Let \(\gamma\in X\) with \(n=\ord(\gamma)\) and let \(x\in R_\gamma\) be non-zero.
As \(x^n\in R_1\) is invertible, we have \(y=(x^{n})^{-1}\cdot x^{n-1}\in R_{\gamma^{n-1}}\) such that \(x y = yx = 1\).
Hence \(x\) is invertible with \(x^{-1}=y\in R_{\gamma^{-1}}\).
Consider \(R_\gamma\) as left \(R_1\)-module.
Then we have \(R_1\)-linear maps \(R_1\to R_\gamma\) and \(R_\gamma\to R_1\) given by \(z \mapsto zx\) and \(z\mapsto zx^{-1}\).
We note that they are mutually inverse, implying that \(R_1\cong_{R_1\catt{-Mod}} R_\gamma\), from which it follows that \(\dim_{R_1} R_\gamma = 1\).
\end{proof}
\end{proposition}

\subsection{Coproduct gradings and joint gradings}\label{sec:grid_prod}

A basic tool in our study of gradings is the ability to construct complex gradings from small building blocks, which we will treat in this section.
The following lemma is an exercise in elementary algebra left to the reader.

\begin{lemma}
Let \(\mathcal{G}=\{G_i\}_{i\in I}\) be a family of grids with \(G_i=(S_i,1_i,D_i,*_i)\). 
Then \(\mathcal{G}\) has a product \(\prod_{i\in I} G_i\), where the underlying set is \(\prod_{i\in I} S_i\) and the product \((a_i)_{i\in I} * (b_i)_{i\in I} := (a_i *_i b_i)_{i\in I}\) is defined precisely when it is for all components.
Additionally \(\mathcal{G}\) has a coproduct \(\coprod_{i\in I} G_i\) where the underlying set is \(S=\{1\}\amalg\coprod_{i\in I} (S_i\setminus\{1_i\})\) with \(1\) the unit.
For \(a,b\in S\setminus\{1\}\) we define \(a*b=a*_i b\) if and only if there is some \(i\in I\) such that \((a,b)\in D_i\). \qed
\end{lemma}

\begin{definition}\label{def:groupification}
Let \(G=(S,1,D,*)\) be a grid. We define the {\em groupification of \(G\)}, symbolically \(G^\catt{grp}\), to be the group with presentation
\[ \langle [g] \text{ for } g\in G \,|\, [1]=1,\, [g]\cdot[h]=[g*h] \text{ for } (g,h)\in D \rangle.\]
\end{definition}

\begin{lemma}\label{lem:groupification_adjunction}
Groupification defines a functor \(\us^{\catt{grp}}:\catt{Grd}\to\catt{Grp}\).
This functor is the left adjoint of the forgetful functor \(\catt{Grp}\to\catt{Grd}\).
\begin{proof}
The proof is straightforward and amounts to showing that for any grid \(G\) and group \(H\) any morphism of grids \(g:G\to H\) factors uniquely through \(G^\catt{grp}\).
Let \(F\) be the free group with set of symbols \(G\), let \(N\subseteq F\) be the normal subgroup such that \(F/N=G^\catt{grp}\) as in Definition~\ref{def:groupification} and let \(f:G\to F/N\) be the natural map.
We combine the universal property of the free group and the homomorphism theorem to obtain the unique \(h:F/N\to H\) as follows
\begin{center}
\begin{tikzcd}[column sep=5em]
G \arrow{rr}{f} \arrow{dr} \arrow[swap]{ddr}{g} &  & F/N \arrow[dashed]{ddl}{\exists! h} \\[-15pt]
& F \arrow{ur} \arrow[dashed]{d}{\exists!i}& \\ 
& H & 
\end{tikzcd}
\end{center}
where we use that always \(N\subseteq\ker(i)\).
\end{proof}
\end{lemma}

Groupification \(\us^\catt{grp}\) is to grids what abelianization \(\us^\catt{ab}\) is to groups.
Since the forgetful functors \(\catt{Ab}\to\catt{Grp}\) and \(\catt{Grp}\to\catt{Grd}\) are right adjoint they commute with limits.
In particular, the product of groups taken in the category of groups is canonically isomorphic to the same product taken in the category of grids instead.
The same holds for products of abelian groups in the category of groups. 
Therefore, when we take products of groups or abelian groups we do not have to specify the underlying category.
However, the coproduct for each category is, in general, fundamentally different, so here we need to be careful.
The most notable difference we encounter is that a coproduct of finitely many finite grids is finite, just like in the category of abelian groups, while a coproduct in the category of groups will only be finite if all but possibly one of the groups is trivial.

\begin{definition}\label{def:coproduct_grading}
Let \(\mathcal{R}=\{\overline{R^i}\}_{i\in I}\) be a finite family of grid-gradings with \(\overline{R^i}=(R_i,G_i,\{R_{i,g}\}_{g\in G_i})\).
Then \(R=\prod_{i\in I} R_i\) has a grading \(\coprod_{i\in I} \overline{R^i}=(R,G,\{R_g\})\) called the {\em coproduct of \(\mathcal{R}\)}, with \(G=\coprod_{i\in I} G_i\) taken in the category of grids and \(R_g=R_{i,g}\) if \(g\in G_i\setminus \{1\}\) and \(R_1=\prod_{i\in I} R_{i,1}\).
\end{definition}

One can verify that the coproduct of efficient gradings is again an efficient grading.
Another construction is the joint grading, which uses the product instead of the coproduct.

\begin{lemma}\label{lem:intersect_grading}
Let \(k\) be a commutative ring and \(R\) be a \(k\)-algebra. 
Let \(I\) be a set and let \(\{\overline{R}_i\}_{i\in I}\) with \(\overline{R}_i = (R,\Gamma_i,\{R_{i,\gamma}\}_{\gamma\in\Gamma_i})\) be a collection of grid-gradings.
Let \(S_{\delta}=\bigcap_{i\in I} R_{i,\delta(i)}\) for all \(\delta\in\Delta=\prod_{i\in I}\Gamma_i\).
If the natural map \(\bigoplus_{\delta\in\Delta} S_{\delta} \to R\) is surjective, then \(\bigcap_{i\in I}\overline{R}_i := (R,\Delta,\{S_{\delta}\}_{\delta\in\Delta})\) is a grid-grading.
\begin{proof}
Clearly the natural map \(\bigoplus_{\delta\in\Delta} S_{\delta} \to R\) is injective and thus an isomorphism.
It suffices to note that \(S_\delta S_{\delta'} \subseteq R_{\delta(i)} R_{\delta'(i)} \subseteq R_{(\delta\delta')(i)}\) for all \(\delta,\delta'\in\Delta\) and \(i\in I\) hence \(S_{\delta} S_{\delta'} \subseteq \bigcap_{i\in I} R_{(\delta\delta')(i)} = S_{\delta\delta'}\).
As \(1\in S_{1}\) we have that \(\bigcap_{i\in I}\overline{R}_i\) is a grading.
\end{proof}
\end{lemma}

\begin{definition}
Using the same notation as in Lemma~\ref{lem:intersect_grading}, we call \(\bigcap_{i\in I}\overline{R}_i\) the {\em joint grading} of \(\{\overline{R}_i\}_{i\in I}\) when it is a grading.
\end{definition}

\begin{remark}\label{rmk:joint_uni_grad_exist}
If \(R\) has a universal grading \(\overline{U}=(R,\Upsilon,\{U_{\upsilon}\}_{\upsilon\in \Upsilon})\), then \(\overline{R}\cap\overline{S}\), with \(\overline{R}\) and \(\overline{S}\) gradings of \(R\), is always a grading. Namely the (unique) maps \(f:\Upsilon\to \Gamma\) and \(g:\Upsilon\to\Delta\) such that \(f_* \overline{U} = \overline{R}\) and \(g_* \overline{U}=\overline{S}\) give a map \(h:\Upsilon\to\Gamma\times\Delta\) and a \(\Gamma\times\Delta\)-grading \(\overline{T}=h_* \overline{U}\). The projection maps \(\pi_\Gamma : \Gamma\times\Delta\to\Gamma\) and \(\pi_\Delta:\Gamma\times\Delta\to\Delta\) then satisfy \(\pi_{\Gamma*} \overline{T} = \overline{R}\) and \(\pi_{\Delta*}\overline{T} = \overline{S}\).
One then easily verifies \(\overline{T}=\overline{R}\cap\overline{S}\).
This argument trivially extends to joint gradings of infinitely many gradings.
To show that \(R\) does not have a universal grading, it suffices to give gradings \(\overline{R}\) and \(\overline{S}\) such that \(\overline{R}\cap\overline{S}\) is not a grading. 
As Proposition~\ref{thm:commutative_implies_universal} will show, under mild finiteness conditions the existence of such \(\overline{R}\) and \(\overline{S}\) is also necessary. 
If we drop these finiteness conditions, Example~\ref{ex:non-compact_counter} provides a counter-example.
\end{remark}

\subsection{Examples}

In this section we consider some examples of gradings. 

\begin{example}\label{ex:group_ring}
Let \(\Gamma\) be a finite group and let \(k\) be a commutative ring which is connected, meaning it is not isomorphic to the product of two non-zero rings.
Consider \(k^\Gamma\), which is a \(k\)-module of which we denote the elements as \(\sum_{\gamma\in\Gamma} x_\gamma \gamma\) for some \(x_\gamma\in k\).
Then \(k^\Gamma\) can be equipped with a \(k\)-algebra structure by defining the multiplication as the \(k\)-linear continuation of the multiplication of \(\Gamma\).
We write \(k[\Gamma]\) for this ring, which we call the {\em group ring of \(\Gamma\)}.
The group ring comes with a natural \(\Gamma\)-grading \(\overline{U}=(k[\Gamma],\Gamma,\{k\cdot \gamma\}_{\gamma\in\Gamma})\).
Note that \(\overline{U}\) cannot be further refined by connectedness of \(k\).
As \(\overline{U}\) is also efficient it must be the universal grading of \(k[\Gamma]\) whenever a universal grading exists.
If \(k=\Z\) and \(\Gamma\) is abelian, then \(k[\Gamma]\) is a reduced order and \(\overline{U}\) must be the universal grading of \(k[\Gamma]\) since it exists by Theorem~\ref{UniGridGradThm}.

If \(\Gamma\) and \(\Delta\) are finite abelian groups, then \(k[\Gamma\times\Delta]\cong k[\Gamma]\tensor_k k[\Delta]\).
Then by the structure theorem for finite abelian groups we may study the structure of \(k[\Gamma]\) through that of \(k[C_{p^n}]\), where \(p\) is a prime, \(n\geq 1\) and \(C_r\) is the cyclic group of order \(r\).
Additionally, \(k[\Gamma]\cong k\tensor_\Z \Z[\Gamma]\).
We have
\[\Q[C_{p^n}]\cong \Q[X]/(X^{p^n}-1)\cong \Q[X]/(X^{p^{n-1}}-1) \times \Q[X]/\bigg( \frac{X^{p^n}-1}{X^{p^{n-1}}-1}\bigg) \cong \Q[C_{p^{n-1}}] \times \Q(\zeta_{p^n}).\]
As \(\Q(\zeta_a)\tensor_\Q\Q(\zeta_b)\cong\Q(\zeta_{ab})\) for \(a,b\in\Z_{>0}\) coprime, we have that for \(k=\Q\) the group rings are isomorphic to products of cyclotomic fields.
\end{example}

By Example~\ref{ex:group_ring} we may consider products of cyclotomic rings as the first generalization of group rings. 
It gives us motivation to study the gradings of cyclotomic fields, which we will do more thoroughly in Section~\ref{sec:cyclo_ring}.

\begin{example}\label{intersection_counter}
Consider the ring \(K=\Q(\zeta_8)\), where \(\zeta_8\) is a primitive \(8\)-th root of unity.
Since it is a field, all its efficient gradings are graded with abelian groups of size dividing \([\Q(\zeta_8):\Q]=4\) by Proposition~\ref{prop:basic_grad_props}.
It has an obvious grading \(\overline{A}=(K,\Z/4\Z,\{A_{n}\}_{n\in\Z/4\Z})\) with \(A_n=\zeta_8^n \cdot \Q\), which is well defined since \(\zeta_8^4\in \Q\).
The natural projection \(\pi:\Z/4\Z\to\Z/2\Z\) then induces the grading \(\overline{B}=\pi_*\overline{A} = (K,\Z/2\Z,\{B_n\}_{n\in\Z/2\Z})\) with \(B_n=\zeta_8^n \cdot \Q(\zeta_4)\).
Two other \(\Z/2\Z\)-gradings of \(K\) are \(\overline{C}\) and \(\overline{D}\) with \(C_n = \sqrt{2^{n}}\cdot \Q(\zeta_4\cdot \sqrt{2})=\zeta_4^n\cdot \Q(\zeta_4\cdot\sqrt{2})\) and \(D_n=\zeta_4^n \cdot \Q(\sqrt{2})\) as \(\zeta_8=(1+\zeta_4)/\sqrt{2}\).
As we will show in Section~\ref{sec:cyclo_ring}, we have in fact listed all efficient gradings of \(K\) with a cyclic group, up to group isomorphism.
Of these we can construct a joint grading \(\overline{E}=\overline{C}\cap\overline{D}=(K,(\Z/2\Z)^2,\{E_{m,n}\}_{(m,n)\in(\Z/2\Z)^2})\) with \(E_{m,n}=\sqrt{2^{m+n}}\cdot \zeta_4^n \cdot \Q\).
However, \(\overline{A}\) and \(\overline{E}\) have no joint grading, as \(\sum_{k,m,n} A_k\cap E_{m,n}=\Q(\zeta_4)\neq K\). 
Thus \(K\) has no universal grading by Remark~\ref{rmk:joint_uni_grad_exist}. 

Now consider \(R=\Z[\zeta_8]\), which is a reduced order.
It certainly has \(\overline{U} := (R,\Z/4\Z,\{A_n\cap R\}_{n\in\Z/4\Z})\) as grading.
On the other hand, \(\zeta_8=(\frac{1}{2}\sqrt{2})+(\frac{1}{2}\zeta_4\sqrt{2})\not\in\bigoplus_{(m,n)\in(\Z/2\Z)^2} (E_{m,n}\cap R)\) as \(\sqrt{2}/2\not\in R\), so \(\overline{E}\) does not similarly produce a grading of \(R\).
Hence the non-existence of a joint grading of \(\overline{A}\) and \(\overline{E}\) does not contradict that \(R\) should have a universal grading.
Moreover, since \(\overline{U}\) cannot be further refined, it must be the universal grading provided by Theorem~\ref{UniGridGradThm}.
\end{example}

The following two examples highlight the difference between grid-gradings, group-gradings and abelian group-gradings.
Specifically, in Example~\ref{ex:grid_not_group} we construct a ring \(A\) and a grid-grading \(\overline{A}\) of \(A\) such that there exists no group-grading \(\overline{B}\) of \(A\) with the same non-zero homogeneous components as \(\overline{A}\).
We do the same in Example~\ref{ex:commutative_non_abelian_group} for groups and abelian groups.

\phantom{X}

\noindent\begin{minipage}[c]{0.7\linewidth}
\begin{example}\label{ex:grid_not_group}
Grids need not be associative.
Consider the following figure consisting of seven points \(S=\{a,\dotsc,g\}\) and \(4\) lines each containing exactly three points.
We define multiplication on \(G=\{1\}\cup S\) such that \(\{1\}\cup L\) is a Klein four-group for each line \(L\), and all other products are undefined.
Then \(G\) is a grid, and \((bd)e=fe=g\) while \(b(de)=bc=a\). Hence \(G\) is not associative.
\end{example}
\end{minipage}
\begin{minipage}[c]{0.3\linewidth}
\begin{center}
\begin{tikzpicture}[main_node/.style={circle,fill=black,draw,minimum size=.3em,inner sep=.3pt]}]
\draw (0,0) -- (2,0) -- (1,1) -- (3,1);
\draw (1,0) -- (2,1);

\coordinate (a) at (0,0);
\coordinate (b) at (1,0);
\coordinate (c) at (2,0);
\coordinate (d) at (1.5,.5);
\coordinate (e) at (1,1);
\coordinate (f) at (2,1);
\coordinate (g) at (3,1);
\node[main_node,label={\(a\)}] (A) at (a) {};
\node[main_node,label={\(b\)}] (B) at (b) {};
\node[main_node,label={\(c\)}] (C) at (c) {};
\node[main_node,label={\(d\)}] (D) at (d) {};
\node[main_node,label={\(e\)}] (E) at (e) {};
\node[main_node,label={\(f\)}] (F) at (f) {};
\node[main_node,label={\(g\)}] (G) at (g) {};

\end{tikzpicture}
\end{center}\vspace{1em}
\end{minipage}

Now consider the efficient grading \(\overline{R}=(\Q(\sqrt{2},\sqrt{3}),(\Z/2\Z)^2,\{ \sqrt{2^u3^v}\cdot \Q \}_{(u,v)\in(\Z/2\Z)^2})\).
For each line \(L\) we have an injective morphism of grids \(\varphi_L:(\Z/2\Z)^2\cong(\{1\}\cup L)\to G\) and the coproduct of these morphisms gives a surjection \(\varphi:\coprod_{L} (\Z/2\Z)^2\to G\). 
Now \(\overline{A}=\varphi_*\coprod_{L} \overline{R}\) gives a \(G\)-grading of \(A=\Q(\sqrt{2},\sqrt{3})^4\) without any zero homogeneous components.
Then for any grading \(\overline{B}=(A,H,\{B_h\}_{h\in H})\) of \(A\) with the same non-zero homogeneous components as \(\overline{A}\) there exists a unique injective morphism \(\psi:G\to H\) such that \(\psi_*\overline{A}=\overline{B}\).
If \(H\) is a group, then \(H\) is certainly associative, so \(\psi(g)=\psi(f)\psi(e)=\psi(b)\psi(d)\psi(e)=\psi(b)\psi(c)=\psi(a)\), a contradiction. Thus the decomposition of \(\overline{A}\) can be turned into a grading with a grid but not with a group.

\begin{example}\label{ex:commutative_non_abelian_group}
A commutative ring can be efficiently graded with a finite non-abelian group, as seen in Example~\ref{ex:why_grid}.
However, for this ring there exists a grading with an abelian group which has the same non-zero homogeneous components.
Here we will construct a grading with a non-abelian group for which this is not the case.

For \(n\in\Z_{\geq 1}\) consider the groups \(N=\{\pm 1\}^{n+1}=\langle z, x_1, \dotsc, x_n\rangle\) and \(G=\{\pm1\}^n=\langle x_{n+1},\dotsc,x_{2n}\rangle\).
On the generators of \(G\) define an action \(\varphi:G\to\Aut(N)\) as
\[ x_{n+i} \mapsto \big( x_j \mapsto x_j z^{\delta_{ij}},\ z \mapsto z \big) \quad \text{for all } 1\leq i,j \leq n, \]
where \(\delta_{ij}=1\) if \(i=j\) and \(\delta_{ij}=0\) otherwise, and define \(\text{ES}_n=N \rtimes_\varphi G\).
Although we do not use this fact, we note that \(\text{ES}_n\) is an {\em extra special \(2\)-group} \cite{gorenstein1986}, as its center is \(\langle z\rangle\cong \F_2\), and \(\text{ES}_n / \langle z \rangle\cong \F_{2}^{2n}\).
Consider the \(\text{ES}_2\)-valued matrix 
\begin{align*}
M = 
\begin{blockarray}{cccc}
1 & 2 & 3 \\
\begin{block}{(ccc)c}
  x_1 & x_2 & x_1x_2 & 4 \\
  x_4 & x_3 & x_3x_4 & 5 \\
  x_1x_4 & x_2x_3 & x_1x_2x_3x_4 & 6 \\
\end{block}
\end{blockarray}
\end{align*}
together with a labeling of its rows and columns. 
Consider the \((\Z/2\Z)^3\)-grading \(\overline{R}\) of \(R=\Q(\sqrt{2},\sqrt{3},\sqrt{5})\) with \(R_{(a,b,c)}=\sqrt{2^a3^b5^c}\cdot \Q\).
For each \(i\) let \(\Delta_i\) be the group generated by \(z\) and the entries of row/column \(i\) of \(M\) and note that we have an isomorphism \(f_i:(\Z/2\Z)^3\to \Delta_i\).
Consider the \(\Delta_i\)-grading \(\overline{R^i}=f_{i*}\overline{R}\) of \(R\).
Then we grade \(A=R^6\) with \(\Delta=\coprod_{i=1}^6 \Delta_i\) as the coproduct (in the category of grids) of \(\overline{R^1},\dotsc,\overline{R^6}\), and let \(\overline{A}\) be the image of this coproduct under the natural map \(g:\Delta\to \text{ES}_2\). 
Since the \(\Delta_i\) generate \(\text{ES}_2\) the grading \(\overline{A}\) is efficient.
Thus the commutative ring \(A\) can be efficiently graded with a non-commutative group.

Assume there is a \(\Gamma\)-grading \(\overline{B}\) of \(A\) with the same non-zero components as \(\overline{A}\) and with \(\Gamma\) an abelian group.
Let \(D=\bigcup_i \Delta_i\) and note that \(A_\delta\neq 0\) for all \(\delta\in D\).
Thus there exists a (unique) injective map \(\psi:D\to\Gamma\) such that \(A_\delta=B_{\psi(\delta)}\) for all \(\delta\in D\).
Moreover, the restriction \(\psi|_{\Delta_i}\) is a morphism of groups for each \(i\), as \(\{0\} \subsetneq A_\delta A_{\delta'} \cap A_{\delta\delta'} = B_{\psi(\delta)} B_{\psi(\delta')} \cap B_{\psi(\delta\delta')}\) implies \(\psi(\delta)\psi(\delta')=\psi(\delta\delta')\) for all \(\delta,\delta'\in\Delta_i\). 
Let \(p_i\) be the product of the three entries of row/column \(i\) of \(M\) for each \(i\), and note \(p_i=1\) for \(i\in\{1,\dotsc,5\}\), while \(p_6=z\).
As \(\Gamma\) is abelian, we get
\[ 1 = \prod_{i=1}^3 \psi(  p_i ) = \prod_{a=1}^3 \prod_{b=1}^3 \psi( M_{ab}) = \prod_{i=4}^6 \psi( p_i ) = \psi(z),  \]
so \(\psi\) is not injective, a contradiction.
Hence there is no abelian group-grading with the same non-zero components as \(\overline{A}\).
\end{example}

\clearpage
\section{Gradings as ring automorphisms}

\subsection{Cyclotomic number fields and Steinitz numbers}

In this section we will define cyclotomic rings and derive some properties.

\begin{definition}
The {\em Steinitz numbers} are the formal products \(\prod_p p^{n_p}\) with \(p\) ranging over the primes and \(n_p\in\Z_{\geq 0}\cup\{\infty\}\) for all \(p\). We denote the set of Steinitz numbers by \(\mathbb{S}\). 
We define multiplication on Steinitz numbers as \( (\prod_p p^{n_p}) \cdot (\prod_p p^{m_p}) := \prod_p p^{n_p+m_p}\), with the convention that \(x+\infty=\infty+x=\infty\) for all \(x\in\Z_{\geq0}\cup\{\infty\}\).
For \(x=\prod_p p^{n_p}\in\mathbb{S}\) we define {\em the order of \(x\) at the prime \(p\)} as \(\ord_p(x):=n_p\).
We say \(x\in \mathbb{S}\) {\em divides} \(y\in\mathbb{S}\), symbolically \(x\divs y\), if there exists some \(z\in\mathbb{S}\) such that \(xz=y\), or equivalently if \(\ord_p(x)\leq\ord_p(y)\) for all primes \(p\).
We also write \(\infty\) for the Steinitz number \(\prod_{p}p^\infty\).
\end{definition}

\begin{remark}
Note that \(\Z_{>0}\) is naturally embedded in \(\mathbb{S}\) as the set of all Steinitz numbers \(\prod_{p}p^{n_p}\) for which only finitely many \(n_p\) are non-zero and none are infinity. 
For this subset of \(\mathbb{S}\) the defined concepts of multiplication, divisibility and order agree with those defined for \(\Z_{>0}\). 
Hence we will interpret this embedding as a true inclusion.
Each definition in this thesis involving Steinitz numbers will generalize similar definitions on \(\Z_{>0}\) if such exist.
We will leave it to the reader to verify this.
\end{remark}

\begin{definition}
Consider the natural projections \(\pi_{mn}:\Z/n\Z\to\Z/m\Z\) for \(m\divs n\) integers.
Note that \(\Z_{>0}\) is a partially ordered set with respect to divisibility, such that \(\{\pi_{mn}\}_{m\divs n \in\Z_{>0}}\) becomes an inverse system of homomorphisms.
Then we define the {\em ring of profinite integers} as the inverse limit \(\hat{\Z}=\lim\limits_{\substack{\longleftarrow\\ n\in\Z_{>0}}} \Z/n\Z\) of this system.
\end{definition}

\begin{remark}
We may interpret Steinitz numbers as special ideals of \(\hat{\Z}\).
Taking \(e=\prod_p p^{n_p}\in\mathbb{S}\), we define through abuse of notation its ideal \(e\hat{\Z}=\bigcap_{d\in\Z_{>0},\ d\divs e} d\hat{\Z}\).
Note that for \(e\in\Z_{>0}\) we have that \(e\hat{\Z}\) is unambiguous, as the ideal is the same when \(e\) is interpreted as a Steinitz number.
We have a Chinese remainder theorem \(\hat{\Z}/(\prod_p p^{n_p})\hat{\Z}\cong_\catt{Rng}\prod_p \hat\Z/p^{n_p}\hat\Z\) and also that for \(n\in\Z_{>0}\) the projection \(\hat\Z\to\Z/n\Z\) induces an isomorphism \(\hat{\Z}/n\hat\Z\cong_\catt{Rng} \Z/n\Z\).
\end{remark}

\begin{definition}
For \(N\subseteq\mathbb{S}\) we define the {\em greatest common divisor} and {\em least common multiple} of the elements of \(N\) as the Steinitz numbers 
\[\gcd(N):=\prod_{p}p^{ \inf\{ \ord_p(x)\,|\, x\in N \} } \quad\text{and}\quad \lcm(N):=\prod_{p}p^{\sup\{ \ord_p(x)\,|\, x\in N \}} \]
respectively. In particular, we have a definition of \(\gcd\) and \(\lcm\) on infinite subsets of \(\Z_{>0}\).
\end{definition}

\begin{definition}
Let \(G\) be a group for which each element has finite order, i.e. \(G\) is {\em torsion}.
We define the {\em exponent} of \(G\) as the Steinitz number \(e(G):=\lcm\{ \ord(g) \,|\, g\in G \}\).
Additionally, when we say \(G\) has exponent \(e\in\mathbb{S}\), we implicitly assume \(G\) is torsion.
\end{definition}

\begin{remark} \label{rem:exp_action}
The profinite integers form a ring, so in particular a multiplicative monoid.
If \(G\) is torsion we have a right action \(G\times\hat{\Z}\to G\) of \(\hat{\Z}\) as monoid given by \((g,\underline{n})\mapsto g^n\), where \(n\) is any integer such that \(n\equiv \underline{n} \mod \ord(g)\).
It satisfies \(1^n=1=g^0\), \(g^1=g\), \(g^{n+m}=g^ng^m\) and \(g^{nm}=(g^n)^m\) for all \(g\in G\) and \(n,m\in\hat{\Z}\).
When \(G\) is abelian \(\hat{\Z}\) even respects the group structure as \((gh)^n=g^nh^n\), making \(G\) a \(\hat{\Z}\)-module.
The map \(\hat{\Z}\to\End_\texttt{Set}(G)\) has kernel \(e(G)\hat{\Z}\), so even \(\hat{\Z}/e(G)\hat{\Z}\) acts on \(G\).
\end{remark}

\begin{definition}
A {\em topological group} is a group \(G\) equipped with a topology such that the multiplication \(G\times G\to G\) and inversion \(G\to G\) are continuous maps, where \(G\times G\) is equipped with the product topology.
We call a topological group {\em locally compact} if the underlying topology is Hausdorff and locally compact.
We write \(\catt{LCA}\) for the category of locally compact abelian groups with as morphisms the continuous group homomorphisms. 
We have a natural inclusion \(\catt{Ab}\to\catt{LCA}\) where we equip each group with the discrete topology.
For topological spaces \(X\) and \(Y\) the {\em compact-open topology} on \(Y^X\) is the topology with subbase \(\{V(K,U)\,|\,K\subseteq X\text{ compact},\, U\subseteq Y\text{ open}\}\), where \(V(K,U)=\{ f\in Y^X \,|\, f(K)\subseteq U\}\).
For any \(\Gamma,\Delta\in\obj(\catt{LCA})\) we equip \(\Hom_\catt{LCA}(\Gamma,\Delta)\) with a group structure given by point-wise multiplication and with the subspace topology from \(\Delta^\Gamma\) with the compact-open topology, making it a topological abelian group.
\end{definition}

\begin{definition}\label{def:dual}
Let \(\mathbb{T}=\{ z\in \C \,|\, |z|=1\}\), which is a compact abelian group.
For \(e\in\mathbb{S}\) we write \(\mu_e\subseteq\mathbb{T}\) for its \(e\)-torsion subgroup, which we equip with the discrete topology.
For a locally compact abelian group \(\Gamma\) we define \(\widehat{\Gamma}=\Hom_\catt{LCA}(\Gamma,\mathbb{T})\) as the {\em dual group} of \(\Gamma\), which is a locally compact abelian group as well (Theorem 1.2.6 of \cite{rudin2017fourier}).
For \(f\in\Hom_\catt{LCA}(\Gamma,\Delta)\) we write \(\widehat{f}\in\Hom_\catt{LCA}(\widehat{\Delta},\widehat{\Gamma})\) for the map given by \(\chi\mapsto \chi\circ f\).
Then \(\hatf:\catt{LCA}\to\catt{LCA}\) is a contravariant functor.
\end{definition}

\begin{remark}
For \(e\in\mathbb{S}\) and \(\chi\in\widehat{\mu}_e\) we have \(\chi(\mu_e)\subseteq\mu_e\).
Hence \(\widehat{\mu}_e\) is closed under composition and this composition is continuous.
This turns \(\widehat{\mu}_e\) into a topological ring where the addition is pointwise multiplication and multiplication is composition.
\end{remark}

\begin{lemma}\label{lem:profinite_dual_isomorphism}
For all \(e\in\mathbb{S}\) the natural action \(\varphi_e:\hat{\Z}/e\hat{\Z}\to\End_\catt{Ab}(\mu_e)\cong\widehat{\mu}_e\) as in Remark~\ref{rem:exp_action} is an isomorphism of rings.
\begin{proof}
First assume \(e\in\Z_{>0}\).
Then \(\mu_e\) is the zero set of \(X^e-1\) in \(\C[X]\), a polynomial which has precisely \(e\) distinct roots, so \(\#\mu_e=e\).
If \(d\divs e\) is the exponent of \(\mu_e\), then \(X^d-1\) has at least \(e\) distinct roots and thus \(d\geq e\) because \(\C\) is a domain. 
Hence \(d=e\) and thus \(\mu_e\) is cyclic.
If we fix any generator \(\zeta\in\mu_e\), then \(\chi\in\widehat{\mu}_e\) is uniquely defined by \(\chi(\zeta)\).
As \(\varphi_e(x)(\zeta)=\zeta^x\) for all \(x\in\hat{\Z}/e\hat\Z\) we have that \(\varphi_e\) is surjective, and since \(\zeta^x=\zeta^y\) implies \(x=y\) for \(x,y\in\hat{\Z}/e\hat{\Z}\) it is also injective.
For general \(e\in\mathbb{S}\) we note that with \(d\divs e\) an integer the diagram
\[\begin{tikzcd}
\hat{\Z}/e\hat{\Z} \arrow{r}\arrow{d}
& \widehat{\mu}_e \arrow{d} \\
\Z/d\Z \arrow{r}{\sim}
& \widehat{\mu}_d
\end{tikzcd}\]
is commutative, hence
\[ \hat{\Z}/e\hat{\Z}=\lim_{\substack{\longleftarrow\\ d\divs e}} \Z/d\Z \cong \lim_{\substack{\longleftarrow\\ d\divs e}} \Hom_\catt{LCA}(\mu_d,\mathbb{T}) \cong \Hom_\catt{LCA}\big( \lim_{\substack{\longrightarrow\\ d\divs e}} \mu_d, \mathbb{T}\big) = \Hom_\catt{LCA}(\mu_e,\mathbb{T}) = \widehat{\mu}_e\]
\end{proof}
\end{lemma}

\begin{corollary}\label{cor:cyclotomic_automorphisms}
Let \(e\in\mathbb{S}\). Then \((\hat{\Z}/e\hat{\Z})^*\cong \Aut_\catt{Ab}(\mu_e) \cong \Aut_\catt{Rng}(\Q(\mu_e)) \cong \Aut_\catt{Rng}(\Z[\mu_e])\).
\begin{proof}
The first isomorphism is a direct consequence of Lemma~\ref{lem:profinite_dual_isomorphism}.
For \(k\in\{\Z,\Q\}\) and \(e\) an integer the isomorphism \(\Aut_\catt{Grp}(\mu_e)\cong\Aut_\catt{Rng}(k[\mu_e])\) is proved in Theorem 3.1 in \cite{algebra}, and the general case follows from \(\mu_e=\lim\limits_{\substack{\longrightarrow \\ d\divs e}} \mu_d\) and \(k[\mu_e]=\lim\limits_{\substack{\longrightarrow \\ d\divs e}} k[\mu_d]\).
\end{proof}
\end{corollary}

\begin{lemma}\label{lem:prime_roots} \label{lem:root_diff_division}
Let \(e\in\Z_{>0}\), \(\zeta\in\mu_e\) primitive and \(\pi=1-\zeta\in\Z[\mu_e]\).
\begin{enumerate}[topsep=0pt,itemsep=-1ex,partopsep=1ex,parsep=1ex,label={\em(\arabic*)}]
\item For all \(a\in \Z\) there exist \(x\in\Z[\mu_e]\) such that \(x \pi = 1-\zeta^a\).
If \(\gcd(a,e)=1\), then \(x\) is a unit.
\item If \(\eta,\theta\in\mu_e\) are distinct, then \(e=x(\eta-\theta)\) for some \(x\in\Z[\mu_e]\). 
\item If \(e=p^n\) for some \(n\in\Z_{>0}\) and prime \(p\), then we may write \(x\pi^m=p\) with \(m=p^{n}-p^{n-1}\) and \(x\in\Z[\mu_{p^n}]^*\). 
\end{enumerate}
\begin{proof}
Recall the following basic algebraic identities in \(\Z[\mu_e][X]\) with \(b,c\in\Z_{>0}\).
\[ \sum_{i=0}^{b-1} X^{ci} = \frac{X^{bc}-1}{X^c-1} = \prod_{\xi\in\mu_{bc}\setminus\mu_c} ( X - \xi ).  \]

(1) Fix \(c=1\). Taking \(b=a\) and \(X=\zeta\) in the above identities we obtain \(x\pi=1-\zeta^a\) with \(x=\sum_{i=0}^{a-1}\zeta^i\).
If \(\gcd(a,e)=1\), then there exists a \(b\in\Z_{>0}\) such that \(ab\equiv 1 \mod e\). 
From the first identity with this \(b\) and \(X=\zeta^a\) it follows that \(x\) is invertible.

(2) Let \(x=\eta^{-1} \cdot \prod_{\xi\in\mu_e\setminus\{1,\eta^{-1}\theta\}} (1-\xi)\).
From the above identities with \(b=e\), \(c=1\) and \(X=1\) we obtain \(x(\eta-\theta)=\prod_{\xi\in\mu_{e}\setminus\mu_1}(X-\xi)=\sum_{i=0}^{e-1}X=e\).

(3) Each \(\xi\in\mu_e\) may be written \(\xi=\zeta^a\) for some \(a\in\Z\), and \(\xi\not\in\mu_{p^{n-1}}\) if and only if \(\gcd(a,e)=1\).
Then by (1) for all \(\xi\in\mu_{p^n}\setminus\mu_{p^{n-1}}\) there exist \(x_\xi\in\Z[\mu_{p^n}]^*\) such that \(1-\xi=\pi x_\xi\).
Taking \(b=p\), \(c=p^{n-1}\) and \(X=1\) it follows that \(x=\prod_{\xi\in\mu_{p^n}\setminus\mu_{p^{n-1}}}x_\xi\) satisfies \(x\pi^m=p\).
\end{proof}
\end{lemma}

\begin{lemma}\label{trace_properties}
For \(n=\prod_{p} p^{n_p} \in \Z_{>0}\) write \(\sqrfree(n)=\prod_{p} p^{\min\{n_p,1\}}\) and \(\omega(n)=\#\{p \,|\, n_p\geq 1\}\).
Then \(\Tr_{\Q(\zeta_n)/\Q}(\zeta_{\sqrfree(n)})=(-1)^{\omega(n)}\cdot n/\sqrfree(n)\) for any primitive \(\sqrfree(n)\)-th root of unity \(\zeta_{\sqrfree(n)}\).
\begin{proof}
Consider the map \(\psi\) given by \(n\mapsto \Tr_{\Q(\zeta_n)/\Q}(\zeta_{\sqrfree(n)})\).
Let \(n,m\in\Z_{>0}\) be coprime. 
We have \((\Z/nm\Z)^*\cong (\Z/n\Z)^*\times(\Z/m\Z)^*\) by the Chinese Remainder Theorem.
Writing \(\zeta_{\sqrfree(nm)}=\zeta_{\sqrfree(n)}\zeta_{\sqrfree(m)}\) and using Corollary~\ref{cor:cyclotomic_automorphisms} we get
\[ \psi(nm) = \sum_{a\in(\Z/nm\Z)^*} \zeta_{\sqrfree(nm)}^a = \sum_{a\in(\Z/nm\Z)^*} \zeta_{\sqrfree(n)}^a\zeta_{\sqrfree(m)}^a= \sum_{\substack{a\in(\Z/n\Z)^*\\b\in(\Z/m\Z)^*}} \zeta_{\sqrfree(n)}^a\zeta_{\sqrfree(m)}^b = \psi(n)\psi(m).  \]
Thus \(\psi\) is a multiplicative function, and so is \(n\mapsto (-1)^{\omega(n)} \cdot n / \sqrfree(n)\).
Because of this we may assume without loss of generality that \(n=p^k\). Here 
\[\psi(n)= \sum_{a\in(\Z/p^k\Z)^*} \zeta_p^a = p^{k-1} \sum_{a\in(\Z/p\Z)^*} \zeta_p = -p^{k-1} = (-1)^{\omega(n)} \cdot n / \sqrfree(n), \]
as was to be shown.
\end{proof}
\end{lemma}

\begin{proposition}\label{invariant_module}
Let \(e\in\mathbb{S}\) and let \(M\) be a \(\Z\)-module on which multiplication by integer divisors of \(e\) is injective.
Define \(N=M\tensor_\Z \Z[\mu_e]\) and consider the action \(\psi\) of \((\hat\Z/e\hat\Z)^*\cong\Aut_\catt{\(\Z\)-Alg}(\Z[\mu_e])\) on \(N\) via the second component.
Then multiplication with integer divisors of \(e\) is injective on \(N\) and the natural map \(M\to N\) is injective with as image the \((\hat\Z/e\hat\Z)^*\)-invariant submodule \(N^{(\hat\Z/e\hat\Z)^*}\) of \(N\).
\begin{proof}
Recall that we have a functor \(M\tensor_\Z\us\) on \(\Z\)-modules. 
Then \(\varinjlim_{d \divs e}(M\tensor_\Z \Z[\mu_d])=N\) and the action of \((\hat\Z/e\hat\Z)^*\) on \(M\tensor_\Z\Z[\mu_d]\) commutes with the direct limit.
Assuming the proposition holds for integers \(d\divs e\) in place of \(e\), we get
\[ N^{(\hat\Z/e\hat\Z)^*} = \big(\lim_{\substack{\longrightarrow\\ d\divs e}} (M\tensor_\Z \Z[\mu_d]) \big)^{(\hat\Z/e\hat\Z)^*} = \lim_{\substack{\longrightarrow\\ d\divs e}} \big( (M\tensor_\Z \Z[\mu_d])^{(\hat\Z/d\hat\Z)^*} \big) = \lim_{\substack{\longrightarrow\\ d\divs e}} M = M. \]
For \(c\divs d \divs e\) with \(c\) and \(d\) integer we have an exact sequence \(0\to M\tensor_\Z\Z[\mu_d]\xrightarrow{c\cdot} c(M\tensor_\Z\Z[\mu_d])\) because multiplication by \(c\) is injective, so \(0\to \varinjlim_{d \divs e} (M\tensor_\Z\Z[\mu_d]) \to \varinjlim_{d \divs e} c(M\tensor_\Z\Z[\mu_d])\) is an exact sequence.
It follows that multiplication by \(c\) is injective on \(N\).
It remains to prove the proposition assuming \(e\in\Z_{>0}\).

Note that \(\Z[\mu_e]\) is a free \(\Z\)-module with a basis \(\mathcal{B}=\{1,\zeta_e,\dotsc,\zeta_e^{\varphi(e)-1}\}\) containing \(1\), so the natural map \(M^\mathcal{B}\to N\) is an isomorphism of \(\Z\)-modules.
In particular, \(M \hookrightarrow M \tensor 1 \hookrightarrow N\) and multiplication by \(e\) is injective on \(M\).
Obviously \(M\tensor 1 \subseteq N^{(\hat\Z/e\hat\Z)^*}\) since \((\hat\Z/e\hat\Z)^*\) acts via the second component.
For the reverse inclusion, let \(x\in N^{(\hat\Z/e\hat\Z)^*}\) be given.
Consider as in Lemma~\ref{trace_properties} the trace function \(\Tr=\Tr_{\Q(\mu_e)/\Q}\) and note that \(\Tr(\Z[\mu_e])\subseteq\Z\).
Write \(x=\sum_{b\in\mathcal{B}} x_b \tensor b\) with \(x_b\in M\) and let \(d=\sqrfree(e)\). Then
\begin{align*}
\Tr(\zeta_d) x
&= \sum_{\tau \in \Aut(\mu_e)} \tau(\zeta_d) x 
= \sum_{\tau \in \Aut(\mu_e)} \tau(\zeta_d x )
= \sum_{\substack{\tau \in \Aut(\mu_e)\\b\in\mathcal{B}}}  \tau( x_b \tensor \zeta_d b ) 
= \sum_{b\in\mathcal{B}} \Big( x_b \tensor \sum_{\tau \in \Aut(\mu_e)} \tau( \zeta_d b ) \Big) \\
&= \sum_{b\in\mathcal{B}}  x_b \tensor \Tr( \zeta_d b ) 
= \Big( \sum_{b\in\mathcal{B}}  x_b \Tr( \zeta_d b ) \Big) \tensor 1 \in M\tensor 1.
\end{align*}
It follows that \(\Tr(\zeta_d) x_b = 0\) for all \(b\neq 1\).
We have that \(\Tr(\zeta_d) \divs e\) by Lemma~\ref{trace_properties}, hence multiplication by \(\Tr(\zeta_d)\) is injective on \(M\).
Thus \(x_b=0\) for all \(b\neq 1\), so \(x=x_1\tensor 1\in M\tensor 1\).
Hence \(N^{(\hat\Z/e\hat\Z)^*} = M\tensor 1\), as was to be shown.
\end{proof}
\end{proposition}

\subsection{Diagonalizable maps}

\begin{definition}
Let \(k\) be a commutative ring, let \(M\) be a \(k\)-module and let \(f\in\text{End}_\catt{\(k\)-Mod}(M)\).
For \(\lambda\in k\) and \(N\subseteq M\) we define the {\em eigenspace of \(f\) in \(N\) at eigenvalue \(\lambda\)} to be \(N(f,\lambda)=\{x\in N\,|\,f(x)=\lambda x\}\).
We say \(f\) is {\em \(S\)-diagonalizable} for some subset \(S\subseteq k\) if the natural map \(\bigoplus_{\lambda\in S} M(f,\lambda)\to M\) is an isomorphism.
\end{definition}

\begin{lemma}\label{com_diag}
Let \(f\in\End_\catt{\(k\)-Mod}(M)\) be \(S\)-diagonalizable and \(g\in\End_\catt{\(k\)-Mod}(M)\) be \(T\)-diagonalizable.
Then \(f\) and \(g\) commute if and only if the natural map
\[ \bigoplus_{\lambda\in S, \mu \in T} ( M(f,\lambda) \cap M(g,\mu) ) \to M \tag{1}\]
is an isomorphism.
\begin{proof}
Assume \(f\) and \(g\) commute.
For all \(y\in M(g,\mu)\) we have \((gf)(y)=(fg)(y)=f(\mu y) = \mu f(y)\), hence \(f(y)\in M(g,\mu)\) and \(f M(g,\mu) \subseteq M(g,\mu)\).
Let \(x\in M(f,\lambda)\) be given.
Since \(g\) is \(T\)-diagonalizable, we may uniquely write \(x=\sum_{\mu\in T} x_\mu\) with \(x_\mu\in M(g,\mu)\) for all \(\mu\in T\).
Now \(\sum_{\mu\in T} \lambda x_\mu = \lambda x = f(x) = \sum_{\mu\in T} f(x_\mu)\) and \(f(x_\mu)\in M(g,\mu)\).
Hence by uniqueness of decomposition we have \(f(x_\mu)=\lambda x_\mu\), so \(x_\mu\in M(f,\lambda) \cap M(g,\mu)\).
Therefore \(M(f,\lambda)=\bigoplus_{\mu\in T} ( M(f,\lambda) \cap M(g,\mu) )\).
That (1) is an isomorphism then follows from \(S\)-diagonalizability of \(f\).

Assume (1) is an isomorphism. Then for \(x\in M(f,\lambda)\cap M(g,\mu)\) we have \(f(g(x))=\lambda\mu x = g(f(x))\), so by \(k\)-linearity we have \(f(g(x))=g(f(x))\) for all \(x\in \sum_{\lambda\in S,\,\mu\in T} (M(f,\lambda)\cap M(g,\mu))=M\). Thus \(f\) and \(g\) commute.
\end{proof}
\end{lemma}

\begin{lemma}\label{lin_ind_eigenspaces}
Let \(M\) be a \(k\)-module and let \(S\subseteq k\) be such that multiplication by \(\lambda-\mu\) is injective on \(M\) for all pairwise distinct \(\lambda,\mu\in S\). 
Then for all \(f\in\End_\catt{\(k\)-Mod}(M)\) the natural map \(\bigoplus_{\lambda \in S} M(f,\lambda) \to M\) is injective.
\begin{proof}
Assume the map is not injective. Then there exists a minimal finite non-empty \(I\subseteq S\) such that \(\sum_{\lambda\in I} m_\lambda=0\) for some non-zero \(m_\lambda\in M(f,\lambda)\).
Let \(\mu\in I\) and consider \(J=I\setminus\{\mu\}\).
Then \(0=\mu m_\mu-f(m_\mu)= \sum_{\lambda\in J} ( f(m_\lambda) - \mu m_\lambda ) = \sum_{\lambda\in J} (\lambda-\mu)m_\lambda \).
By minimality of \(I\) we must have that \((\lambda-\mu)m_\lambda = 0\) for some \(\lambda \in J\), hence \(m_\lambda = 0\), a contradiction.
Thus the map is injective.
\end{proof}
\end{lemma}

\begin{corollary}\label{cor:composition_diagonalizable}
Let \(f,g\in \End_\catt{\(k\)-Mod}(M)\) be respectively \(S\)-diagonalizable and \(T\)-diagonalizable and let \(U=\{\lambda\mu\,|\,\lambda\in S,\,\mu\in T\}\).
If \(f\) and \(g\) commute and \(U\) satisfies the conditions to Lemma~\ref{lin_ind_eigenspaces}, then \(f\circ g\) is \(U\)-diagonalizable.
\begin{proof}
By Lemma~\ref{lin_ind_eigenspaces} the natural map \(\varphi:\bigoplus_{\nu \in U} M(fg,\nu) \to M\) is injective.
Note \(M(f,\lambda)\cap M(g,\mu)\subseteq M(fg,\lambda\mu)\) for all \(\lambda\in S\) and \(\mu\in T\), so by Lemma~\ref{com_diag} we have \(\sum_{\nu\in U} M(fg,\nu) \supseteq \sum_{\lambda\in S,\, \lambda\in T}  M(f,\lambda)\cap M(g,\mu) = M \), so \(\varphi\) is surjective.
We conclude that \(fg\) is \(U\)-diagonalizable.
\end{proof}
\end{corollary}

\begin{remark}\label{rem:diag_is_decomp}
Each \(S\)-diagonalizable \(f\in\End_{\catt{\(k\)-Mod}}(M)\) gives a decomposition \(\{ M(f,s) \}_{s\in S}\) of \(M\) and conversely each decomposition \(\{E_s\}_{s\in S}\) of \(M\) gives a morphism that sends \(x\in E_s\) to \(sx\).
If we assume multiplication by \(\lambda-\mu\) is injective for all distinct \(\lambda,\mu\in S\subseteq k\), then this is a bijective correspondence.
The monoid \(G=\End_\texttt{Set}(S)\) acts on the set of \(S\)-decompositions of \(M\) given by \((g,\{E_s\}_{s\in S})\mapsto g_*\{E_s\}_{s\in S}=\{\sum_{t\in g^{-1}s} E_t\}_{s\in S}\) and consequently \(G\) acts on the set \(D\subseteq \End_{\catt{\(k\)-Mod}}(M)\) of \(S\)-diagonalizable morphisms.
If \(S=\mu_e\), then \(\hat{\Z}/e\hat{\Z}\hookrightarrow G\) acts on \(D\), which we will write as \(f^a\) for \(f\in D\) and \(a\in \hat{\Z}/e\hat{\Z}\).
\end{remark}

\subsection{Cyclic group gradings}\label{sec:cyclic_group_grad}

In this section we take \(e\in\mathbb{S}\) and let \(R\) be a not necessarily commutative algebra over a commutative ring \(k\) such that multiplication in \(R\) by integer divisors of \(e\) is injective.
If \(k\) is a \(\Z[\mu_e]\)-algebra we can interpret \(\mu_e\)-gradings of \(R\) as \(\mu_e\)-diagonalizable \(k\)-algebra automorphisms as Proposition~\ref{prop:tensored_correspondence} will show.
In the remainder of this section we work towards a generalization for when \(k\) is not necessarily a \(\Z[\mu_e]\)-algebra.

\begin{proposition}\label{prop:tensored_correspondence}
Assume \(k\) is a \(\Z[\mu_e]\)-algebra. Then we have mutually inverse bijections
\begin{align*}
\mkern-48mu
\begin{array}{lrcl}
&\{\text{\em \(\mu_e\)-gradings of \(R\) as a \(k\)-algebra}\}&\leftrightarrow& \{\sigma\in\Aut_\catt{\(k\)-Alg}(R)\,|\, \text{\em \(\sigma\) is \(\mu_e\)-diagonalizable}\}  \\
\text{given by} & (R,\mu_e,\{R_\zeta\}_{\zeta\in\mu_e}) &\mapsto& \bigg(  r_\zeta \in R_\zeta \mapsto \zeta r_\zeta \bigg) \\
\text{and} & (R,\mu_e,\{R(\sigma,\zeta)\}_{\zeta\in\mu_e}) &\mapsfrom& \sigma.
\end{array}
\end{align*}
\begin{proof}
First we show both maps are well-defined.

(\(\leftarrow\)) By \(\mu_e\)-diagonalizability of \(\sigma\) we have that \(\{R(\sigma,\zeta)\}_{\zeta\in\mu_e}\) is a decomposition of \(R\) as a \(k\)-module.
Because \(\sigma(1)=1\) we have \(1\in R(\sigma,1)\). 
For \(x\in R(\sigma,\zeta)\) and \(y\in R(\sigma,\xi)\) we have \(\sigma(xy)=\sigma(x)\sigma(y)=(\zeta\xi)(xy)\) so also \(R(\sigma,\zeta)\cdot R(\sigma,\xi)\subseteq R(\sigma,\zeta\xi)\).
Hence \((R,\mu_e,\{R(\sigma,\zeta)\}_{\zeta\in\mu_e})\) is a grading of \(R\) as a \(k\)-algebra.

(\(\rightarrow\)) Let \(\sigma\) be the \(k\)-module homomorphism given by \(r_\zeta\mapsto \zeta r_\zeta\) for all \(r_\zeta \in R_\zeta\), which is well defined since \(\{R_\zeta\}_{\zeta\in\mu_e}\) is a decomposition of \(R\).
For \(r_\zeta\in R_\zeta\) and \(r_\xi\in R_\xi\) we have \(r_\zeta r_\xi\in R_{\zeta\xi}\) so \(\sigma(r_\zeta r_\xi)=\zeta\xi r_\zeta r_\xi = \sigma(r_\zeta)\sigma(r_\xi)\) and \(\sigma(1)=1\) as \(1\in R_1\). Hence \(\sigma\) is a \(k\)-algebra homomorphism by linearity.
For distinct \(\zeta,\xi\in \mu_e\) multiplication by \(\zeta-\xi\) is injective on \(R\) by Lemma~\ref{lem:root_diff_division}.2.
Then by Lemma~\ref{lin_ind_eigenspaces} the natural map \(\bigoplus_{\zeta\in\mu_e} R(\sigma,\zeta) \to R\) is injective.
Because \(R_\zeta\subseteq R(\sigma,\zeta)\) for all \(\zeta\) we have \(R=\sum_{\zeta\in\mu_e} R_\zeta \subseteq \sum_{\zeta\in\mu_e} R(\sigma,\zeta)\), so the map must in fact be an isomorphism.
We conclude that \(\sigma\) is \(\mu_e\)-diagonalizable.

It follows trivially from Remark~\ref{rem:diag_is_decomp} that the given maps are each other's inverse.
\end{proof}
\end{proposition}

We now consider the ring \(R'=R\tensor_\Z\Z[\mu_e]\), which is a \(k'=k\tensor_\Z\Z[\mu_e]\)-algebra.
Multiplication by integer divisors of \(e\) is injective on \(R'\) by Proposition~\ref{invariant_module}, so \(R'\) satisfies the conditions to Proposition~\ref{prop:tensored_correspondence}.
For \(a\in(\hat{\Z}/e\hat{\Z})^*\) let \(\tau_a\) be the the image of \(a\) in \(\Aut_{k\texttt{-Alg}}(R')\) under the action defined in Proposition~\ref{invariant_module}.
Now \((\hat{\Z}/e\hat{\Z})^*\) has a left action on \(\Aut_{k'\texttt{-Alg}}(R')\) given by \((a,\sigma)\mapsto {}^a\sigma=\tau_a \sigma \tau_a^{-1}\).
Furthermore, \((\hat{\Z}/e\hat{\Z})^*\) has a right action on the set of \(\mu_e\)-diagonalizable \(k'\)-module homomorphisms by Remark~\ref{rem:diag_is_decomp}.
Hence we may define the following. 

\begin{definition}\label{def:X}
We define 
\[X_e(R) = \{ \sigma\in\Aut_{k'\texttt{-Alg}}(R') \,|\, \sigma\text{ is \(\mu_e\)-diagonalizable and } (\forall a\in (\hat{\Z}/e\hat{\Z})^*)\, {}^a\sigma = \sigma^a \}.\]
\end{definition}

We will show that the \(\sigma\in X_e(R)\), under the bijection of Proposition~\ref{prop:tensored_correspondence} applied to \(R'\), correspond to the \(k'\)-algebra gradings of \(R'\) obtained from gradings \((R,\mu_e,\{R_\zeta\}_{\zeta\in\mu_e})\) of \(R\) as \((R',\mu_e,\{R_\zeta\tensor_\Z\Z[\mu_e]\}_{\zeta\in\mu_e})\).

\begin{lemma}\label{extended_grading}
Let \(\sigma\in X_e(R)\) be given. Then 
\begin{enumerate}[topsep=0pt,itemsep=-1ex,partopsep=1ex,parsep=1ex,label={\em(\arabic*)}]
\item For each \(a\in(\hat{\Z}/e\hat{\Z})^*\) and \(\zeta\in\mu_e\) we have \(\tau_a R'(\sigma,\zeta) = R'(\sigma,\zeta)\).
\item For each \(b\in\hat{\Z}/e\hat{\Z}\) also \(\sigma^b\in X_e(R)\).
\end{enumerate}
\begin{proof}
(1) Let \(x\in R'(\sigma,\zeta)\) and let \(d\divs e\) be an integer such that \(x,\tau_a(x)\in D = R'(\sigma^d,1)\).
Let \(b\in\Z_{>0}\) be such that \(ab\equiv 1 \mod d\). Then
\[ \zeta \tau_a(x) = \tau_a(\zeta^{b}x) = \tau_a \sigma^b (x) = \sigma^a \tau_a \sigma^{b-1}(x) = \dotsm = \sigma^{ab} \tau_a(x) = \sigma \tau_a(x)\]
as \(\sigma^{ab}|_D=\sigma|_D\).
Thus \(\tau_a(x)\in R'(\sigma,\zeta)\) and \(\tau_a R'(\sigma,\zeta)\subseteq R'(\sigma,\zeta)\).
From this fact applied to \(\tau_a^{-1}=\tau_{a^{-1}}\) instead we obtain \(R'(\sigma,\zeta)=\tau_a\tau_{a^{-1}} R'(\sigma,\zeta) \subseteq \tau_a R'(\sigma,\zeta)\), proving equality.

(2) Now let \(b\in\hat{\Z}/e\hat{\Z}\) be given.
By definition \(\sigma^b\) is diagonalizable with \(R'(\sigma^b,\zeta)=\sum_{\xi : \xi^b=\zeta} R'(\sigma,\xi) \) for all \(\zeta\in\mu_e\).
For all \(\zeta\in\mu_e\), \(a\in(\hat{\Z}/e\hat{\Z})^*\) and \(x\in R'(\sigma,\zeta)\) we have using (1) that \(\tau_a \sigma^b(x) = \tau_a( \zeta^b x ) = \zeta^{ab} \tau_a(x) = \sigma^{ab} \tau_a(x)\), so \(\tau_a\sigma^b=\sigma^{ab}\tau_a\) by linearity. Hence \(\sigma^b\in X_e(R)\).
\end{proof}
\end{lemma}

A consequence of (2) which we will use later is that \(X_e(R)\) is closed under exponentiation and taking inverses. We may now prove the following.

\begin{lemma}\label{lem:grading_lifts_to_tensor}
Let \(Y\) be the set of \(k'\)-algebra gradings \((R',\Gamma,\{R_\gamma'\}_{\gamma\in\Gamma})\) such that \(\tau_a R_\gamma'=R_\gamma'\) for all \(a\in(\hat\Z/e\hat\Z)^*\).
Then we have mutually inverse bijections
\begin{align*} \mkern-150mu
\begin{array}{lrcl}
&\{\text{\em gradings of \(R\) as \(k\)-algebra}\}&\leftrightarrow& Y  \\
\text{given by}\mkern53mu & (R,\Gamma,\{R_\gamma\}_{\gamma\in\Gamma}) &\mapsto& (R',\Gamma,\{R_\gamma\tensor_\Z\Z[\mu_e]\}_{\gamma\in\Gamma})\\
\text{and}&(R,\Gamma,\{R_\gamma'\cap R\}_{\zeta\in\mu_e}) &\mapsfrom& (R',\Gamma,\{R_\gamma'\}_{\gamma\in\Gamma}).
\end{array}
\end{align*}
\begin{proof}
We show both maps are well-defined.

(\(\rightarrow\)) For each grading \(\overline{R}\) its image \(\overline{R'}\) is a grading of \(R'\) simply by the distributivity of the tensor product over the direct sum. That it satisfies the additional conditions follows from the fact that \(\tau_a\) is an isomorphism and acts on the second component of each homogeneous component. Hence the mapping is well-defined.

(\(\leftarrow\)) 
Let \(G=(\hat\Z/e\hat\Z)^*\). For each grading \(\overline{R'}=(R',\Gamma,\{R_\gamma'\}_{\gamma\in\Gamma})\in Y\) the group \(G\) respects the homogeneous components, so by Proposition~\ref{invariant_module} we have
\[ R = (R')^G = \big( \bigoplus_{\gamma\in\Gamma} R_\gamma' \big)^G = \bigoplus_{\gamma\in\Gamma} (R_\gamma')^G = \bigoplus_{\gamma\in\Gamma} (R_\gamma'\cap R).\]
Thus the image of \(\overline{R'}\) is in fact a \(k\)-algebra grading of \(R\) and the mapping is well-defined.

That the maps are each other's inverse follows readily from the fact that the natural map \(R\to R'\) is injective by Proposition~\ref{invariant_module}.
\end{proof}
\end{lemma}

\begin{theorem}\label{cyclic_correspondence}
Let \(e\in\mathbb{S}\), \(k\) a commutative ring and \(R\) a \(k\)-algebra in which multiplication by integer divisors of \(e\) is injective.
Then we have mutually inverse bijections
\begin{align*} \mkern-65mu
\begin{array}{lrcl}
&\{\text{\em \(\mu_e\)-gradings of \(R\) as \(k\)-algebra}\}&\leftrightarrow& X_e(R) \quad\quad\quad\quad\quad\quad\quad\quad  \\
\text{given by}\mkern43mu & (R,\mu_e,\{R_\zeta\}_{\zeta\in\mu_e}) &\mapsto& \bigg(  r_\zeta\tensor 1 \in R_\zeta\tensor_\Z \Z[\mu_e] \mapsto  r_\zeta\tensor \zeta \bigg) \\
\text{and}&(R,\mu_e,\{R(\sigma,\zeta)\}_{\zeta\in\mu_e}) &\mapsfrom& \sigma.
\end{array}
\end{align*}
\begin{proof}
Let \(Y\) be as in Lemma~\ref{lem:grading_lifts_to_tensor}.
Using Lemma~\ref{extended_grading}, one easily verifies that the image of \(Y\) in \(\Aut_{k\tensor_\Z\Z[\mu_e]\catt{-Alg}}(R\tensor_\Z\Z[\mu_e])\) under the correspondence of Proposition~\ref{prop:tensored_correspondence} is precisely \(X_e(R)\).
Now the given maps are simply a composition of those of Proposition~\ref{prop:tensored_correspondence} and Lemma~\ref{lem:grading_lifts_to_tensor}, making them well-defined and each other's inverse.
\end{proof}
\end{theorem}

\begin{corollary}\label{cor:one_index_2_grad}
Assume \(k\) is a field with \(\characteristic(k)\neq 2\). 
If \(k\subseteq L\) is a field extension of degree \(2\), then \(L\) has exactly one non-trivial efficient \(k\)-algebra grading.
\begin{proof}
By Proposition~\ref{prop:basic_grad_props} the \(k\)-algebra \(L\) can (up to unique isomorphism) only be non-trivially efficiently graded with \(\mu_2\).
Since \(\characteristic(k)\neq 2\) the extension \(k\subseteq L\) is Galois.
Hence \(X_2(L)=\Aut_{k\catt{-Alg}}(L)\) contains \([L:k]=2\) elements, one of which induces the trivial grading. 
\end{proof}
\end{corollary}

\begin{example}
The correspondence from Theorem~\ref{cyclic_correspondence} is not powerful enough to find \(\mu_p\)-gradings of fields of characteristic \(p\).
Let \(\F_{p^n}\) be the unique finite field extension of degree \(n\) of \(\F_p=\Z/p\Z\).
If \((\F_{p^n},\mu_p,\{R_\zeta\}_{\zeta\in\mu_p})\) is a non-trivial grading, then \(p\divs n\) and \(R_1=\F_{p^{n/p}}\) by Proposition~\ref{prop:basic_grad_props}.
For \(\alpha\in R_\zeta\) we have \(\alpha^p\in R_1=\F_{p^{n/p}}\), hence \(\alpha\in \F_{p^{n/p}}\) as the Frobenius automorphism \(x\mapsto x^p\) of \(\F_{p^n}\) sends \(\F_{p^{n/p}}\) to itself. 
Thus \(\alpha\in R_\zeta\cap R_1\) implies \(\alpha=0\) if \(\zeta\neq 1\), hence \(R_\zeta=0\), which is a contradiction.
Thus no non-trivial \(\mu_p\)-grading of \(\F_{p^n}\) exists.
\end{example}

\begin{example}
Consider \(K_n=\Q(\!\sqrt[n]{p})\) for \(p\) prime and \(n\geq 1\) and note that \([K_n:\Q]=n\) since \(X^n-p\) is an irreducible polynomial over \(\Q\) by Eisenstein's Criterion.
Hence to compute all cyclic gradings of \(K_n\) it suffices to consider \(\mu_n\)-gradings.

Assume \(n\) is odd or \(\sqrt{p}\not\in \Q(\zeta_n)\). 
We then have that \(K_n':=K_n\tensor_\Z \Z[\mu_n]\cong \Q(\sqrt[n]{p},\zeta_n)\) by Theorem A of Jacobson and V\'elez \cite{Jacobson1990}.
Hence \(\Aut_{\Z[\mu_n]\catt{-Alg}}(K_n')=\{ \sigma_i = (\sqrt[n]{p} \mapsto \zeta_n^i \cdot \sqrt[n]{p}) \,|\, i\in\Z/n\Z \}\).
The grading corresponding to \(\sigma_1\) by Theorem~\ref{cyclic_correspondence} is \(\overline{R}=(K_n,\mu_n,\{ p^{i/n}\cdot \Q \}_{\zeta_n^i\in\mu_n})\) and the grading corresponding to \(\sigma_i\) is simply \(f_*\overline{R}\) with \(f:\mu_n\to\mu_n\) given by \(\zeta\mapsto\zeta^i\).
It follows that these are all efficient gradings of \(K_n\) and thus \(\overline{R}\) is the universal abelian group grading.
In fact, by Proposition~\ref{prop:basic_grad_props} it is the universal grid-grading.
\end{example}

\subsection{Abelian torsion group gradings}

In this section we generalize Theorem~\ref{cyclic_correspondence} to gradings with arbitrary abelian \(e\)-torsion groups.
This relies on a theorem of Pontryagin on dual groups, as defined in Definition~\ref{def:dual}.

\begin{theorem}[Pontryagin, Theorem 1.7.2 of \cite{rudin2017fourier}]\label{thm:pontryagin}
There is a natural isomorphism \(\Phi:\id_\catt{LCA}\to\hhatf\) with the component \(\Phi(\Gamma)\) at \(\Gamma\in\obj(\catt{LCA})\) given by \(\gamma\mapsto \ev_\gamma\), where \(\ev_\gamma\in\widehhat{\Gamma}\) is evaluation map \(\chi\mapsto\chi(\gamma)\). \qed
\end{theorem}

\begin{corollary}\label{cor:dualization_is_hom_isomorphism}
The functor \(\hatf\) is self-adjoint, i.e. there is a natural isomorphism between the two bifunctors \(\catt{LCA}^2\to\catt{Set}\) given by \((\Gamma,\Delta)\mapsto\Hom_\catt{LCA}(\widehat{\Gamma},\Delta)\) respectively \((\Gamma,\Delta)\mapsto\Hom_\catt{LCA}(\widehat{\Delta},\Gamma)\). \qed
\end{corollary}

\begin{remark}\label{rem:def_thm}
Let \(e\in\mathbb{S}\), let \(k\) be a commutative ring and \(R\) be a \(k\)-algebra on which multiplication by integer divisors of \(e\) is injective.
We then equip \(R'=R\tensor_\Z\Z[\mu_e]\) with the discrete topology and \(A=\Aut_{k\tensor_\Z\Z[\mu_e]\catt{-Alg}}(R')\) with the compact-open topology, making it a topological group.
Note that since \(R'\) is discrete, \(\{v(x,y)\,|\,x,y\in R'\}\) with \(v(x,y)=\{f\in A\,|\, f(x)=y\}\) is a subbase of \(A\).
\end{remark}

\begin{proposition}\label{prop:direct_grading}
Let \(e\), \(k\) and \(R\) be as in Remark~\ref{rem:def_thm}.
Let \(\Delta\subseteq X_e(R)\) be a compact abelian subgroup of \(A\).
Then \(\overline{R}=(R,\widehat{\Delta},\{R_\chi\}_{\chi\in\widehat{\Delta}})\) with \(R_\chi = \bigcap_{\delta\in\Delta} R(\delta,\chi(\delta))\) is a grading.
\begin{proof}
If we replace \(\widehat{\Delta}\) by \(\mu_e^\Delta\) in the definition of \(\overline{R}\), then the resulting grading is simply the joint grading of the gradings corresponding to the \(\delta\in\Delta\) under the bijection of Theorem~\ref{cyclic_correspondence}.
To show \(\overline{R}\) is a grading it suffices by Lemma~\ref{lem:intersect_grading} to verify that the natural map \(f:\bigoplus_{\chi\in\widehat{\Delta}} R_\chi' \to R'\) is surjective, with \(R_\chi'=\bigcap_{\delta\in\Delta} R'(\delta,\chi(\delta))\).

Let \(x\in R'\) be given. 
Then \(\{v(x,y)\,|\, y\in R'\}\) is an open cover of \(\Delta\) so by compactness it has a finite subcover.
Hence \(\Delta x\) is finite and the evaluation map \(\ev_x:\Delta \to \Delta x\) has only finitely many fibers.
For each \(y\in\Delta x\) take a representative \(\delta_y\in \ev_x^{-1}\{y\}\) and note that since the \(\delta_y\) commute pair-wise we may write \(x=\sum_{a\in\mu_e^{\Delta x}} x_a \) with \(x_a\in S_a=\bigcap_{y\in\Delta x}R'(\delta_y,a(y))\) by Lemma~\ref{com_diag}.
It now suffices to show that \(x_a\in \sum_{\chi\in\widehat{\Delta}}R_\chi'\).
Fix some \(a\in\mu_e^{\Delta x}\) such that \(x_a\neq 0\) and let \(\chi=a\circ\ev_x:\Delta\to\mu_e\).
We have \((\chi(\delta\gamma)-\chi(\delta)\chi(\gamma))x_a=(\delta\gamma)(x_a)-\delta(\gamma(x_a))=0\) hence \(\chi(\delta\gamma)=\chi(\delta)\chi(\gamma)\) by Lemma~\ref{lem:root_diff_division}.2.
As \(\chi(1)=1\) we have that \(\chi\) is a group homomorphism.
Similarly, since \(\chi^{-1}\{\zeta\}=v(x_a,\zeta x_a)\cap\Delta\) is open for all \(\zeta\in\mu_e\) we even have \(\chi\in\widehat{\Delta}\).
Each \(\delta\in\Delta\) commutes with all representatives \(\delta_v\) so \(\delta(x_a)\in S_a\).
Then from \(\sum_{a\in \mu_e^{\Delta x}} \delta(x_a) =\delta(x)=\delta_{\delta(x)}(x) = \sum_{a\in\mu_e^{\Delta x}} \chi(\delta) \cdot x_a\) it follows that \(x_a\in R_\chi\). Hence \(f\) is surjective. We conclude that \(\overline{R}\) is a grading.
\end{proof}
\end{proposition}

\begin{remark}\label{rem:cyclic_morphism}
Theorem~\ref{cyclic_correspondence} can be considered a special case of Proposition~\ref{prop:direct_grading}.
Let \(e\), \(k\) and \(R\) be as in Remark~\ref{rem:def_thm}.
For any \(\sigma\in X_e(R)\) consider \(\Delta = \{ \sigma^a \,|\, a\in\hat{\Z}/e\hat{\Z} \}\), which is a compact abelian subgroup of \(X_e(R)\) by Lemma~\ref{extended_grading}.2.
Then the map \(f_\sigma\in\Hom_\catt{LCA}(\widehat{\mu}_e,\Delta)\) given by \((\zeta\mapsto\zeta^a)\mapsto \sigma^a\) is continuous and surjective.
Corollary~\ref{cor:dualization_is_hom_isomorphism} applied to \(\widehat{f_\sigma}\) gives a injective morphism of groups \(\varphi:\widehat{\Delta}\to \mu_e\).
The \(\mu_e\)-grading corresponding to \(\sigma\) by Theorem~\ref{cyclic_correspondence} is simply \(\varphi_*\overline{R}\), with \(\overline{R}\) the \(\widehat{\Delta}\)-grading corresponding to \(\Delta\) by Proposition~\ref{prop:direct_grading}.
\end{remark}

\begin{theorem}\label{thm:grad_is_morph}
Let \(e\), \(k\), \(R\) and \(A\) be as in Remark~\ref{rem:def_thm}.
Write \(\catt{Ab}_e\subseteq\catt{Ab}\) for the full subcategory of \(e\)-torsion groups.
Let \(G:\catt{Ab}_e\to\catt{Set}\) be the functor that sends \(\Gamma\) to the set of \(\Gamma\)-gradings and \(H:\catt{Ab}_e\to\catt{Set}\) the functor that sends \(\Gamma\) to \(\{\varphi\in\Hom_\catt{TGrp}(\widehat{\Gamma},A)\,|\, \im(\varphi)\subseteq X_e(R)\}\).
Then we have a natural isomorphism \(\Phi:G\to H\) where the component at \(\Gamma\in\obj(\catt{Ab}_e)\) is given by
\begin{align*}
\Phi(\Gamma):\{\text{\(\Gamma\)-gradings of \(R\) as \(k\)-algebra}\} &\to \{\varphi\in\Hom_\catt{TGrp}(\widehat{\Gamma},A)\,|\, \im(\varphi)\subseteq X_e(R)\}\\
\overline{R}=(R,\Gamma,\{R_\gamma\}_{\gamma\in\Gamma}) &\mapsto \Big( \chi \mapsto (  x\in R_\gamma \mapsto \chi(\gamma)\cdot x ) \Big) \\
\bigg(R,\Gamma,\Big\{  \bigcap_{\chi\in\widehat{\Gamma}} R(\varphi(\chi),\chi(\gamma))  \Big\}_{\gamma\in\Gamma} \bigg) &\mapsfrom \varphi.
\end{align*}
\begin{proof}
We show that for fixed \(\Gamma\in\obj(\catt{Ab}_e)\) the maps \(G(\Gamma)\to H(\Gamma)\) and \(H(\Gamma)\to G(\Gamma)\) are well-defined, as it then easily follows that the maps are mutually inverse and that \(\Phi\) is a natural transformation.

(\(\rightarrow\)) Note that since \(\Gamma\) is \(e\)-torsion we have \(\chi(\Gamma)\subseteq \mu_e\) for all \(\chi\in\widehat{\Gamma}\).
For a \(\Gamma\)-grading \(\overline{R}\) it is then clear that the corresponding \(\varphi\) is a well-defined \(k\tensor\Z[\mu_e]\)-algebra isomorphism with \(\varphi(\chi)\in X_e(R)\) for all \(\chi\in\widehat{\Gamma}\).
It remains to show that it is continuous. 
By Remark~\ref{rem:def_thm} it suffices to show that \(\varphi^{-1}v(x,y)\) is open for all \(x,y\in R'\). 
Write \(x=\sum_{\gamma\in\Gamma} x_\gamma\tensor u_\gamma\) with \(x_\gamma\in R_\gamma\) and \(u_\gamma\in\Z[\mu_e]\) for all \(\gamma\in\Gamma\).
If \(y\not\in \sum_{\gamma\in\Gamma} \mu_e\cdot (x_\gamma \tensor u_\gamma)\), then \(\varphi^{-1}v(x,y)=\emptyset\) is open.
Otherwise we may write \(y=\sum_{\gamma\in\Gamma} x_\gamma \tensor (\zeta_\gamma \cdot u_\gamma) \) for some \(\zeta_\gamma\in\mu_e\) and \(\varphi^{-1}v(x,y)=\bigcap_{\gamma: x_\gamma\neq 0} \varphi^{-1}v(x_\gamma\tensor u_\gamma,x_\gamma\tensor(\zeta_\gamma u_\gamma))\).
Thus without loss of generality we have \(x=x_\gamma\tensor u_\gamma\neq 0\) and \(y=\zeta_\gamma x\).
By Pontryagin the evaluation map \(\ev_\gamma:\widehat{\Gamma}\to \mathbb{T}\) is continuous.
Since \(x\neq 0\) and multiplication by \(\chi(\gamma)-\zeta_\gamma\) is injective when \(\chi(\gamma)\neq\zeta_\gamma\) for all \(\chi\in\widehat{\Gamma}\) by Lemma~\ref{lem:root_diff_division}.2 we get that
\[\varphi^{-1}(x,\zeta_\gamma x)=\{\chi\in\widehat{\Gamma}\,|\, \varphi(\chi)(x)=\zeta_\gamma x\}=\{\chi\in\widehat{\Gamma}\,|\, (\chi(\gamma)-\zeta_\gamma)x=0\} = \ev_\gamma^{-1}(\zeta_\gamma)\]
is open. Hence \(\varphi\) is continuous, as was to be shown.

(\(\leftarrow\)) Let \(\varphi\in H(\Gamma)\) be given and let \(\Delta=\im(\varphi)\), which is a compact abelian group.
Since \(\varphi\) factors through \(\Delta\) we get a morphism \(f:\widehat{\Gamma}\to\Delta\) and thus by Pontryagin a morphism \(g:\widehat{\Delta}\to\Gamma\).
Then \(\overline{S}=(R,\widehat{\Delta},S_\chi)\) with \(S_\chi=\bigcap_{\delta\in\Delta} R(\delta,\chi(\delta))\) is a grading by Proposition~\ref{prop:direct_grading}, and the grading corresponding to \(\varphi\) is precisely \(g_*\overline{S}\).
Hence \(\Phi(\Gamma)\) is well-defined and bijective, as was to be shown.
\end{proof}
\end{theorem}

Theorem~\ref{thm:grad_is_morph} essentially states that each grading of \(R\) with an abelian \(e\)-torsion group can be realized as the joint gradings of several \(\mu_e\)-gradings.
These joint gradings exist precisely when all involved \(\sigma\in X_e(R)\) commute pair-wise, motivating the following theorem.

\begin{theorem}\label{thm:commutative_implies_universal}
Let \(e\in\mathbb{S}\), let \(k\) be a commutative ring and let \(R\) be a \(k\)-algebra on which multiplication by integer divisors of \(e\) is injective.
Then \(R\) has a universal abelian \(e\)-torsion grading if and only if \(X_e(R)\) is a compact abelian group.
If these equivalent conditions hold, then the universal grading is given by Proposition~\ref{prop:direct_grading} applied to \(\Delta=X_e(R)\).
\begin{proof}
By Theorem~\ref{thm:grad_is_morph} the functors \(G\) and \(H\) as defined there are naturally isomorphic.
Thus \(R\) has a universal grading if and only if \(H\) is representable.

Assume \(X_e(R)\) is a compact abelian group.
Then \(\Upsilon = \widehat{X_e(R)}\) is a discrete group such that \(\widehat{\Upsilon}\cong X_e(R)\) by Pontryagin.
Corollary~\ref{cor:dualization_is_hom_isomorphism} now provides the natural isomorphism \(\Hom_\catt{Ab}(\Upsilon,\us)\cong \Hom_\catt{LCA}(\widehat{\us},X_e(R))=H\).
Hence \(H\) is representable.

Assume \(H\) is representable, so that we have an \(\Upsilon\in\obj(\catt{Ab})\) and natural isomorphism \(\Phi:\Hom_\catt{LCA}(\Upsilon,\us)\to H\).
By Yoneda's lemma \(\Phi\) corresponds to \(\varphi=\Phi(\Upsilon)(\id_\Upsilon)\in H(\Upsilon)\), and since \(\Phi\) is an isomorphism we have that for all \(\Gamma\in\obj(\catt{Ab})\) and \(\psi\in H(\Gamma)\) there exists a unique \(f\in\Hom_\catt{Ab}(\Upsilon,\Gamma)\) such that \(\psi = \Phi(\Gamma)(f) (\varphi)\).
Let \(\sigma\in X_e(R)\) be given and define \(\psi:\widehat{\mu_e}\to X_e(R)\) by \((\zeta\mapsto \zeta^a)\mapsto \sigma^a\) as in Remark~\ref{rem:cyclic_morphism}.
Then \(\psi\in H(\mu_e)\), so there exists a \(f\in\Hom_\catt{Ab}(\Upsilon,\mu_e)=\widehat{\Upsilon}\) such that \(\psi=\widehat{f}_*(\varphi)\).
Hence \(\sigma=\psi(\id)=\varphi(\widehat{f}(\id))=\varphi(f)\), so \(X_e(R)= \varphi \widehat{\Upsilon}\).
Since \(\varphi\) is a morphism of topological groups and \(\widehat{\Upsilon}\) is compact, we have that \(X_e(R)\) is a compact abelian group.
\end{proof}
\end{theorem}

If \(\Aut_{k'\catt{-Alg}}(R')\) is finite, as is the case when \(R\) is a finite product of number fields as will follow from Proposition~\ref{aut_wreath_iso}, then it certainly is compact.
In this case \(R\) has a universal abelian group grading if and only if the elements of \(X_\infty(R)\) commute pair-wise.
In Example~\ref{ex:non-compact_counter} we construct a ring \(R\) such that \(X_\infty(R)\) is abelian though no universal abelian torsion group-grading of \(R\) exists.

\subsection{Computing joint gradings}\label{sec:X_to_grad_bij_computation}

The correspondence from Theorem~\ref{thm:grad_is_morph} reduces the problem of Theorem~\ref{AlgAbGrpThm}, finding all gradings with cyclic groups of prime power order, to the problem of finding a specific subset of some automorphism group. 
For this reduction to be of practical use we need to show that computing the grading corresponding to the automorphism can be done in polynomial time.

Let \(k\) be either \(\Z\) or \(\Q\).
If \(R\) and \(S\) are free \(k\)-modules represented by a basis \(\mathcal{A}\) respectively \(\mathcal{B}\), then we have a canonical basis \(\{a\tensor b \,|\, a\in\mathcal{A},\,b\in\mathcal{B}\}\) for \(R\tensor_k S\).
Computationally, if \(A=(a_{hij})_{h,i,j\in I}\) and \(B=(b_{hij})_{h,i,j\in J}\) are structure constants for commutative \(k\)-algebras \(R\) respectively \(S\), we construct \(R\tensor_k S\) as the \(k\)-algebra defined by the structure constants \(A\tensor B := ( a_{hij}\cdot b_{h'i'j'} )_{(h,h'),(i,i'),(j,j')\in I\times J}\).
Note that the length of \(A\tensor B\) is bounded by a polynomial in the length of \((A,B)\).
If \(1\) is in the basis of \(S\), or equivalently for some \(h'\in J\) we have \(b_{h'i'j'}=\mathbbm{1}(i'=j')\), then we can easily construct the natural morphism \(\epsilon:R\to R\tensor_k S\) as the map sending the \(h\)-th basis vector \(e_h\) to \(e_{(h,h')}\) for all \(h\in I\).
Similarly, we can construct a left inverse \(k\)-module homomorphism \(\pi:R\tensor_k S\to R\) of \(\epsilon\) by sending \(e_{(i,i')}\) to \(e_i\) and \(e_{(i,j')}\) to zero when \(j'\neq i'\). 
In particular we can do this for the natural map \(R\to R\tensor_\Z \Z[\mu_e]\) when we are free to choose \(\{1,\zeta,\zeta^2,\dotsc,\zeta^{\varphi(e)-1}\}\) as basis of \(\Z[\mu_e]\).

Using the kernel algorithm from \cite{KernelAlgorithm}, we may compute a basis of the kernel of any \(k\)-module homomorphism \(\sigma\) in polynomial time with respect to the length of the matrix \(\sigma\).
If \(\sigma_1,\dotsc,\sigma_l:M\to N\) are \(k\)-module homomorphism we may compute \(\bigcap_{i=1}^l\ker(\sigma_i)\) as the kernel of the map \(x\mapsto (\sigma_1(x),\dotsc,\sigma_l(x))\), which can then be done in polynomial time with respect to the length of \((\sigma_1,\dotsc,\sigma_l)\).

\begin{proposition}\label{prop:computing_eigenspaces}
Let \(e\in\Z_{\geq 1}\), let \(k\in\{\Z,\Q\}\) and let \(R\) be a free \(k\)-algebra represented by structure constants.
Let \(\sigma_1,\dotsc,\sigma_l\) be a sequence of \(k[\mu_e]\)-algebra endomorphisms of \(R'=R\tensor_kk[\mu_e]\) represented by \(k\)-valued matrices.
Write \(K_{(\zeta_1,\dotsc,\zeta_l)}=\bigcap_{i=1}^l R(\sigma_i,\zeta_i)\) for \(\zeta_1,\dotsc,\zeta_l\in\mu_e^l\).
We can compute the sets \(Z=\{ z\in \mu_e^l \,|\, K_z \neq 0\}\) and \(\mathcal{K}=\{K_z\}_{z\in Z}\) and verify that \(\mathcal{K}\) is a decomposition of \(R\) in polynomial time with respect to the length of the input \((R,e,\sigma_1,\dotsc,\sigma_l)\).
\begin{proof}
Let \(\epsilon:R\to R'\) and \(\pi:R'\to R\) be as before and for \(\zeta\in\mu_e\) let \(A_\zeta:R'\to R'\) be the multiplication by \(\zeta\), all which can be computed in polynomial time.
Then for \(i\leq l\) and \(\zeta\in\mu_e\) define \(M_{i,\zeta}=(\sigma_i-A_\zeta)\cdot \epsilon\) and note that \(\ker(M_{i,\zeta})=R(\sigma_i,\zeta)\).
Let \(Z_m=\{(\zeta_1,\dotsc,\zeta_m)\in\mu_e^m\,|\, \bigcap_{i=1}^{m} R(\sigma_i,\zeta_i) \neq 0 \}\).
Note that \(\# Z_m \) is at most the rank of \(R\), which in turn is bounded by the length of its encoding.
We can thus compute \(Z_{m+1}\) recursively from \(Z_m\) by computing \(\bigcap_{i=1}^{m+1} R(\sigma_i,\zeta_i)\) for all \((\zeta_1,\dotsc,\zeta_m)\in Z_m\) and \(\zeta_{m+1}\in\mu_e\).
Computation of \(\bigcap_{i=1}^{m+1} R(\sigma_i,\zeta_i)\) can be done in polynomial time as discussed before using the kernel algorithm applied to the map \(x\mapsto (M_{i,\zeta_1} x, \dotsc,M_{m+1,\zeta_{m+1}}x)\), so we may compute \(Z_{m+1}\) from \(Z_m\) in polynomial time.
Hence we may compute \(Z=Z_l\) in polynomial time, as well as \(K_z\) for all \(z\in Z\), as was to be shown.
For each \(K_z\) we have a basis, together giving a basis of \(\sum_{z\in Z} K_z\).
To verify \(\sum_{z\in Z} K_z = R\) we simply check whether the matrix associated to the map \(B:\bigoplus_{z\in Z} K_z \to R\) has \(\det(B)\in k^*\), which we may compute using Theorem 6.6 from \cite{determinant}.
\end{proof}
\end{proposition}

\clearpage
\section{Automorphisms and wreaths}

In order to understand gradings of rings in terms of Theorem~\ref{cyclic_correspondence}, we first need to understand the automorphism group \(\Aut_{\Z[\mu_e]}(R\tensor_\Z \Z[\mu_e])\). 
We describe this group using a general construction involving wreaths in the case \(R\) is a product of fields.

\subsection{Wreaths}\label{sec:wreaths}

In this subsection we assume \(\mathcal{C}\) to be a small groupoid.
We will define morphisms on subsets of \(\obj(\mathcal{C})\), which collectively form a new small groupoid, called the power groupoid.
This gives us a group \(\wreath(\mathcal{C})\) which later turns up as the automorphism group of a product of fields.

\begin{definition}\label{power_cat_def}
We define the {\em power category \(2^\mathcal{C}\) of \(\mathcal{C}\)} as follows. Take \(\obj(2^\mathcal{C})=\{X\,|\, X\subseteq \obj(\mathcal{C})\}\) the power set of \(\obj(\mathcal{C})\) and for \(X,Y\in\obj(2^\mathcal{C})\) let
\[ \Hom_{2^\mathcal{C}}(X,Y) = \big\{ \big((\sigma_K)_{K\in X},s\big) \,\big|\, s\in\Iso_\texttt{Set}(X,Y),\ (\forall K\in X)\, \sigma_K\in \Hom_{\mathcal{C}}(K,s(K)) \big\}. \]
For \(\rho=((\rho_K)_{K\in X},r)\in \Hom_{2^\mathcal{C}}(X,Y)\) and \(\sigma=((\sigma_K)_{K\in Y},s)\in \Hom_{2^\mathcal{C}}(Y,Z)\) with \(X,Y,Z\in\obj(2^\mathcal{C})\) define the composition 
\[ \sigma \circ \rho := \big((\sigma_{r(K)}\circ \rho_K)_{K\in X}, s\circ r\big) \in \Hom_{2^\mathcal{C}}(X,Z). \]
For \(X\in \obj(2^\mathcal{C})\) we have \(\id_X=( (\id_K)_{K\in X}, \id )\in \Hom_{2^\mathcal{C}}(X,X)\) as identity.
For \(W,X,Y\in \obj(2^\mathcal{C})\) such that \(W\subseteq X\) and \(\sigma=((\sigma_K)_{K\in X},s)\in\Hom_{2^\mathcal{C}}(X,Y)\) we define \(\sigma|_W = ((\sigma_K)_{K\in W}, s|_W)\in\Hom_{2^\mathcal{C}}(W,s(W))\) as the {\em restriction of \(\sigma\) to \(W\)}, where \(s|_W\in\Iso_\catt{Set}(W,s(W))\) is the restriction of \(s\) to \(W\).
\end{definition}

It is easy to verify that \(2^\mathcal{C}\) is a category. 

\begin{lemma}\label{lem:power_is_groupoid}
The category \(2^\mathcal{C}\) is a groupoid.
Specifically, each \(\sigma=((\sigma_K)_K,s)\in \Hom_{2^\mathcal{C}}(X,Y)\) has an inverse 
\begin{align*}\sigma^{-1} = \big( ((\sigma_{s^{-1}(L)})^{-1} )_{L\in Y}, s^{-1} \big) \in \Hom_{2^\mathcal{C}}(Y,X)\tag*{\qed}\end{align*}
\end{lemma}

If we let \(\catt{SGd}\) be the category of small groupoids with as morphisms the functors which are bijective on the set of objects, then \(\catt{SGd}\) is a groupoid and \(2^{\us}:\catt{SGd}\to\catt{SGd}\) is a functor. Namely, for each \(F\in\Iso_\catt{SGd}(\mathcal{C},\mathcal{D})\) we obtain \(2^F:2^\mathcal{C}\to2^\mathcal{D}\) by applying \(F\) `coordinate-wise'.

\begin{definition}
Define the {\em wreath of \(\mathcal{C}\)} as \(\wreath(\mathcal{C})=\End_{2^\mathcal{C}}(\obj(\mathcal{C}))\), which is a group by Lemma~\ref{lem:power_is_groupoid}.
\end{definition}

For the sake of clarity, we will use the Latin alphabet for permutations in \(\Aut_\catt{Set}(\obj(\mathcal{C}))\) and the Greek alphabet for elements of \(\hom(\mathcal{C})\) in the context of \(\wreath(\mathcal{C})\), as in Definition~\ref{power_cat_def}. For the sake of brevity, we will assume the following questionable correspondence between characters
\[\begin{tabu}{l|lllllll}
\text{Latin} & a & b & o & r & s & t & u \\ \hline
\text{Greek} & \alpha & \beta & \omega & \rho & \sigma & \tau & \upsilon
\end{tabu}\]
such that for example \(\alpha=((\alpha_K)_{K\in\obj(\mathcal{C})},a)\) is implicit when we take \(\alpha\in\wreath(\mathcal{C})\).

\begin{lemma}
If \(\mathcal{C}\) is connected, then for any \(K\in X=\obj(\mathcal{C})\) we have an isomorphism of groups \(\wreath(\mathcal{C})\cong \Aut_{\mathcal{C}}(K)^X \rtimes \Aut_{\catt{Set}}(X)\), where the action of \(\Aut_\catt{Set}( X)\) on \(\Aut_\mathcal{C}(K)^X\) is given by \({}^s(\sigma_L)_{L\in X} = (\sigma_{s^{-1}(L)})_{L\in X}\).
\begin{proof}
Since all \(L\in X\) are isomorphic, we can choose isomorphisms \(f_L\in\Hom_{\mathcal{C}}(L,K)\) for all \(L\in X\).
Then define the map
\begin{align*}
F:\Aut_{\mathcal{C}}(K)^{ X} \rtimes \Aut_\catt{Set}( X )&\to  \wreath(\mathcal{C})\\
 (\sigma_L)_{L\in X} \cdot s &\mapsto \big( ( f_{s(L)}^{-1} \sigma_{s(L)} f_L )_{L\in X}, s \big).
\end{align*}
It is routine verification that this is in fact an isomorphism.
\end{proof}
\end{lemma}

In the general case, we can decompose \(\mathcal{C}\) into connected components.

\begin{proposition}\label{prop:wreath_is_semi-direct_product}
If the elements of \(S\subseteq\obj(2^\mathcal{C})\) are pair-wise disconnected, then \(\Aut_{2^\mathcal{C}}(\bigcup_{X\in S} X)\cong \prod_{X\in S} \Aut_{2^\mathcal{C}}(X)\).
If \(S=\obj(\mathcal{C})/{\cong}\) is the set of isomorphism classes of \(\mathcal{C}\), then
\begin{align*} \wreath(\mathcal{C}) \cong \prod_{X\in S} \Aut_{2^\mathcal{C}}(X) \cong \prod_{[K]\in S} \bigg( \Aut_\mathcal{C}(K)^{[K]} \rtimes \Aut_{\catt{Set}}([K]) \bigg).\tag*{\qed}\end{align*}
\end{proposition}

By definition, an action of a group \(\Gamma\) on an object \(A\) in some category \(\mathcal{D}\) is simply a homomorphism \(\Gamma\to\Aut_{\mathcal{D}}(A)\) in the category of groups. If we take \(A=\mathcal{C}\) to be an object in the category \(\mathcal{D}=\catt{Cat}\) of small categories, then \(\Aut_{\catt{Cat}}(\mathcal{C})\) is the set of all invertible functors \(F:\mathcal{C}\to\mathcal{C}\).
In particular, a group can act on a category.
Note that each \(F\in\Aut(\mathcal{C})\) induces a unique pair of permutations \(\varphi\) of \(\hom(\mathcal{C})\) and \(f\) of \(\obj(\mathcal{C})\).
Hence if \(\Gamma\) acts on \(\mathcal{C}\), it induces an action on the set of objects of \(\mathcal{C}\).
As we will show, a nice property of \(\wreath(\mathcal{C})\) is that it acts on \(\mathcal{C}\), and that any action of \(\Gamma\) on \(\mathcal{C}\) induces an action of \(\Gamma\) on \(\wreath(\mathcal{C})\).

\begin{lemma}\label{lem:inner_morph}
There is a natural action \(\wreath(\mathcal{C})\to\Aut_\catt{SGd}(\mathcal{C})\) given by sending \(\sigma\in\wreath(\mathcal{C})\) to the functor defined by
\[ A \mapsto s(A) \quad\text{and}\quad f \mapsto \sigma_{B} \circ f \circ (\sigma_A)^{-1} \]
for all \(A,B\in\obj(\mathcal{C})\) and \(f\in\Hom_{\mathcal{C}}(A,B)\).\qed
\end{lemma}

If we take \(\mathcal{C}\) to be a group, which is a groupoid with a single object, then \(\wreath(\mathcal{C})\cong\aut(\mathcal{C})\) and the action of Lemma~\ref{lem:inner_morph} is simply conjugation, such that the image of \(\wreath(\mathcal{C})\) in \(\Aut_\catt{SGd}(\mathcal{C})\) equals the group of inner automorphisms of \(\mathcal{C}\).

\begin{remark}\label{rem:induced_action_on_wreath}
If we have a group action \(\varphi:\Gamma\to\Aut_\catt{SGd}(\mathcal{C})\), then we also have an action \((2^{\us})_{\mathcal{C},\mathcal{C}}\circ\varphi\) of \(\Gamma\) on \(2^\mathcal{C}\), with \((2^{\us})_{\mathcal{C},\mathcal{C}}\) the restriction of \(2^{\us}\) to \(\Hom_\catt{SGd}(\mathcal{C},\mathcal{C})\).
This induces an action \(\Gamma\to\wreath(\mathcal{C})\) and in particular ensures \(\wreath(\mathcal{C})\rtimes\Gamma\) is defined.
\end{remark}

Instead of letting \(\Gamma\) act on \(\mathcal{C}\) after taking its wreath, we can also combine \(\mathcal{C}\) and \(\Gamma\) into a single category first.

\begin{definition}
Let \(\Gamma\) be a group with an action on \(\mathcal{C}\). We define the {\em semi-direct product} \(\mathcal{C}\rtimes \Gamma\) of \(\mathcal{C}\) and \(\Gamma\) as the category with \(\obj(\mathcal{C}\rtimes \Gamma)=\obj(\mathcal{C})\) and 
\[ \Hom_{\mathcal{C}\rtimes\Gamma}(A,B) = \{ (\sigma,\gamma) \,|\, \gamma\in\Gamma,\, \sigma\in\Hom_{\mathcal{C}}(\gamma(A),B) \}. \]
Composition for \( (\rho, \gamma) \in \Hom_{\mathcal{C}\rtimes\Gamma}(A,B)\) and \((\sigma, \delta) \in \Hom_{\mathcal{C}\rtimes\Gamma}(B,C)\) is defined as
\[ (\sigma,\delta) \circ (\rho,\gamma) := \big(\sigma \cdot {{}^\delta}\rho, \delta\gamma\big) \in \Hom_{\mathcal{C}\rtimes\Gamma}( A, C). \]
\end{definition}

\begin{remark}\label{rem:wreath_hierarchy_prop}
Let \(\Gamma\) be a group acting on \(\mathcal{C}\).
We have an monomorphism of groupoids \(F:\mathcal{C}\to\mathcal{C}\rtimes \Gamma\) given by \(A\mapsto A\) and \(\sigma\mapsto (\sigma,1)\) for all \(A\in\obj(\mathcal{C})\) and \(\sigma\in\hom(\mathcal{C})\).
By functoriality of \(2^{\us}\) we get a monomorphism \(2^F:2^\mathcal{C}\to2^{\mathcal{C}\rtimes\Gamma}\) and thus a monomorphism \(f=(2^F)_{\obj(\mathcal{C}),\obj(\mathcal{C})}:\wreath(\mathcal{C})\to\wreath(\mathcal{C}\rtimes\Gamma)\).
The action of \(\Gamma\) gives a map \(G:\Gamma\to\Aut_\catt{SGd}(\mathcal{C})\) which induces a map \(g:\Gamma\to\wreath(\mathcal{C}\rtimes\Gamma)\) given by \(\gamma\mapsto ( ((\id_{G(\gamma)(K)},\gamma))_{K\in\obj(\mathcal{C}\rtimes\Gamma)}, G(\gamma)|_{\obj(\mathcal{C}\rtimes\Gamma)} )\).
The maps \(f\) and \(g\) induce a map \(h:\wreath(\mathcal{C})\rtimes\Gamma\to\wreath(\mathcal{C}\rtimes\Gamma)\) such that the following diagram commutes.
\begin{center}
\begin{tikzcd}
	& \wreath(\mathcal{C})\rtimes \Gamma \arrow[dashed]{dd}{h} & \\
\wreath(\mathcal{C}) \arrow[tail]{ur} \arrow[tail]{dr}{f} & & \Gamma \arrow[tail]{ul} \arrow{dl} \\
	& \wreath(\mathcal{C}\rtimes\Gamma) \arrow[no head]{ur}{g}&
\end{tikzcd} 
\end{center}
When \(\mathcal{C}\) is not the empty groupoid \(g\) and thus \(h\) is injective.
\end{remark}

We can also use actions to define invariant categories.

\begin{definition}\label{def:invariant_category}
Let \(\mathcal{C}\) be a small groupoid with \(\Gamma\) a group acting on \(\mathcal{C}\) and consider the induced action of \(\Gamma\) on \(2^\mathcal{C}\).
We define the {\em \(\Gamma\)-invariant of \(\mathcal{C}\)} as the groupoid \(\mathcal{C}^\Gamma\) with \(\obj(\mathcal{C}^\Gamma)=\obj(\mathcal{C})/\Gamma\) and for \(\Gamma A,\Gamma B\in\obj(\mathcal{C}/\Gamma)\)
\[\Hom_{\mathcal{C}^\Gamma}(\Gamma A,\Gamma B)=\{ \sigma\in \Hom_{2^\mathcal{C}}(\Gamma A,\Gamma B) \,|\, (\forall\gamma\in\Gamma)\ {}^\gamma \sigma = \sigma \}.\]
\end{definition}

It follows from the definition that we have a natural inclusion \(\wreath(\mathcal{C}^\Gamma)\rightarrow \wreath(\mathcal{C})\), and its image is precisely \(\wreath(\mathcal{C})^\Gamma\).

\subsection{Automorphisms of products of fields}

Now we show how we can consider the automorphism group of a product of fields as a wreath.

\begin{lemma}\label{basic_comalg_wreath}
Let \(k\) be a commutative ring and let \(E\) be a commutative \(k\)-algebra.
\begin{enumerate}[topsep=0pt,itemsep=-1ex,partopsep=1ex,parsep=1ex,label={\em(\arabic*)}]
\item If \(k\) is a field and \(\dim_k E < \infty\), then \(\minspec E = \spec E = \maxspec E\) and \(\#\spec E <\infty\).
\item If \(E\) is reduced, then the natural \(k\)-algebra homomorphism
\[ \varphi: E \to \prod_{\p\in\minspec E} E/\p\]
is injective.
\item If \(E\) is reduced, \(\minspec E=\maxspec E\) and \(\#\spec E <\infty\), then \(\varphi\) is an isomorphism.
\end{enumerate}
\begin{proof}
(1) Note that \(E\) is an Artinian ring being finitely generated as module over the Artinian ring \(k\). Hence this is a special case of Theorem~2.14 in \cite{Eisenbud}.
(2) The kernel of \(\varphi\) is \(\bigcap_{\p\in\minspec E} \p\), which is the nilradical of \(E\) by Corollary~2.12 in \cite{Eisenbud}, which is trivial when \(E\) is reduced. 
(3) All minimal prime ideals are maximal and hence pair-wise coprime. Thus \(\varphi\) is an isomorphism by (2) and by the Chinese remainder theorem for rings, which is Exercise~2.6 in \cite{Eisenbud}.
\end{proof}
\end{lemma}

Let \(k\) be a commutative ring.
If a \(k\)-algebra \(E\) is isomorphic to a finite product of fields \(\varphi:E\xrightarrow{\sim}\prod_{i\in I} K_i\), then this family \(\{K_i\}_{i\in I}\) is determined up to isomorphism.
Namely, there is a bijection \(f:I\to\spec E\) such that \(E/f(i) \cong K_i\) for all \(i\in I\).
Hence when we say \(E\) is a finite product of fields we simply mean that \(\#\minspec E=\#\spec E=\#\maxspec E<\infty\) and that the natural map \(E\to\prod_{\p\in\spec E} E/\p\) is an isomorphism.
Similarly, each automorphism of \(E\) induces a permutation on \(\spec E\).

Now to such \(E\) we may associate the following groupoid \(\mathcal{C}\): Take \(\obj(\mathcal{C})=\spec E\) and \(\Hom_\mathcal{C}(\m,\n)=\Iso_{k\texttt{-Alg}}(E/\m,E/\n)\) for all \(\m,\n\in\spec E\). 
Each \(\varphi\in\Aut_\catt{\(k\)-Alg}(E)\) induces \(\varphi_\p\in\Iso_\catt{\(k\)-Alg}(E/\p,E/\varphi(\p))\) and \(f\in\Aut_\catt{Set}(\spec E)\).
Hence, we have a map \(\Phi:\Aut_\catt{\(k\)-Alg}(E)\to\wreath(\mathcal{C})\) given by \(\varphi\mapsto ((\varphi_\p)_{\p\in\obj(\mathcal{C})},f)\).
Conversely, each element of \(\wreath(\mathcal{C})\) induces an automorphism of \(\prod_{\p\in\spec E}E/\p\) using the universal property of the direct product, which in turn gives an automorphism of \(E\) by Lemma~\ref{basic_comalg_wreath}.
Thus we get a map \(\Psi:\wreath(\mathcal{C})\to \Aut_\catt{\(k\)-Alg}(E)\).
We leave the routine verification as an exercise to the reader.

\begin{proposition}\label{aut_wreath_iso}
The maps \(\Phi\) and \(\Psi\) are mutually inverse group homomorphisms, giving a natural isomorphism \(\Aut_{\catt{\(k\)-Alg}}(E)\cong \wreath(\mathcal{C})\). \qed
\end{proposition}

\begin{example}\label{ex:aut_wreath}
A special case we consider will be the following.
Let \(E\neq 0\) be a finite dimensional reduced \(\Q\)-algebra and let \(\Q'=\Q(\mu_e)\) for some \(e\in\mathbb{S}\).
Then \(E'=E\tensor_\Q \Q'\) for \(e\in\mathbb{S}\) is both a \(\Q\)-algebra and a \(\Q'\)-algebra.
Moreover, \(E'\) is a finite product of fields by Lemma~\ref{basic_comalg_wreath} as it is finite dimensional over \(\Q'\) and reduced by Theorem~A1.3 in \cite{Eisenbud}.

In the context of Definition~\ref{def:X} we have three relevant groups: \(A=\Aut_{\Q\catt{-Alg}}(E')\), \(B=\{\sigma\in A\,|\,\sigma\Q'=\Q'\}\) and \(C=\Aut_{\Q'\catt{-Alg}}(E')\).
Let \(\Gamma=\Aut_{\Q\catt{-Alg}}(\Q')\) and consider its natural action \(\varphi:\Gamma\to A\) on \(E'\). 
Here \(X_e(E)\subseteq C\) while \(C\) is generally too small to contain \(\varphi(\Gamma)\).
However, we do have \(\varphi(\Gamma)\subseteq B\).
We even have that \(B\cong C\rtimes \Gamma\).

Now construct the groupoids \(\mathcal{C}\), \(\mathcal{D}\) and \(\mathcal{E}\) as in Proposition~\ref{aut_wreath_iso} for \(E'\) as \(\Q'\)-algebra, \(E'\) as \(\Q\)-algebra and \(E\) as \(\Q\)-algebra respectively.
Then \(\mathcal{C}\) is a wide subgroupoid of \(\mathcal{D}\).
We may also construct \(\mathcal{D}\) and \(\mathcal{E}\) using \(\mathcal{C}\), as \(\mathcal{D}\cong\mathcal{C}\rtimes\Gamma\) and \(\mathcal{E}\cong\mathcal{C}^\Gamma\).
Under the same natural inclusions as in Remark~\ref{rem:wreath_hierarchy_prop} for the wreaths we get the following commutative diagram.
\begin{center}
\begin{tikzcd}
C \arrow[tail]{r} \arrow[leftrightarrow]{d} & B \arrow[tail]{r}\arrow[leftrightarrow]{d} & A \arrow[leftrightarrow]{d} \\
\wreath(\mathcal{C}) \arrow[tail]{r} & \wreath(\mathcal{C})\rtimes \Gamma \arrow[tail]{r} & \wreath(\mathcal{C}\rtimes\Gamma)
\end{tikzcd}
\end{center}
\end{example}

\subsection{Conjugacy in wreaths}

In this section we generalize a well-known theorem (Theorem~\ref{basic_conj_thm}), relating conjugacy of permutations to their cycle type, to wreaths.
For the remainder of this section \(\mathcal{C}\) will be a groupoid with only finitely many objects and \(\mathcal{D}\) will be a wide subgroupoid of \(\mathcal{C}\).

\begin{definition}\label{def:old_labeling}
Let \(\Omega\) be a finite set and let \(s\in\Aut_\texttt{Set}(\Omega)\) be a permutation of \(\Omega\).
Define \(c_k(s)\) to be the number of orbits of \(\langle s\rangle\) of length \(k\).
Then we call \(c(s):=(c_1(s),c_2(s),\dotsc)\) the {\em cycle type} of \(s\).
\end{definition}

\begin{theorem}\label{basic_conj_thm}
Let \(r,s\in\Aut_\catt{Set}(\Omega)\).
Then \(r\) and \(s\) are conjugate in \(\Aut_\catt{Set}(\Omega)\) if and only if \(c(r)=c(s)\).
Furthermore, if those equivalent conditions hold the number of \(a\in\Aut_\catt{Set}(\Omega)\) such that \(ar=sa\) is 
\[ \prod_{k\geq 1} \big( c_k(s)! \cdot k^{c_k(s)} \big) \leq (\#\Omega)^{\#(\Omega/\langle s\rangle)}. \tag*{\(\square\)} \]
\end{theorem}

The idea of the proof is the following: For each orbit of \(\langle r\rangle\) one chooses an orbit of \(\langle s\rangle\) to send it to, which must necessarily be of the same length \(k\), and one notes that this can be done in exactly \(k\) ways. 
In a wreath, it is generally not sufficient to consider only the length of the cycle to be able to decide where to send it. 

\begin{definition}\label{def:conjugacy}
For \(K,L\in\obj(\mathcal{C})\) we call \(\rho\in\Aut_\mathcal{C}(K)\) and \(\sigma\in\Aut_\mathcal{C}(L)\) {\em \(\mathcal{D}\)-conjugate}, written \(\rho\sim_\mathcal{D}\sigma\), if there exists an \(\alpha\in\Hom_\mathcal{D}(K,L)\) such that \(\alpha\rho=\sigma\alpha\).
\end{definition}

Note that if \(K=L\) and \(\mathcal{C}=\mathcal{D}\) in the above definition, our definition of conjugacy coincides with that of conjugacy in the group \(\Aut_\mathcal{C}(K)\), hence this is a generalization. For this generalization it is easy to see that \(\sim_\mathcal{D}\) is also an equivalence relation. 

\begin{definition}\label{def:lambda_sigma}
For \(K\in X\subseteq\obj(\mathcal{C})\) and \(\sigma=((\sigma_L)_{L\in X},s)\in\Aut_{2^\mathcal{C}}(X)\) define
\[ \lambda_\sigma(K) = (\sigma^{\#(\langle\sigma\rangle K)})_K= \sigma_{s^{-1}(K)} \circ \sigma_{s^{-2}(K)} \circ \dotsm\circ \sigma_{s(K)}\circ\sigma_K \in \Aut_\mathcal{C}(K).\]
\end{definition}


\begin{lemma}\label{conj_bauble}
Let \(K,L\in X\subseteq\obj(\mathcal{C})\) and \(\sigma\in\Aut_{2^\mathcal{C}}(X)\).
If \(\langle\sigma\rangle K = \langle\sigma\rangle L\), then \(\lambda_\sigma(K)\sim_\mathcal{C}\lambda_\sigma(L)\).
\begin{proof}
Let \(n=\#(\langle\sigma\rangle K)\) and let \(m\geq 0\) be such that \(s^m(K)=L\).
Consider \(\alpha=(\sigma^m)_K \in \Hom_\mathcal{C}(K,L)\).
Then \(\lambda_\sigma(L)\circ\alpha = (\sigma^{n+m})_K=\alpha \circ \lambda_\sigma(K)\), so \(\lambda_\sigma(K)\sim_\mathcal{C}\lambda_\sigma(L)\).
\end{proof}
\end{lemma}

\begin{definition}
We define \(\digamma(\mathcal{C})=\aut(\mathcal{C})/{\sim_\mathcal{C}}\).
For \(X\subseteq \obj(\mathcal{C})\) and \(\sigma\in\Aut_{2^\mathcal{C}}(X)\) the map \(X\xrightarrow{\lambda_\sigma} \aut(\mathcal{C}) \to \digamma(\mathcal{C})\) factors through \(X/\langle\sigma\rangle\) by Lemma~\ref{conj_bauble} and we write \(\Lambda_\sigma\) for the factor \(X/\langle \sigma\rangle\to\digamma(\mathcal{C})\).
\end{definition}

For a fixed \(\sigma\) we treat the map \(\Lambda_\sigma\) as a labeling of the orbits of \(\sigma\). Then analogously to Definition~\ref{def:old_labeling} we may define the type of an element \(\sigma\in\wreath(\mathcal{C})\).

\begin{definition}
For \(\sigma\in\wreath(\mathcal{C})\), \(k\in\Z_{\geq 1}\) and \(\gamma\in\digamma(\mathcal{C})\), let \(c(\sigma,k,\gamma)\) be the number of orbits \(O\) of \(\langle s \rangle\) of length \(k\) and label \(\Lambda_\sigma(O)=\gamma\).
We then define \(c(\sigma)=(c(\sigma,k,\gamma))_{k\in\Z_{\geq1},\gamma\in\digamma(\mathcal{C})}\) to be the {\em type} of \(\sigma\).
\end{definition}

This generalization of cycle type allows us to formulate the right analogue to Theorem~\ref{basic_conj_thm}.
For \(K\in\obj(\mathcal{C})\) and \(\sigma\in\Aut_{\mathcal{C}}(K)\) write \(C_\mathcal{D}(\sigma)=\{\alpha\in\Aut_\mathcal{D}(K)\,|\, \alpha\sigma=\sigma\alpha\}\) for the centralizer of \(\gamma\) in \(\Aut_\mathcal{D}(K)\).

\begin{theorem}\label{conjugacy_theorem}
Let \(\mathcal{C}\) be a finite groupoid and let \(\rho,\sigma\in\wreath(\mathcal{C})\) be given.
Then \(\rho\) and \(\sigma\) are conjugate in \(\wreath(\mathcal{C})\) if and only if \(c(\rho)=c(\sigma)\).
If these equivalent conditions hold, then the number of \(\alpha\in\wreath(\mathcal{C})\) such that \(\alpha\rho=\sigma\alpha\) is
\[\prod_{\substack{k\in\Z_{\geq 1} \\ [\gamma] \in \digamma}} \bigg( c(\rho,k,[\gamma])!\cdot \big(k\cdot\# C_\mathcal{C}(\gamma)\big)^{c(\rho,k,[\gamma])} \bigg) .\]
With \(\mathcal{D}\) a wide subgroupoid of \(\mathcal{C}\), the number of \(\alpha\in\wreath(\mathcal{D})\) such that \(\alpha\rho=\sigma\alpha\) is at most \((\#\aut(\mathcal{D}))^{\#(\obj(\mathcal{D})/\langle \rho\rangle)}\).
\end{theorem}

We prove this theorem using the following lemma.

\begin{lemma}\label{conj_help_lemma}
Let \(X,Y\subseteq\obj(\mathcal{C})\) and let \(\rho\in \Aut_{2^\mathcal{C}}(X)\) and \(\sigma\in \Aut_{2^\mathcal{C}}(Y)\) be such that \(\langle \rho\rangle\) and \(\langle\sigma\rangle\) act transitively on \(X\) respectively \(Y\).
If \(\rho\sim_{2^\mathcal{D}}\sigma\), then \(A=\{\alpha\in \Hom_{2^\mathcal{D}}(X,Y)\,|\, \alpha\rho = \sigma\alpha \}\) contains at most \(\sum_{L\in Y} \# C_\mathcal{D}(\lambda_\rho(L))\) elements, with equality when \(\mathcal{C}=\mathcal{D}\).
\begin{proof}
Assume \(\rho\sim_{2^\mathcal{D}}\sigma\), hence \(\# X=\# Y =:n\).
Fix \(K\in X\) and define 
\[B=\{(L,\omega)\,|\,L\in Y,\, \omega\in \Hom_{\mathcal{D}}(K,L),\,  \omega \circ \lambda_\rho(K) = \lambda_\sigma(L) \circ \omega \}.\]
If \(\lambda_\rho(K) \sim_\mathcal{D} \lambda_\sigma(L)\) for some \(L\in Y\), then there are precisely \(\#C_{\mathcal{D}}(\lambda_\rho(L))\) elements \(\omega\) such that \((L,\omega)\in B\), and otherwise there are none.
Thus \(\# B\leq \sum_{L\in Y} \# C_\mathcal{D}(\lambda_\rho(L))\), with equality when \(\mathcal{C}=\mathcal{D}\) by Lemma~\ref{conj_bauble}.
Now consider the map \(\Phi:A\to B\) given by \(\alpha\mapsto (a(K),\alpha_K)\).
It is well-defined, since for each \(\alpha\in A\) we have \(\alpha_K\lambda_\rho(K)\alpha_K^{-1}=(\alpha\rho^n\alpha^{-1})_L=(\sigma^n)_K=\lambda_\sigma(L)\).
It suffices to show \(\Phi\) is injective in general and a bijection when \(\mathcal{C}=\mathcal{D}\).
For all \(\alpha\in A\) and \(k\geq 0\) we have the relations
\[ a(r^k(K))=s^k(a(K)) \quad \text{and} \quad \alpha_{r^k(K)} = (\sigma^k)_{a(K)} \circ \alpha_K \cdot (\rho^{-k})_{r^k(K)} \tag{1}\]
We claim that given \((L,\omega)\in B\) the equations of (1) with \(k\in\{0,\dotsc,n-1\}\) uniquely define an element \(\alpha\) of \(\Hom_{2^\mathcal{C}}(X,Y)\) satisfying \(a(K)=L\), \(\alpha_K=\omega\) and \(\alpha\rho=\sigma\alpha\), giving a left inverse \(\Psi:B\to\Hom_{2^\mathcal{C}}(X,Y)\) of \(\Phi\).
By transitivity of \(\langle r\rangle\) the map \(a\) is determined entirely in terms of \(a(K)=L\), meaning \(a\) is unique, and by transitivity of \(\langle s\rangle\) in fact \(a\in\Iso_\catt{Set}(X,Y)\).
From \(\alpha_K=\omega\) and \(\langle r\rangle\) being transitive it follows analogously that \(\alpha\) is uniquely and well-defined.
Thus \(\Psi\) is a well-defined injection, proving the general case.

Assume \(\mathcal{C}=\mathcal{D}\). To prove \(\Psi\) is a right inverse of \(\Phi\) it remains to be shown that \(\im(\Psi)\subseteq A\).
Let \((L,\omega)\in B\) and let \(\alpha=\Psi(L,\omega)\). 
Let \(J\in X\), which we may write as \(J=r^m(K)\) for some \(m\in\{0,\dotsc,n-1\}\). 
When \(m\neq n-1\) the relations \(ar(J)=sa(J)\) and \(\alpha_{r(J)}\circ\rho_J = \sigma_{a(J)} \circ \alpha_J\) follow directly from the defining relations (1) with \(k\in\{0,\dotsc,n-1\}\).
If \(m=n-1\) we also have \(ar(J)=a(K)=s^na(K)=sar^m(K)=sa(J)\), hence \(ar=sa\).
Then,
\begin{align*} 
\alpha_{r(J)} \circ \rho_J 
&= \omega \circ \rho_J = \omega \circ \lambda_\rho(K) \circ (\rho^{m})_{r^{m}(K)} = \lambda_\sigma(L) \circ \omega \circ (\rho^{m})_{r^{m}(K)} \\ 
&= \sigma_{a(J)} \circ (\sigma^{m})_{a(K)} \circ \omega \circ (\rho^{m})_{r^{m}(K)} = \sigma_{a(J)} \circ \alpha_{J}.
\end{align*}
Hence \(\alpha\rho=\sigma\alpha\) and \(\alpha\in A\), as was to be shown. Thus \(\Phi\) and \(\Psi\) are mutually inverse.
\end{proof}
\end{lemma}

\noindent \textbf{Proof of Theorem~\ref{conjugacy_theorem}:} 

\noindent
Assume \(\rho=\alpha\sigma\alpha^{-1}\) for some \(\alpha\in\wreath(\mathcal{C})\), then we certainly have \(r=asa^{-1}\).
Hence \(a\) induces a bijection on the orbits of \(r\) and \(s\). By basic algebra paired orbits are equally large, and by Lemma~\ref{conj_bauble} these orbits must have conjugate labels. In particular, \(c(\rho,k,\gamma)=c(\sigma,k,\gamma)\) for all \(k\in\Z_{\geq 0}\) and \(\gamma\in\digamma\), so \(\rho\) and \(\sigma\) are of the same type.

If \(\rho\) and \(\sigma\) are of the same type, we may choose a bijection \(f\) between the set of orbits \(\mathcal{O}\) of \(\rho\) and those of \(\sigma\) that preserves length and label.
Then by Lemma~\ref{conj_help_lemma}, for each \(O\in\mathcal{O}\) there exists \(\alpha_O\in \Hom_{2^\mathcal{C}}(O,f(O))\) such that \(\sigma|_{f(O)}\circ \alpha_O=\alpha_O \circ \rho|_O\).
Then the concatenation \(\alpha:=\prod_{O\in\mathcal{O}} \alpha_O \in \wreath(\mathcal{C})\) satisfies \(\alpha\rho=\sigma\alpha\), so \(\rho\sim_{2^\mathcal{C}}\sigma\).
Note that for different choices of \(f\) and \(\alpha_O\) we obtain different \(\alpha\).
Conversely, each \(\alpha\) induces a bijection \(f\) and parts \(\alpha_O\in \Hom_{2^\mathcal{C}}(O,f(O))\) for all \(O\in\mathcal{O}\).
Basic combinatorics and Lemma~\ref{conj_help_lemma} then show that the number of \(\alpha\) is
\[ \Bigg(\prod_{\substack{k\in\Z_{\geq 0} \\ [\gamma] \in \digamma}} c(\rho,k,[\gamma])! \Bigg) \cdot \Bigg( \prod_{\substack{k\in\Z_{\geq 0} \\ [\gamma] \in \digamma}} ( k \cdot \#C_\mathcal{C}(\gamma))^{c(\rho,k,[\gamma])} \Bigg), \tag{2}\]
as was to be shown. 

Consider a wide subgroupoid \(\mathcal{D}\) of \(\mathcal{C}\). 
Note that \(\sum_{k\in\Z_{\geq 0},\, [\gamma] \in \digamma} c(\rho,k,[\gamma])=\#\mathcal{O}=\#\obj(\mathcal{D})/\langle\rho\rangle\).
Then from Lemma~\ref{conj_help_lemma} we get a upper bound on the number of \(\alpha\in\wreath(\mathcal{D})\) such that \(\alpha\rho=\sigma\alpha\) similarly as for (1), namely
\[ \bigg(\prod_{\substack{k\in\Z_{\geq 0} \\ [\gamma] \in \digamma}} c(\rho,k,[\gamma])! \bigg)  \cdot \bigg( \prod_{O\in\mathcal{O}} \sum_{K\in O} \# \Aut_\mathcal{D}(K) \bigg) \leq \bigg( \sum_{K\in\obj(\mathcal{D})}\#\Aut_\mathcal{D}(K)\bigg)^{\#(\obj(\mathcal{D})/\langle\rho\rangle)},\]
as was to be shown.
\qed 

\begin{corollary}\label{cor:semidirect_conjugacy_theorem}
Let \(\mathcal{D}\) be a finite groupoid and \(\Gamma\) a group acting on \(\mathcal{D}\).
If \(\sigma\gamma,\tau\gamma'\in \wreath(\mathcal{D})\rtimes \Gamma\subseteq \wreath(\mathcal{D}\rtimes \Gamma)\) are \(\wreath(\mathcal{D})\)-conjugate, then \(\gamma=\gamma'\) and the number of \(\alpha\in\wreath(\mathcal{D})\) such that \(\alpha\sigma\gamma=\rho\gamma\alpha\) is at most \((\#\aut(\mathcal{D}))^{\#(\obj(\mathcal{D})/ \langle \sigma\gamma\rangle)}\).
\begin{proof}
If \(\alpha\sigma\gamma=\rho\gamma'\alpha\) for some \(\alpha\in\wreath(\mathcal{D})\), then \(\gamma=\pi_\Gamma(\alpha\sigma\gamma)=\pi_\Gamma(\rho\gamma'\alpha)=\gamma'\).
For the second claim we simply apply Theorem~\ref{conjugacy_theorem} to the categories \(\mathcal{C}=\mathcal{D}\rtimes \Gamma\) and \(\mathcal{D}\subseteq\mathcal{C}\).
\end{proof}
\end{corollary}

\subsection{Computation} \label{sec:wreath_conjugacy_computation}

If we want to do computations on groupoids and wreaths we first need to specify an encoding.
In our applications, the objects of \(\mathcal{C}\) are number fields as in Example~\ref{ex:aut_wreath}.
For both the objects and morphisms of \(\mathcal{C}\) we already made a choice of how to encode them in Section~\ref{sec:introduction}, namely by structure constants and matrices respectively.
We then encode \(\mathcal{C}\) as a sequence \(K_1,\dotsc,K_n\) of objects followed by \(n\times n\) sequences of morphisms \(\sigma_{ij1},\dotsc,\sigma_{ijm_{ij}}\in \Hom_\mathcal{C}(K_i,K_j)\).
Then \(\sigma\in\wreath(\mathcal{C})\) can be encoded as a permutation \(s\) of \(\{1,\dotsc,n\}\) and a sequence of \(n\) morphisms in the obvious way. 
There may exist some more practical or efficient encodings, we may pre-compute all inverses of the matrices for example or work with multiplication tables instead of actual matrices. 
This simple encoding however already satisfies our requirements, namely that composition and inversion in both \(\mathcal{C}\) and \(\wreath(\mathcal{C})\) can be done in polynomial time with respect to the length of the encoding.

One may verify that Lemma~\ref{conj_help_lemma} implicitly gives a polynomial time algorithm to compute \(A\).
Thus we may turn Theorem~\ref{conjugacy_theorem} into an algorithm as well.

\begin{corollary}
There is an algorithm that takes as input a quadruple \((\mathcal{C},\mathcal{D},\rho,\sigma)\), where \(\mathcal{C}\) is a groupoid, \(\mathcal{D}\subseteq \mathcal{C}\) is a wide subgroupoid and \(\sigma,\rho\in\wreath(\mathcal{C})\), and outputs all \(\alpha\in\wreath(\mathcal{D})\) such that \(\alpha\rho=\sigma\alpha\) in time \(n^{O(\#\obj(\mathcal{D})/\langle\rho\rangle)}\), where \(n\) is the length of the input. \qed
\end{corollary}

\newpage\section{Gradings of cyclotomic number fields}\label{sec:cyclo_ring}

In this section we will use the acquired theory to compute the gradings of cyclotomic fields \(\Q(\mu_{e})\) with \(e\in\Z_{>0}\) a prime power.
We simply list the gradings of \(\Q(\mu_e)\) with the torsion subgroup \(\mu_\infty\) of \(\C^*\) and later prove that these must be all such gradings.
The sections beyond this one do not depend on any of the results proven here and thus this section can be skipped by a reader that is only interested in the proof Theorem~\ref{AlgGridThm}.

\begin{theorem}\label{thm:cyclic_cyclotomic_gradings}
Let \(p\) be prime and \(k\geq 1\). Then the \(\mu_{\infty}\)-gradings of \(\Q(\mu_{p^k})\) are precisely the gradings in Example~\ref{2-grading_example}.
\end{theorem}

For now, we not only fix the ring \(\Q(\mu_{p^k})\) we grade, but also the group \(\mu_\infty\) we grade it with.
By nature of the definitions we use, if \(\overline{E}\) is a \(\mu_\infty\)-grading and \(\varphi\in\End(\mu_\infty)\), then we have \(\varphi_* \overline{E}=\overline{E}\) only when \(\varphi\) restricted to \(\langle \zeta\in\mu_{\infty} \,|\, E_\zeta \neq 0 \rangle\) is the identity, not simply when this restriction is an isomorphism. This is an important aspect to the counting argument we use in the proof of the theorem.

\begin{lemma}\label{lem:ugly_cite}
Let \(p\) be prime and \(k\geq 1\).
Then \(1\to 1+p\Z/p^k\Z\to(\Z/p^k\Z)^*\to(\Z/p\Z)^*\to 1\) is an exact sequence.
When \(p\neq 2\) the group \(1+p^n\Z/p^k\Z=\langle 1+p^n\rangle\) is a cyclic group of order \(p^{k-n}\) for \(1\leq n\leq k\).
If \(p=2\) and \(k\geq 2\), then  \(1+p\Z/p^k\Z=\langle-1\rangle\times(1+4\Z/2^k\Z)\), and \(1+2^n\Z/2^k\Z=\langle1+2^n\rangle\) is a cyclic group of order \(2^{k-n}\) for \(2\leq n\leq k\).
\begin{proof}
Standard. It follows from Proposition~2' in \cite{oude_algebra}.
\end{proof}
\end{lemma}

\begin{corollary}\label{cor:subgroup_diagram}
For each \(n\geq 1\) let \(\zeta_n\) be a primitive \(n\)-th root of unity and let \(\eta_{n}=\zeta_{n}+\zeta_{n}^{-1}\).
For \(k\geq 2\), the subgroups of \((\Z/2^k\Z)^*\) and subfields of \(\Q(\mu_{2^k})\) are as follows, where lines indicate subgroups of index \(2\) respectively field extensions of degree \(2\).
\begin{center}
\begin{tikzcd}[every arrow/.append style={dash}, every label/.append style = {font = \small},column sep={7em,between origins},row sep={3em,between origins}]
& & & & \langle 1 \rangle & \\
& & & \langle  1+2^{k-1} \rangle \arrow{ur} & \langle -1+2^{k-1} \rangle \arrow{u} &  \langle -1 \rangle \arrow{ul}\\ 
& & \langle1+2^4 \rangle \arrow[dotted]{ur} & & \langle -1,1+2^{k-1} \rangle \arrow{ul}\arrow{u}\arrow{ur} \\
& \langle 1+2^{3}\rangle \arrow{ur} & \langle -1 +2^{3} \rangle \arrow{u} & \langle -1,1+2^{4} \rangle \arrow{ul} \arrow[dotted]{ur} \\
\langle 1+2^2 \rangle \arrow{ur}& \langle -1+2^2 \rangle \arrow{u} & \langle -1,1+2^3 \rangle \arrow{ul} \arrow{u} \arrow{ur} & & \Q(\zeta_{2^k}) &\\
& \langle -1, 1+2^2 \rangle \arrow{ul}\arrow{u}\arrow{ur}& & \Q(\zeta_{2^{k-1}}) \arrow{ur} & \Q(\zeta_{2^{k-1}}\cdot\eta_{2^k}) \arrow{u} & \Q(\eta_{2^k})  \arrow{ul} \\
& & \Q(\zeta_{2^4}) \arrow[dotted]{ur}  & & \Q(\eta_{2^{k-1}})  \arrow{ul}\arrow{u}\arrow{ur}  \\
& \Q(\zeta_{2^{3}}) \arrow{ur} & \Q(\zeta_{2^3}\cdot\eta_{2^4}) \arrow{u} & \Q(\eta_{2^{4}})  \arrow{ul} \arrow[dotted]{ur} \\
\Q(\zeta_{2^2}) \arrow{ur} & \Q(\zeta_{2^2}\cdot\eta_{2^3}) \arrow{u} & \Q(\eta_{2^3}) \arrow{ul}\arrow{u}\arrow{ur} \\
& \Q \arrow{ul} \arrow{u} \arrow{ur}
\end{tikzcd}
\end{center}

\begin{proof}
The subgroups follow directly from Lemma~\ref{lem:ugly_cite}. For the fields we note that the extension \(\Q(\zeta_{2^k})/\Q\) is Galois, so we apply Galois correspondence and Corollary~\ref{cor:cyclotomic_automorphisms} to the diagram of groups.
\end{proof}
\end{corollary}

\begin{example}\label{2-grading_example}
Let \(E=\Q(\mu_{p^k})\) with \(p\) prime and \(k\geq 1\), let \(\zeta\in\mu_{p^k}\) be primitive and let \(G=\Z/p^{k-1}\Z\).

\textbf{a.} Note that \(\{\zeta^0,\dotsc,\zeta^{p^{k-1}-1}\}\) forms a \(\Q(\mu_p)\)-basis of \(E\) and that \(\zeta^{p^{k-1}}\in\Q(\mu_p)\). With \(C_n= \zeta^n \cdot \Q(\mu_p)\) for all \(n\in G\), we get a grading \(\overline{C}=(E,G,\{C_n\}_{n\in G})\). 
There are \(p^{k-1}\) morphisms \(G\to \mu_\infty\), giving \(p^{k-1}\) distinct \(\mu_\infty\)-gradings of \(E\), including the trivial grading.

\textbf{b.}
Assume \(p\neq 2\) or \(k\geq 3\).
Consider \(R_1=\Q(\alpha)\) and \(R_{-1}=\beta\cdot \Q(\alpha)\) with \(\alpha=\zeta+\zeta^{-1}\) and \(\beta=\zeta-\zeta^{-1}\).
By Corollary~\ref{cor:subgroup_diagram} we have \([E:R_1]=2\). Because \(\zeta=(\alpha+\beta)/2\in R_1+R_{-1}\), we must have \(E=R_1\oplus R_{-1}\).
Furthermore, we have \(\beta^2 = \zeta^2 + \zeta^{-2} - 2  = -4 + \alpha^2\in R_1\), hence \(R_{-1}R_{-1}\subseteq R_1\).
We conclude \(\overline{R}=(E,\mu_2,\{R_s\}_{s\in\mu_2})\) is a grading, which in turn gives a new \(\mu_\infty\)-grading of \(E\) under the natural inclusion \(\mu_2\to\mu_\infty\).
If \(p=2\), we can additionally decompose \(E\) into \(I_1=\Q(\beta)\) and \(I_{-1}=\alpha\cdot\Q(\beta)\) yielding another \(\mu_2\)-grading \(\overline{I}\), using that \(\Q(\beta)\) is the \(\langle-1+2^{k-1}\rangle\)-invariant subfield of \(E\).

\textbf{c.} Assume \(p=2\) and \(k\geq 4\).
Note that \(\sqrt{2}=\zeta_8+\zeta_8^{-1}\in E\) and consider \(\alpha=\zeta\sqrt{2}\in E\).
We have \(\Q(\sqrt{2})\subseteq\Q(\zeta_8)\subseteq\Q(\alpha^2)\), so \(\zeta=\alpha/\sqrt{2}\in\Q(\alpha)\) and \(E=\Q(\alpha)\).
Moreover, \(\alpha^{2^{k-1}}=-2^{2^{k-2}}\in\Q\). 
Thus with \(S_n=\alpha^n\cdot \Q\) we get a grading \(\overline{S}=(E,G,\{S_n\}_{n\in G})\).
Note that \(\alpha\) and \(\zeta^i\) are \(\Q\)-linearly independent for each \(i\in\Z\).
Hence for any of the \(2^{k-2}\) injections \(\varphi:G\to\mu_\infty\) the grading \(\varphi_*\overline{S}\) is not listed in (a), nor is it listed in (b).
For any non-injective \(\varphi\in\End(G)\) we have \(\varphi_*\overline{S}=\varphi_*\overline{C}\).
Thus we get exactly \(2^{k-2}\) new \(\mu_\infty\)-gradings of \(E\).

\textbf{d.} Assume \(p=2\) and \(k\geq 4\).
Consider \(\alpha=\zeta^2+\zeta^{-2}\) and \(\beta=\zeta_8(\zeta+\zeta^{-1})\).
By Corollary~\ref{cor:subgroup_diagram} the field \(\Q(\alpha)\) is the \(\langle-1,1+2^{k-1}\rangle\)-invariant subfield of \(E\) and we have \([E:\Q(\alpha)]= 4\).
Since \(\beta\) is neither invariant under \(\zeta\mapsto\zeta^{-1}\), \(\zeta\mapsto\zeta^{-1+2^{k-1}}\) nor \(\zeta\mapsto\zeta^{1+2^{k-1}}\), we have that \(\Q(\alpha,\beta)=E\).
Moreover \(\beta^4=(\zeta_4(2+\alpha))^2=-(2+\alpha)^2\in\Q(\alpha)\).
Then with \(W_n=\beta^n\cdot \Q(\alpha)\) we get a grading \(\overline{W}=(E,\Z/4\Z,\{W_n\}_{\Z/4\Z})\).
The injective morphisms \(\Z/4\Z\to \mu_\infty\) give two new \(\mu_\infty\)-gradings of \(E\).
Since \(\Q(\alpha,\beta^2)=\Q(\zeta^2)\) the others coincide with gradings of (a) by Corollary~\ref{cor:one_index_2_grad}.
\end{example}

First we show that besides a single possible exception, all \(\mu_\infty\)-gradings of \(\Q(\mu_{p^k})\) are \(\mu_{p^\infty}\)-gradings. 

\begin{lemma}\label{lem:cyclotomic_gradings_group_order}
Let \(p\) and \(q\) be prime and let \(E=\Q(\mu_{p^k})\) for some \(k\geq 1\).
If \(E\) has an efficient \(\mu_{q^n}\)-grading for some \(n\geq 1\), then \(q^n=2\) or \(p=q\).
Moreover, if \(p\neq 2\) then \(E\) has exactly one efficient \(\mu_2\)-grading and no efficient \(\mu_{2p}\)-gradings.
\begin{proof}
Let \(q\) be such that \(E\) has an efficient \(\mu_{q^n}\)-grading \(\overline{E}\) and assume \(p\neq q\).
Using Theorem~\ref{cyclic_correspondence}, this grading corresponds to some non-trivial \(\sigma\in\Aut_{\Z[\mu_{q^n}]}(E \tensor_\Z \Z[\mu_{q^n}] )\cong \Aut_{\Q(\mu_{q^n})}( \Q(\mu_{q^np^k}))\) such that \(\sigma^{q^n}=1\) and \(\tau_a \sigma=\sigma^a \tau_a\) for all \(a\in(\Z/q^n\Z)^*\), where \(\tau_a\) is the image of \(a\) under \((\Z/q^n\Z)^*\to\Aut(\Z[\mu_{q^n}])\to\Aut(E\tensor_\Z\Z[\mu_{q^n}])\). 
Note that \(\tau_a\) and \(\sigma\) commute as \(\Aut(\Q(\mu_{q^np^k}))\) is abelian. 
Then for all \(a\in(\Z/q^n\Z)^*\) we have \(\sigma^a = \tau_a \sigma \tau_a^{-1} = \tau_a \tau_a^{-1} \sigma = \sigma\), so \(a=1\).
Thus \((\Z/q^n\Z)^*=\{1\}\) and \(q^n=2\).
It follows that \(E\) can only have a non-trivial efficient \(\mu_{q^n}\)-grading if \(q^n=2\) or \(p=q\).

Assume \(p\neq 2\). Then \(\Aut_{\Z[\mu_2]}(E\tensor_\Z \Z[\mu_2]) \cong \Aut(E) \cong (\Z/p^k\Z)^*\) is cyclic of even order by Lemma~\ref{lem:ugly_cite}.
Hence it contains exactly one \(\sigma\) of order \(2\), and it trivially satisfies \(\tau_a\sigma=\sigma^a\tau_a\) for all \(a\in (\Z/2\Z)^*\).
The corresponding efficient \(\mu_2\)-grading of \(E\) is then unique.
If \(\overline{E}\) is an efficient \(\mu_{2p}\)-grading of \(E\), then we get an efficient \(\mu_p\)-grading of the unique subfield \(E'=\sum_{\zeta\in\mu_p} E_\zeta\) of \(E\) of degree two. 
In particular, we may write \(E'=F(\alpha)\) for some \(\alpha\in E'\) and subfield \(F\subseteq E'\) of degree \(p\) such that \(\alpha^p\in F\). 
The minimal polynomial of \(\alpha\) must then be \(f=X^p-\alpha^p\), but \(\zeta_p\alpha\) is a root of \(f\) while \(\zeta_p=\zeta_p\alpha/\alpha\not\in E'\), so \(E'/F\) is not normal.
However, \(E/\Q\) is a Galois extension, so this is a contradiction and no efficient \(2p\)-grading exists.
\end{proof}
\end{lemma}

Lemma~\ref{lem:cyclotomic_gradings_group_order} shows that, besides a single \(\mu_2\)-grading in the case \(p\neq 2\), all \(\mu_\infty\)-gradings of \(\Q(\mu_{p^k})\) are \(\mu_{p^\infty}\)-gradings.
Since the order of \(\langle \gamma\in\Gamma \,|\, E_\gamma\neq 0\rangle \) divides \([\Q(\mu_{p^k}):\Q]=(p-1)p^{k-1}\) for any \(\Gamma\)-grading \(\overline{E}\) of \(\Q(\mu_{p^k})\) by Proposition~\ref{prop:basic_grad_props}, we may even restrict to the group \(\mu_{p^k}\).

%

\begin{lemma}\label{reduction_to_subgroups}
Let \(p^k>1\) be a prime power and let \(A=\Z/p^k\Z\) and \(M=(\Z/p^k\Z)^*\).
Then \(M\) acts on \(A\) by multiplication and we let \(G=A\rtimes M\) be the corresponding semi-direct product.
Then there is a bijection
\[ \{\mu_{p^k}\text{-gradings of }\Q(\mu_{p^k})\} \leftrightarrow Z = \{ H \subseteq G \,|\, H\text{ a group},\, \text{the natural map }M\to G/H\text{ is a bijection} \}.\]
\begin{proof}
By Theorem~\ref{cyclic_correspondence}, the set of \(\mu_{p^k}\)-gradings of \(E=\Q(\mu_{p^k})\) is in bijection with \(X_{p^k}(E)\) as defined in Definition~\ref{def:X}.
Then by Proposition~\ref{aut_wreath_iso}, we have \(\Aut_{\Z[\mu_{p^k}]\catt{-Alg}}(E')\cong \wreath(\mathcal{C})\) with \(E'=E\tensor_\Z\Z[\mu_{p^k}]\) and \(\mathcal{C}\) the category with \(\obj(\mathcal{C})=\spec E'\) and \(\Hom_{\mathcal{C}}(\m,\n)=\Hom_{\Z[\mu_{p^k}]\catt{-Alg}}(E'/\m,E'/\n)\) for all \(\m,\n\in\obj(\mathcal{C})\).
We have an injection \(\Q(\mu_{p^k})\hookrightarrow E'/\m\) since \(E'/\m\) is a \(\Q(\mu_{p^k})\)-algebra and morphisms of fields are injective and an injection of \(E'/\m\) into a compositum of the factors of \(E'\), which is \(\Q(\mu_{p^k})\).
Therefore \(E'/\m\cong\Q(\mu_{p^k})\) for all \(\m\in\spec E'\), so \(\wreath(\mathcal{C})\cong \Aut_\catt{Set}(\obj(\mathcal{C}))\) by Proposition~\ref{prop:wreath_is_semi-direct_product}.
Moreover, \(M\cong\Aut(\Q(\mu_{p^k}))\) acts regularly on \(\obj(\mathcal{C})\), so we have a canonical isomorphism \(\Aut_\catt{Set}(\obj(\mathcal{C}))\cong \Aut_\catt{Set}(M)\) such that the induced action \(C_M:M\to\Aut_\catt{Set}(M)\) is the Cayley-action.
Then \(X_{p^k}(\Q(\mu_{p^k}))\) is in natural bijection with \(X=\{ \gamma\in\Aut_\catt{Set}(M) \,|\, \gamma^{p^k} = 1,\ (\forall a\in M)\, C_M(a)\cdot \gamma \cdot C_M(a^{-1}) = \gamma^a \}\).
With \(Y=\{\varphi\in\Hom(G,\Aut_\catt{Set}(M))\,|\, \varphi|_M=C_M\}\) we have a bijection \(X\to Y\) given by \(\gamma\mapsto [ \sigma^n\tau_a \mapsto \gamma^n C_M(a) ]\), where \(\sigma\) and \(\tau_a\) are the image of \(1\in A\) respectively \(a\in M\) in \(G\).

Now consider the map \(\Phi:Y\to Z\) given by \(\varphi\mapsto\{ g\in G\,|\,\varphi(g)(1)=1 \}\).
To show that it is well-defined, consider \(\varphi\in Y\) and let \(H=\Phi(\varphi)\), which is clearly a group.
If \(\tau_a\in H\) then \(1=\varphi(\tau_a)(1)=a\) hence \(\tau_a=1\), so \(M\to G/H\) is injective.
For each \(\sigma^n\in A\) we may take \(a=\varphi(\sigma^n)(1)\) such that \(\varphi(\tau_a^{-1}\sigma^n)(1)=1\) hence \(\tau_a^{-1}\sigma^n\in H\).
Thus \(\# A \leq \# H\) and \(\# G/H \leq \# M\), so \(M\to G/H\) is also surjective.
Hence \(H\in Z\) and \(\Phi\) is well-defined.

For \(H\in Z\) we have a natural action \(G\to\Aut_\catt{Set}(G/H)\) and using the natural bijection \(G/H\to M\) an action \(\varphi_H:G\to\Aut_\catt{Set}(M)\). Note that \(\varphi_H|_M=C_M\) and \(\Phi(\varphi_H)=H\).
Hence the map \(\Psi:Z\to Y\) given by \(H\mapsto \varphi_H\) is a right inverse of \(\Phi\).
For \(\varphi\in Y\), \(g\in G\) and \(a\in M\) we may uniquely write \(g\tau_a=\tau_bh\) for some \(h\in\Phi(\varphi)\) and \(b\in M\), hence \((\Psi\circ\Phi)(\varphi)(g)(a)=(\Psi\circ\Phi)(\tau_bh)(1)=b=\varphi(\tau_b h)(1)=\varphi(g)(a)\), hence \(\Psi\) is a left inverse of \(\Phi\).
Thus \(\Phi\) is bijective and the lemma follows.
\end{proof}
\end{lemma}

\begin{lemma}\label{lem:cyclo_help}
Let \(p\) prime, \(0 < n < k\) and consider \(G=A\rtimes M\) as in Lemma~\ref{reduction_to_subgroups}. 
Assume \(\alpha,\beta\in A\) and \(\langle\tau\rangle=1+p^{k-n}\Z/p^k\Z\subseteq M\).
\begin{enumerate}[topsep=0pt,itemsep=-1ex,partopsep=1ex,parsep=1ex]
\item If \(\alpha^{p^{n-1}}\in\langle\beta\rangle\), then \(\langle \beta, \alpha\tau\rangle \cap M \neq 1\).
\item If \(p=2\), \(\beta^2\neq 1\) and \(2\leq n\), then \(\langle \beta, \alpha\eta\tau \rangle \cap M \neq 1\) with \(\eta=-1\in M\).
\end{enumerate}
\begin{proof}
(1) Consider \(S_r = \sum_{i=0}^{r-1} \tau^i\) and note that \((\alpha\tau)^r=\alpha^{S_r}\tau^r\) for all \(r\geq 0\). 
We have
\[S_{p^{n-1}} = \sum_{i=0}^{p^{n-1}-1} \tau^i \equiv \sum_{i=0}^{p^{n-1}-1}(1+ip^{k-n}) \equiv p^{n-1}\bigg( 1+\frac{p^{k-n}(p^{n-1}-1)}{2} \bigg)\mod p^{k-1},\]
using that both sums range over all elements of \(1+p^{k-n}\Z/p^{k-1}\Z\) and \(n<k\). 
Hence \(p^{n-1}\divs S_{p^{n-1}}\).
When \(\alpha^{p^{n-1}}\in\langle\beta\rangle\) we have \(\alpha^{S_{p^{n-1}}} = \beta^r \) for some \(r\), so \(1\neq\tau^{p^{n-1}}=\beta^{-r}(\alpha\tau)^{p^{n-1}}\in\langle \beta,\alpha\tau\rangle\), as was to be shown.

(2) Assume \(p=2\), \(\beta\neq 1\) and \(2\leq n\). 
Then \((\alpha\eta\tau)^2=\alpha^{1-\tau}\tau^2\) and note that \(1\neq\tau^2\) generates \(1+2^{k-n+1}\Z/2^k\Z\) since \(2\leq n\).
We have \((\alpha^{1-\tau})^{2^{n-2}}\in 2^{k-2}\Z/2^k\Z\subseteq\langle\beta\rangle\) because \(\beta^2\neq 1\).
The claim now follows from (1) applied to \(\alpha^{1-\tau}\) and \(\tau^2\).
\end{proof}
\end{lemma}

\noindent\textbf{Proof of Theorem~\ref{thm:cyclic_cyclotomic_gradings}.}
Using Lemma~\ref{lem:cyclotomic_gradings_group_order}, with the same notation as in Lemma~\ref{reduction_to_subgroups}, it suffices to bound \(|Z|\) from above by the number of \(\mu_{p^k}\)-gradings of \(\Q(\mu_{p^k})\) listed in Example~\ref{2-grading_example}.
The group \(G\) induces an exact sequence \(1\to A \to G\to M \to 1\), and any subgroup \(H\subseteq G\) induces an exact sequence \(1\to H_A \to H\to H_M\to 1\) where \(H_A=H\cap A\) and \(H_M=\pi_M(H)\) is the projection to \(M\). Since \(A=\langle \sigma\rangle\) is cyclic of order \(p^k\), \(H_A=\langle\sigma^{p^n}\rangle\) is uniquely determined by its size \(p^{k-n}\). 
If \(H\in Z\), then by \(p^k=|H|=|H_A|\cdot|H_M|=p^{k-n} |H_M|\) we have \(|H_M|=p^n\). In particular, \(H_M\subseteq 1+p\Z/p^k\Z\) and \(n<k\) using Lemma~\ref{lem:ugly_cite}.
By the same lemma, if \(p\neq 2\), then \(H_M=1+p^{k-n}\Z/p^k\Z=\langle 1+p^{k-n} \rangle\), which is cyclic.
From Corollary~\ref{cor:subgroup_diagram} it follows that, in general, we may distinguish between the following mutually exclusive cases.
\begin{align*}
\begin{array}{lllr}
\text{(1)} & H_M = \langle 1+p^{k-n} \rangle &\quad \text{for } 0\leq n \leq k-1 \text{ and if \(p=2\) also } n \leq k-2 \\
\text{(2)} & H_M = \langle -1+p^{k-n} \rangle & \quad\text{for } p = 2 \text{ and } 1 \leq n \leq k-2 \\
\text{(3)} & H_M = \langle 1+p^{k-n+1}, -1 \rangle & \quad\text{for } p = 2 \text{ and } 1 \leq n \leq k-1 
\end{array}
\end{align*}

{\em Case }(1):
Here \(H=\langle \sigma^{p^n}, \sigma^r \tau \rangle\) for some \(r\in\Z/p^n\Z\), where \(\tau\) is the image of \(1+p^{k-n}\in M\) in \(G\).
If \(n=0\), then \(H=A\), which is a single solution.
Now assume \(n>0\).
If \(p\divs r\), then \(H\cap M \neq 1\) by Lemma~\ref{lem:cyclo_help}.1.
Thus \(r\in(\Z/p^n\Z)^*\), meaning that given \(n\) there are at most \(p^n-p^{n-1}\) possible groups \(H\) in this case.

{\em Case }(2): Here \(H=\langle \sigma^{2^n}, \sigma^r \eta\tau\rangle\) for some \(r\in\Z/2^n\Z\), \(\eta=-1\) and \(\tau=1+2^{k-n}\).
If \(n=1\), then \(r=1\) as \(\eta\tau\not\in H\), which yields one possible \(H\).
For larger \(n\) we have \(\sigma^{2^{n+1}}\neq 1\) hence there are no possible subgroups \(H\) by Lemma~\ref{lem:cyclo_help}.2.

{\em Case }(3): Here \(H=\langle \sigma^{2^n}, \sigma^r \tau, \sigma^s \eta \rangle\) for \(r,s\in\Z/2^n\Z\), \(\eta=-1\) and \(\tau=1+2^{k-n+1}\). 
If \(n=1\), then \(r=0\) and \(s=1\), giving a single subgroup.
If \(n=2\), then \(r,s\neq 0\) and \(r\neq s\) as \(\sigma^{s-t}\eta\tau\in H\) while \(1\neq \eta\tau \in M\). 
Moreover \((\sigma^r\tau)^2=\sigma^{r(2+2^{k-1})}\in H\) implies \(2\divs r\). 
Thus \((r,s)\in\{(2,1),(2,3)\}\) gives two possible subgroups.
If \(3\leq n<k-1\), we consider \(\sigma^{s-r}\eta\tau=\sigma^s\eta \sigma^r\tau\in H\) and note that by Lemma~\ref{lem:cyclo_help}.2 no subgroups exist.
Hence assume \(n=k-1\). Then \(\sigma^{2r+4s}=[\sigma^r\tau,\sigma^s\eta]\in H\), meaning that \(2^{k-2}\divs r+2s\).
We have that \(4\nmid r\) by applying Lemma~\ref{lem:cyclo_help} to \(\sigma^r\tau\), so \(2\nmid s\).
There are then \(2^{k-2}\) possible values of \(s\), and for each only two possible \(r\), giving at most \(2^{k-1}\) subgroups \(H\).

If \(p\neq 2\), then summing over all \(n\) in case (1) gives \(|Z|\leq 1+\sum_{i=1}^{k-1} (p^i-p^{i-1}) = p^{k-1}\), the right hand side of which equals the number of \(\mu_{p^k}\)-gradings of \(\Q(\mu_{p^k})\) listed in Example~\ref{2-grading_example}.
If \(p=2\), then \(|Z|\leq 2\) if \(k=2\), \(|Z|\leq 6\) if \(k=3\) and \(|Z|\leq 4+2^{k-2}+2^{k-1}\) if \(4\leq k\), which is again precisely the number of distinct \(\mu_{2^k}\)-gradings of \(\Q(\mu_{2^k})\) listed in Example~\ref{2-grading_example}, as was to be shown. \qed

\clearpage
\section{Existence of universal gradings}

In this section we present two proofs for the existence of universal gradings of reduced orders.

\subsection{Primitive idempotents and factorization of gradings}

\begin{definition}\label{def:idempotent}
Let \(R\) be a commutative ring. 
Then \(e\in R\) is called {\em idempotent} if \(e^2=e\).
We call \(e\) a {\em primitive idempotent} if \(e\neq 0\) and \(ee'\in\{e,0\}\) for all idempotents \(e'\). 
We write \(\text{prid}(R)\) for the set of primitive idempotents of \(R\).
We say \(R\) is {\em connected} when \(1\in \prid R\).
\end{definition}

Note that for each commutative ring \(R\) and idempotent \(e\in R\) we get a ring \(eR\).

\begin{lemma}\label{lem:idempotent_isomorphism}
Let \(R\) be a commutative ring. 
If \(\#\minspec R<\infty\), then \(R\) has only finitely many idempotents.
If \(R\) has only finitely many idempotents, then \(\sum_{e\in\prid R} e = 1\) and the natural morphism of abelian groups \(\prod_{e\in\prid R} eR \to R\) is an isomorphism of rings.
\begin{proof}
The natural map \(R/\sqrt{0}\to\prod_{\p\in\minspec R} R/\p\) is injective by Lemma~\ref{basic_comalg_wreath}.
The idempotents of \(\prod_{\p} R/\p\) are \(\{0,1\}^{\minspec R}\) since each \(R/\p\) is a domain.
If \(\#\minspec R<\infty\), then \(\prod_\p R/\p\) and therefore \(R/\sqrt{0}\) have only finitely many idempotents.
Assume \(e,f\in R\) are idempotents such that \(e-f\in\sqrt{0}\).
Then for some \(n\in\Z_{\geq 0}\) we have 
\[0=(e-f)^{2n+1}=\sum_{i=0}^{2n+1}\binom{2n+1}{i}e^i(-f)^{2n+1-i} =  e-f + ef\sum_{i=1}^{2n}\binom{2n}{i}(-1)^i = e-f, \]
hence \(e=f\). Hence \(R\) has only finitely many idempotents as well. 

Assume \(R\) has only finitely many idempotents.
Note that the primitive idempotents are simply the minimal elements among the non-zero idempotents of \(R\) under the partial order \(\preceq\) given by \(e\preceq e'\) if and only if \(ee'=e\).
Now let \(x=\sum_{e\in\prid R} e\) and assume \(x\neq 1\). Then \(y=1-x\neq 0\) is an idempotent, so there exists some \(e'\in\prid R\) such that \(e'\preceq y\) by finiteness.
Then \(e'=e'y=e'-\sum_{e\in\prid R} ee'=e'-e'=0\) since \(ee'=0\) for all \(e'\neq e\in\prid R\), which is a contradiction.
Hence \(x=1\).

Now consider the natural map \(\varphi:\prod_{e\in\prid R} eR\to R\) and the map \(\psi:R\to \prod_{e\in\prid R} eR\) given by \(r\mapsto (er)_{e\in\prid R}\).
Clearly \(\psi\varphi=\id\) since \(ee'=0\) for distinct \(e,e'\in\prid R\), while \(\varphi\psi=\id\) additionally requires that \(\sum_{e\in\prid R} e= 1\).
Hence \(\varphi\) is an isomorphism.
\end{proof}
\end{lemma}

\begin{corollary}\label{cor:spectrum_and_idempotents}
Let \(R\) be a commutative reduced \(\Q\)-algebra of finite dimension as \(\Q\)-module. 
Then \(R\) is connected if and only if it is a field, and \(\spec R = \{ (1-e)R \,|\, e\in\prid(R) \}\).
\begin{proof}
Use the facts that \(eR\cong R/(1-e)R\) for idempotents \(e\), that \(R\cong \prod_{\m\in\spec R}R/\m\) and that \(\prid \prod_{\m} R/\m=\{e_\n\,|\, \n\in\spec R\}\) with \(e_\n=(\mathbbm{1}_{\m=\n})_{\m\in\spec R}\).
\end{proof}
\end{corollary}

In Lemma~\ref{lem:idempotent_isomorphism} we used the primitive idempotents to factor our rings.
We can apply the same technique to gradings.
For this, recall that we have defined a coproduct of gradings in Definition~\ref{def:coproduct_grading}.

\begin{definition}\label{def:grading_factorization}
Let \(R\) be a commutative ring with only finitely many idempotents and let \(\overline{R}=(R,\Gamma,\{R_\gamma\}_{\gamma\in\Gamma})\) be a grid-grading.
For each \(e\in\prid R_1\) we define \(\Gamma_e=\langle \gamma\in\Gamma\,|\, eR_\gamma\neq 0\rangle\) and \(\overline{eR}=(eR,\Gamma_e,\{eR_\gamma\}_{\gamma\in\Gamma_e})\).
\end{definition}

\begin{lemma}\label{lem:coproduct_idempotent_intermediate}
Let \(R\) be a commutative ring with only finitely many idempotents and let \(\overline{R}=(R,\Gamma,\{R_\gamma\}_{\gamma\in\Gamma})\) be a grid-grading.
Then \(\overline{eR}\) and \(\overline{R'}:=\coprod_{e\in\prid R_1} \overline{eR}\) are efficient gradings of \(R\) and the natural map \(\varphi:\coprod_{e\in\prid R_1} \Gamma_e\to\Gamma\) satisfies \(\varphi_* \overline{R'} = \overline{R}\).
\begin{proof}
Clearly \(1_{eR}=e\in eR_1\) and \(eR_\gamma \cdot e R_{\gamma'}=e R_\gamma R_{\gamma'}\subseteq e R_{\gamma\gamma'}\), so \(\overline{eR}\) is a grading. 
Then \(\overline{R'}\) is a grading of \(R\) by Lemma~\ref{lem:idempotent_isomorphism}, and \(\varphi_* \overline{R'}=\overline{R}\) follows trivially.
Additionally \(\overline{eR}\) and \(\overline{R'}\) are efficient by construction.
\end{proof}
\end{lemma}

Lemma~\ref{lem:coproduct_idempotent_intermediate} will allow us to assume that the trivial components of our gradings are connected in some of the proofs, most notably Proposition~\ref{prop:orthogonal_homogenous}.
In the case of finite-dimensional reduced \(\Q\)-algebras we even have that each factor \(\overline{eR}\) of \(\overline{R}\) is an abelian group-grading by Corollary~\ref{cor:spectrum_and_idempotents} and Proposition~\ref{prop:basic_grad_props}.

\begin{example}\label{ex:non-compact_counter}
In this example we construct a ring \(R\) for which \(X_\infty(R)\) as defined in Definition~\ref{def:X} is an abelian group but is not compact, meaning that any pair of abelian group-gradings of \(R\) have a joint grading as follows from Theorem~\ref{thm:grad_is_morph}, while no universal grading exists by Theorem~\ref{thm:commutative_implies_universal}.
Consider \(R=\Z[\sqrt{2}]^\Z\) and let \(p\) be prime.
Note that \(\sqrt{2}\in\Q(\zeta_8) \setminus \Q(\zeta_4)\) hence \(\sqrt{2}\not\in\Q(\zeta_{p^2})\).
Thus \(\Z[\sqrt{2}]\tensor_\Z\Z[\zeta_{p^2}]\cong\Z[\sqrt{2},\zeta_{p^2}]\) is a domain.
Now consider \(R'=\Z[\sqrt{2},\zeta_{p^2}]^\Z\) and note that we have a natural map \(R\tensor_\Z\Z[\zeta_{p^2}]\to R'\).
Since \(\Z[\zeta_{p^2}]\) is free of finite rank as \(\Z\)-module, we even have that this map is an isomorphism.
We have that each automorphism of \(R'\) must permute the primitive idempotents, which are the elements \(e_i=(\mathbbm{1}(j=i))_{j\in\Z}\) of \(R'\).
Thus we have a natural map \(\Aut_{\Z[\zeta_{p^2}]\catt{-Alg}}(R')\to \Aut_{\Z[\zeta_{p^2}]\catt{-Alg}}(\Z[\sqrt{2},\zeta_{p^2}])^\Z\rtimes\Aut_\catt{Set}(\Z)\).
Using the universal property of the product one can show that this map is an isomorphism as well.

Now let \(\sigma\in X_{p^2}(R)\) and \(i\in\Z\).
Let \(n=\#(\langle\sigma\rangle e_i)\) and assume \(n\neq 1\).
Because \(n\divs \ord(\sigma)\divs p^2\) we have \(\zeta_n\in\Z[\zeta_{p^2}]\).
For \(0\leq m<n\) let \(a_m=\sum_{j=0}^{n-1}\zeta_n^{-mj}\sigma^j(e_i)\) and note that \(a_m\in R'(\sigma,\zeta_n^{m})\).
Then
\[\sum_{m=0}^{n-1} a_m = \sum_{j=0}^{n-1}\sigma^j(e_i)\sum_{m=0}^{n-1}\zeta_n^{-mj} = \sum_{j=0}^{n-1}\sigma^j(e_i) \cdot n\mathbbm{1}(j=0) = ne_i.\]
We have \(a_m/n\not\in R'\) since \(\zeta_n/n\not\in\Z[\zeta_n]\), but \(e_i\in\bigoplus_{j=0}^{p^2-1} R'(\sigma,\zeta_{p^2}^j)\) by \(\mu_{p^2}\)-diagonalizability of \(\sigma\), a contradiction.
Thus \(n=1\).
It follows that \(X_{p^2}(R)\subseteq\Aut_{\Z[\zeta_{p^2}]\catt{-Alg}}(\Z[\sqrt{2},\zeta_{p^2}])^\Z\).
Since the latter has exponent \(2\), we only obtain non-trivial gradings when \(p=2\), and the only efficient gradings of \(R\) with \(\mu_e\) for \(e\in\mathbb{S}\) have \(e\divs 2\).
We conclude that \(X_\infty(R)\cong X_2(R)\cong\Aut(\Z[\sqrt{2}])^\Z\), which is abelian.
However, it is not compact since \((\sqrt{2})_{i\in\Z}\) has an infinite orbit under \(X_\infty(R)\).
Hence for any two abelian torsion group-gradings of \(R\) a joint-grading exists, while no universal abelian torsion group-grading of \(R\) exists.
\end{example}

\subsection{Hierarchy of categories}

We consider gradings of a ring \(R\) with a grid \(G\).
In this section we show that if \(R\) has a universal grid-grading, then it has a universal group-grading.
Similarly the existence of a universal abelian group-grading follows from the existence of a universal group-grading.

\begin{proposition}\label{prop:hierarchy}
Let \(k\) be a commutative ring, let \(R\) be a \(k\)-algebra and let \(\overline{U}=(R,\Upsilon,\mathcal{R})\) be a grading.
If \(\overline{U}\) is a universal grid-grading, then \(f_*\overline{U}\) is a universal group-grading, with \(f:\Upsilon\to \Upsilon^\catt{grp}\) the groupification map as in Lemma~\ref{lem:groupification_adjunction}.
If \(\overline{U}\) is a universal group-grading, then \(f_*\overline{U}\) is a universal group-grading, with \(f:\Upsilon\to \Upsilon^\catt{ab}\) the abelianization map.
\begin{proof}
We have forgetful functors \(M_\catt{Grp}:\catt{Grp}\to\catt{Grd}\) and \(M_\catt{Ab}:\catt{Ab}\to\catt{Grp}\) which are right adjoints of \(\us^\catt{grp}\) (Lemma~\ref{lem:groupification_adjunction}) respectively \(\us^\catt{ab}\).
For \(\mathcal{C}\in\{\catt{Grd},\catt{Grp},\catt{Ab}\}\), there exists a universal \(\mathcal{C}\)-grading of \(R\) if and only if the functor \(F_\mathcal{C}:\mathcal{C}\to\catt{Set}\) given by \(\Gamma\mapsto\{\Gamma\text{-gradings of }R\}\) is representable.
Note that \(F_\catt{Ab}\cong F_\catt{Grp}\circ M_\catt{Ab}\) and \(F_\catt{Grp} \cong F_\catt{Grd}\circ M_\catt{Grp}\).
If we have categories \(\mathcal{C}\) and \(\mathcal{D}\), a right adjoint functor \(M:\mathcal{D}\to\mathcal{C}\) and a representable functor \(F:\mathcal{C}\to\catt{Set}\), then \(F\circ M\) is representable.
Namely, let \(L\) be a left adjoint of \(M\) and assume \(F\) is representable.
Then there exists an \(A\in\obj(\mathcal{C})\) and natural isomorphism \(\Hom_\mathcal{C}(A,\us)\cong F\).
Hence \(F\circ M \cong \Hom_\mathcal{D}(A,M(\us)) \cong \Hom_\mathcal{C}(L(A),\us)\) so \(F\circ R\) is representable.
Moreover, we obtain a universal element of \(F\circ M\) from a universal element of \(F\) by applying \(L\) to it, from which the lemma follows.
\end{proof}
\end{proposition}

\begin{proposition}\label{prop:refinement_universal}
Let \(k\) be a commutative ring and \(R\) be a \(k\)-algebra with a decomposition \(\mathcal{M}\) not containing \(\{0\}\).
Then the functor \(F:\catt{Grd}\to\catt{Set}\) given by 
\[\Gamma\mapsto \{ \text{grid-gradings } (R,\Gamma,\mathcal{R}) \text{ such that } \mathcal{M}\preceq\mathcal{R} \}\]
is representable, where \(\preceq\) is as defined in Definition~\ref{def:decomposition}.
\begin{proof}
Write \(\mathcal{M}=\{M_i\}_{i\in I}\) and let \(X\) be the set of equivalence relations \(\sim\) on \(I\) such that \(\{N_j\}_{j\in I/{\sim}}:=\pi^\sim_{*} \mathcal{M}\), where \(\pi^\sim:I\to I/{\sim}\) is the quotient map, is a pre-grading of \(R\) such that \(1\in N_e\) for some unique \(e\in I/{\sim}\) and if \(N_j N_j \subseteq N_j\) for \(j\in I/{\sim}\) we have \(j=e\).
Now define an equivalence relation \(\simeq\) on \(I\) such that for \(i,j\in I\) we have \(i\simeq j\) if and only if \(i\sim j\) for all \({\sim}\in X\).
It is a straightforward verification that \({\simeq}\in X\). 
Let \(S=I/{\simeq}\) and \(\mathcal{U}=\{U_s\}_{s\in S}:=\pi^{\simeq}_*\mathcal{M}\).

If \(R=0\) then \(I=\emptyset\) and \(F\cong \Hom_{\catt{Grd}}(\{1\},\us)\) is representable.
Hence assume \(R\neq 0\).
Let \(e\in S\) be the unique element such that \(1\in U_e\). 
Let \(\Upsilon=(S,e,D,*)\) with \(D=\{ (s,t)\in S \,|\, U_s\cdot U_t \neq 0 \}\) and \(s*t\) the unique \(u\in S\) such that \(U_s\cdot U_t \subseteq U_u\) for all \((s,t)\in D\).
Note that \((e,s),(s,e)\in D\) and \(s=s*1=1*s\) for all \(s\in S\) since \(0\neq 1\in U_e\) and \(U_s\neq 0\) so \(U_s\cap (U_e U_s) \cap (U_s U_e) \neq 0\). 
Then \(\Upsilon\) is a grid and \(\overline{U}=(R,\Upsilon,\mathcal{U})\) is a loose grid-grading.

Now let \(\overline{R}=(R,\Gamma,\mathcal{R})\) be any grid-grading such that \(\mathcal{M}\preceq\mathcal{R}=:\{R_\gamma\}_{\gamma\in\Gamma}\).
Then there exists a map \(f:I\to\Gamma\) such that \(f_*\mathcal{M}=\mathcal{R}\), and \(f\) induces an equivalence relation \(\sim\) on \(I\) where \(i\sim j\) if and only if \(f(i)=f(j)\).
Note that \({\sim}\in X\), hence \(\pi^{\sim}=q\circ \pi^{\simeq}\) for some unique \(q\) by definition of \(\simeq\).
Then \(f=g\circ \pi^{\simeq}\) for some unique \(g:S\to \Gamma\).
We note that \(g\) as map of grids \(G\to \Gamma\) is in fact a morphism, hence \(\overline{R}=g_*\overline{U}\).
Thus \(\overline{U}\) is a universal object of \(F\) and \(F\) is representable.
\end{proof}
\end{proposition}

\subsection{A proof using lattices}\label{sec:UniGridGradThm}

We prove the existence of universal grid gradings using a generalization of the proof in \cite{Lenstra2018}, Section 9.

\begin{definition}\label{hom_inner_product}
Let \(R\) be a subring of a finite product of number fields. We then define the map \(R^2\to\C\) given by
\[(x,y)\mapsto\langle x, y\rangle_R := \sum_{\sigma\in\Hom_\catt{Rng}(R,\C)} \sigma(x)\overline{\sigma(y)}. \]
\end{definition}

\begin{remark}
Note that \(\langle\us,\us\rangle_R\) as in Definition~\ref{hom_inner_product} is in fact real valued.
For \(\sigma\in\Hom_{\catt{Rng}}(R,\C)\) the map \(x\mapsto\overline{\sigma(x)}\) is also an element of \(\Hom_{\catt{Rng}}(R,\C)\), so complex conjugation of \(\langle x,y\rangle_R\) for \(x,y\in R\) simply permutes the terms of the sum defining \(\langle\us,\us\rangle_R\).
Hence \(\overline{\langle x,y\rangle}_R=\langle x,y\rangle_R\) must be real.
Additionally, for all \(x,y\in R\) we have the following properties: (1) \(\langle\us,y\rangle_R\) is \(\Z\)-linear; (2) \(\langle x,y\rangle_R=\langle y,x\rangle_R\) and (3) \(\langle x,x\rangle_R\in\R_{\geq0}\) with \(\langle x,x\rangle_R=0\) if and only if \(x=0\).
If \(R\) is a \(\Q\)-module, we would call \(\langle\us,\us\rangle_R\) an inner product. 
\end{remark}

\begin{lemma}\label{lem:hom_to_C}
Let \(R\) be a finite product of number fields. Then for all \(n\in\Z_{\geq0}\), \(x,y\in R\) and \(e\in R\) idempotent we have \(\langle ex,ey \rangle_R=\langle ex,ey\rangle_{eR}\) and \(\langle x,y\rangle_{R\tensor_\Q\Q(\mu_n)}=\varphi(n)\cdot\langle x,y\rangle_R\), where \(\varphi\) is Euler's totient function.
\begin{proof}
For each idempotent \(e\in R\) a map \(eR\to\C\) extends to a map \(R\to\C\) by composing with the projection.
We have an isomorphism \(R\cong eR\times (1-e)R\) and from this it follows the natural map \(\Hom_\catt{Rng}(eR,\C) \sqcup \Hom_{\catt{Rng}}((1-e)R,\C)\to\Hom_\catt{Rng}(R,\C)\) is a bijection.
Then \(\langle ex,ey \rangle_R=\langle ex,ey\rangle_{eR}\) for all \(x,y\in R\) follows trivially.

Let \(n\in\Z_{\geq0}\). By the universal property of the tensor product we have mutually inverse maps 
\begin{align*}
\Hom_{\Q\catt{-Alg}}(R,\C)\times\Hom_{\Q\catt{-Alg}}(\Q(\mu_n),\C)&\to\Hom_{\Q\catt{-Alg}}(R\tensor_\Q\Q(\mu_n),\C) \phantom{XXXXXX}\\
(f,g) &\mapsto f\tensor g = \big[ x\tensor y \mapsto f(x)\tensor g(y) \big] \\
\big[ x \mapsto f(x\tensor 1), \ y \mapsto f(1\tensor y) \big]   &\mapsfrom f
\end{align*}
From Corollary~\ref{cor:cyclotomic_automorphisms} it follows that \(\#\Hom_{\Q\catt{-Alg}}(\Q(\mu_n),\C)=\varphi(n)\).
Let \(\sigma\in\Hom_\catt{Rng}(R,\C)\).
Now \((\sigma\tensor\chi)(x)=\sigma(x)\) for all \(\chi\in\Hom_\catt{Rng}(\Q(\mu_n),\C)\) and \(x\in R\).
Hence 
\[\sum_{\chi\in\Hom_\catt{Rng}(\Q(\mu_n),\C)} (\sigma\tensor\chi)(x)\overline{(\sigma\tensor\chi)(y)}=\varphi(n) \cdot \sigma(x)\overline{\sigma(y)}\]
for all \(x,y\in R\), from which it follows that \(\langle x,y\rangle_{R\tensor_\Q\Q(\mu_n)}=\varphi(n)\cdot\langle x,y\rangle_R\).
\end{proof}
\end{lemma}

\begin{proposition}\label{prop:orthogonal_homogenous}
Let \(R\) be a product of number fields with \(\dim_\Q(R)<\infty\), let \(\Gamma\) be a grid and let \(\overline{R}\) be a \(\Gamma\)-grading of \(R\). 
Then \(\langle R_\gamma, R_\delta\rangle_R=0\) if and only if \(\gamma\neq\delta\) or \(R_\gamma=0\).
\begin{proof}
(\(\Rightarrow\)) 
If \(R_\gamma\neq 0\) and \(\gamma=\delta\), then for any non-zero \(x\in R_\gamma\) we have \(0<\langle x, x\rangle \in \langle R_\gamma,R_\delta\rangle\).

(\(\Leftarrow\)) 
Let \(J=\prid(R_1)\) as in Definition~\ref{def:idempotent} and note that the elements of \(J\) are orthogonal and sum to \(1\) by Lemma~\ref{lem:idempotent_isomorphism}.
Thus for all \(x,y\in R\) we have \(\langle x, y\rangle_R = \sum_{e\in J}\langle e x,ey\rangle_R\).
Hence to show \(\langle R_\gamma, R_\delta \rangle_R=0\) it suffices to show that \(\langle eR_\gamma, eR_\delta \rangle_{R} = 0\) for all \(e\in J\).
By Lemma~\ref{lem:hom_to_C} we have \(\langle eR_\gamma,eR_\delta\rangle_R=\langle eR_\gamma, eR_\delta \rangle_{eR}\).
Thus we may assume without loss of generality that \(\overline{R}=\overline{eR}\).
Then \(R_1\) is connected, so it is a field by Corollary~\ref{cor:spectrum_and_idempotents}. 
Thus \(\Gamma\) is a finite abelian group by Proposition~\ref{prop:basic_grad_props}.2.
If \(R_\gamma=0\) then \(\langle R_\gamma,R_\delta\rangle = 0\) and we are done.
Hence assume \(\delta\neq\gamma\), so there exists a character \(\chi\in\widehat{\Gamma}\) such that \(\chi(\gamma)\neq\chi(\delta)\) by Theorem~\ref{thm:pontryagin}.
Let \(n\) be the exponent of \(\Gamma\) and write \(R_\alpha'=R_\alpha\tensor_\Q\Q(\mu_n)\) for all \(\alpha\in\Gamma\).
Note that \(\widehat{\Gamma}\) acts on \(R'=R\tensor_\Q \Q(\mu_n)\), where \(\chi\in\widehat{\Gamma}\) acts by sending \(x\in R_\alpha'\) to \(\chi(\alpha)x\) for all \(\alpha\in\Gamma\).
Then \(\langle \chi * x, \chi*y\rangle_{R'} = \langle x,y\rangle_{R'}\) for all \(x,y\in R'\), as applying \(\chi\) just results in a reordering of the sum defining the inner product.
Then for any \(x\in R_\gamma'\) and \(y\in R_\delta'\), we have
\[ \langle x, y \rangle_{R'} = \langle \chi * x, \chi * y \rangle_{R'} = \langle \chi(\gamma) \cdot x, \chi(\delta) \cdot y \rangle_{R'} =  \langle \chi(\gamma\delta^{-1}) \cdot x, y \rangle_{R'}. \]
As \(s=1-\chi(\gamma\delta^{-1})\neq 0\) there exists some \(t\in\Z[\mu_n]\) such that \(st=n\) by Lemma~\ref{lem:root_diff_division}.2.
Thus
\[ 0 = \langle sR_\gamma', R_\delta'\rangle_{R'} \supseteq \langle st R_\gamma',R_\delta'\rangle_{R'} = n\langle R_\gamma', R_\delta' \rangle_{R'}\]
Then from Lemma~\ref{lem:hom_to_C} it follows that \(\langle R_\gamma,R_\delta\rangle_R = 0\), as was to be shown.
\end{proof}
\end{proposition}

\begin{corollary} \label{cor:orthogonal_homogenous}
Let \(R\) be a reduced order and let \(\overline{R}\) be a \(\Gamma\)-grading of \(R\). Then \(\langle R_\gamma, R_\delta\rangle_R=0\) if and only if \(\gamma\neq\delta\) or \(R_\gamma=0\).
\begin{proof}
Since \(R\) is a reduced order, \(R\tensor_\Z\Q\) is a product of number fields and \(\overline{E}=\overline{R}\tensor_\Z\Q\) is a \(\Gamma\)-grading of \(E\) such that \(E_\gamma\cap R= R_\gamma\). 
Any \(\sigma:R\to \C\) naturally extends to a morphism \(\sigma':E\to\C\) and each \(\rho:E\to\C\) restricts to a map \(\rho':R\to\C\), hence \(\langle x,y\rangle_E=\langle x,y\rangle_R\) for all \(x,y\in R\).
The claim then follows from Proposition~\ref{prop:orthogonal_homogenous}.
\end{proof}
\end{corollary}

\noindent\textbf{Proof of Theorem~\ref{UniGridGradThm}:}
Since \(R\) is a reduced order, it has a lattice structure with inner product \(\langle \us,\us\rangle_R\) as in Definition~\ref{hom_inner_product}.
Then Theorem 2 in \cite{Eichler1952} states that \(R\) has a decomposition \(\mathcal{E}=\{E_s\}_{s\in S}\) such that \(\langle E_s,E_t\rangle=0\) when \(s,t\in S\) are distinct, and such that \(\mathcal{E}\) is universal with this property, i.e.\ for any other such decomposition \(\mathcal{F}=\{F_t\}_{t\in T}\) of \(R\) there exists a unique map \(f:S\to T\) such that \(f_*\mathcal{E}=\mathcal{F}\).
We then obtain a universal grid grading of \(R\) by applying Proposition~\ref{prop:refinement_universal} to \(\mathcal{E}\), because Corollary~\ref{cor:orthogonal_homogenous} ensures all gradings of \(R\) have \(\mathcal{E}\) as refinement.
A universal group grading and abelian group grading then exist by Proposition~\ref{prop:hierarchy}.

\qed

\subsection{A proof using automorphism groups}

Instead of generalizing groups to grids, we will now focus on generalizing reduced orders to a broader class of rings.
To show \(R\) has a universal abelian group-grading, it suffices to show that \(X_\infty(R)\) as defined in Definition~\ref{def:X} is a finite abelian group by Theorem~\ref{thm:commutative_implies_universal}.

\begin{definition}
Let \(k\) be a commutative ring with \(p\in k\) and let \(M\) be a \(k\)-module.
A \(k\)-subalgebra \(A\subseteq \End_{k\catt{-Mod}}(M)\) is called {\em \(p\)-good} if \(p(1-ps)\) is injective for all \(s\in A\).
\end{definition}

\begin{lemma}\label{lem:root_p_good}
Let \(k\) be a commutative ring, let \(M\) be a \(k\)-module and let \(A\subseteq \End_{k\catt{-Mod}}(M)\) be a \(k\)-subalgebra.
Let \(p\in k\), \(u\in k^*\) and \(n\in\Z_{>0}\).
If \(A\) is \(up^n\)-good, then \(A\) is \(p\)-good.
\begin{proof}
Since multiplication by \(up^n\) is injective, so is multiplication \(p\).
If \((1-ps)(x)=0\) for \(s\in A\) and \(x\in M\), then \(x=ps(x)\) so \(x=(ps)^n(x)=p^ns^n(x)\).
Then \(r=u^{-1}s^n\in A\) and \(A\) is \(up^n\)-good.
Because \(up^n(1-up^nr)(x)=0\) we get that \(x=0\) and thus \(A\) is \(p\)-good.
\end{proof}
\end{lemma}

\begin{proposition}\label{prop:qp-proof}
Let \(p\) and \(q\) be (not necessarily distinct) primes, let \(k=\Z[\zeta_p,\zeta_q]\) and let \(M\) be a \(k\)-module.
Write \(E=\End_{k\catt{-Mod}}(M)\), \(\pi=1-\zeta_p\) and \(\tau=1-\zeta_q\) and let \(A\subseteq E\) be a \(p\)-good and \(q\)-good \(k\)-subalgebra.
\begin{enumerate}[topsep=0pt,itemsep=-1ex,partopsep=1ex,parsep=1ex,label={\em(\arabic*)}]
\item If \(\langle1+\pi\tau s\rangle\subseteq A^*\) for some \(s\in A\) is torsion, then \(s=0\).
\item If \(\langle1+\pi s,1+\tau r\rangle\subseteq A^*\) for some \(s,r\in A\) is torsion, then \(1+\pi s\) and \(1+\tau r\) commute.
\end{enumerate}
\begin{proof}
Note that by Lemma~\ref{lem:prime_roots}.1 and Lemma~\ref{lem:root_p_good} a \(p\)-good algebra is \(\pi\)-good.

(1) Write \(T=1+\tau k[X]\subseteq k[X]\).
Recall that \(\pi^{p-1}\in pk\) by Lemma~\ref{lem:prime_roots} and that \(p\divs\binom{p}{i}\) for \(0<i<p\).
Hence  
\[(1+\pi\tau T X)^p \subseteq 1 + p\pi\tau T \left( 1+\tau \sum_{i=2}^p \frac{\binom{p}{i}\pi^{i-1}\tau^{i-2}}{p} (XT)^{i-1} \right) X \subseteq 1+p\pi\tau TX. \tag{\ref{prop:qp-proof}.1}\]
Then inductively we have \((1+\pi\tau TX)^{p^n} \subseteq ( 1+p^{n-1}\pi\tau TX )^p \subseteq 1+p^n\pi\tau TX\), where we obtain the second inclusion by substituting \(p^{n-1}X\) for \(X\) in (\ref{prop:qp-proof}.1).

Consider \(\sigma=1+\pi\tau s\). 
By assumption we have \(\sigma^n=1\) for some \(n\in\Z_{>0}\).
If \(n=p^m\) for some \(m\in\Z_{\geq 0}\), then \(1=\sigma^n=1+p^m\pi\tau(1+\tau e)s\) for some \(e\in A\) by what we have shown before.
We have \(p^m\pi\tau(1+\tau e)s=0\) while \(p^m\pi\tau(1+\tau e)\) is injective by \(p\)-goodness and \(q\)-goodness, hence \(s=0\).
Thus we assume \(n\) and \(p\) are coprime, in which case there exists some \(m>0\) such that \(n\divs p^m-1\).
Then \(1+\pi\tau s=\sigma=\sigma^{p^m}=1+p^m\pi\tau(1+\tau e)s\) and thus \(\pi\tau(1-p^m(1+\tau e))s=0\).
As \(\pi\tau(1-p^m(1+\tau e))\) is injective we again have \(s=0\), as was to be shown.

(2) Let \(\sigma=1+\pi s\) and \(\rho=1+\tau r\).
By assumption \(\sigma^{-1}\in A\) hence \(s'=-\sigma^{-1}s\in A\) and \(\sigma(1+\pi s')=1=(1+\pi s')\sigma\), so \(\sigma^{-1}=1+\pi s'\).
As \(1=\sigma\sigma^{-1}=1+\pi(s+s'+\pi ss')\) we get \(\pi(s+s'+\pi ss')=0\).
Similarly \(\tau(r+r'+\tau rr')=0\) for \(\rho^{-1}=1+\tau r'\).
Now
\[[\sigma,\rho]\equiv (1+\pi s)(1+\tau r)(1+\pi s')(1+\tau r')\equiv 1+\pi(s+s'+\pi ss')+\tau(r+r'+\tau rr')\equiv 1 \mod \pi\tau.\]
The group \(\langle[\sigma,\rho]\rangle\subseteq\langle\sigma,\rho\rangle\) is torsion, hence we have \([\sigma,\rho]=1\) by (1).
\end{proof}
\end{proposition}

\begin{proposition}\label{prop:pgood_implies_commutative}
Let \(e\in\mathbb{S}\), let \(k=\Z[\mu_e]\), let \(M\) be a \(k\)-module and let \(A\subseteq E=\End_{k\catt{-Mod}}(M)\) be \(p\)-good \(k\)-subalgebra such that \((1-\zeta_p) A=A\cap (1-\zeta_p)E\) for all primes \(p\divs e\).
If \(\rho,\sigma\in A\) are \(\mu_e\)-diagonalizable and \(\langle\rho,\sigma\rangle\) is torsion, then \(\rho\) and \(\sigma\) commute.
\begin{proof}
By the Chinese remainder theorem we may write \(\sigma=\sigma_1\dotsm \sigma_s\) with \([\sigma_i,\sigma_j]=1\) and \(\sigma_i\in\langle\sigma\rangle\) of prime power order for all \(1\leq i,j\leq s\), and similarly we write \(\rho=\rho_1\dotsm \rho_r\). 
To show \([\sigma,\rho]=1\) it suffices to show that \([\sigma_i,\rho_j]=1\) for all \(i\) and \(j\). 
Thus we may assume without loss of generality that \(\ord(\sigma)=p^a\) and \(\ord(\rho)=q^b\) for \(a,b\in\Z_{\geq0}\) and primes \(p,q\divs e\).

We apply induction to \(n=a+b\).
First assume \(a,b\leq 1\).
Consider \(\pi=1-\zeta_p\) and recall that \(1-\zeta_p^i\equiv 0 \mod \pi\) for all \(i\in\Z\) by Lemma~\ref{lem:prime_roots}.
Thus \(1-\sigma\equiv 0 \mod \pi\) by \(\mu_p\)-diagonalizability.
Hence \(\sigma=1+\pi s\) for some \(s\in E\) and thus \(s\in A\) by \(\pi A=A\cap\pi E\).
Similarly \(\rho=1+\tau r\) with \(\tau=1-\zeta_q\), hence \([\sigma,\rho]=1\) by Proposition~\ref{prop:qp-proof}.2.
Otherwise we may assume without loss of generality that \(a\geq 2\).
Then \([\sigma^p,\rho]=1\) by the induction hypothesis.
Fix \(\zeta\in\mu_{p^a}\) and consider \(N=M(\sigma^p,\zeta^p)\) and the \(k\)-subalgebra \(B=\{\sigma\in A\,|\, \sigma N\subseteq N\}\subseteq A\).
Then \(\sigma\in B\) by construction of \(N\) and \(\rho\in B\) by Lemma~\ref{com_diag}.
Since \(N\) is a direct summand of \(M\) we get that \(C=B|_N\subseteq\End_{k\catt{-Mod}}(N)\) is \(r\)-good and \((1-\zeta_r)C=C\cap(1-\zeta_r)\End_{k\catt{-Mod}}(N)\) for all primes \(r\divs e\).
Let \(\alpha=\rho|_N\) and \(\beta=(\zeta^{-1} \sigma)|_N\), which are respectively \(\mu_{p^a}\)-diagonalizable and \(\mu_q\)-diagonalizable.
Moreover, \(\langle \alpha,\beta\rangle\subseteq C^*\) is a torsion subgroup. 
We have \(\ord(\alpha)=1<\ord(\sigma)\) and \(\ord(\beta)\leq\ord(\rho)\), so we may apply the induction hypothesis to conclude \(\alpha\) and \(\beta\) commute.
In particular, \(\rho\) and \(\sigma\) commute when restricted to \(M(\sigma^p,\zeta^p)\) for each \(\zeta\), so \(\sigma\) and \(\rho\) commute on \(M\) by \(\mu_{p^{a-1}}\)-diagonalizability of \(\sigma^p\). 
Hence \(\langle \rho,\sigma\rangle\) is abelian when it is torsion.
\end{proof}
\end{proposition}

\begin{lemma}\label{lem:tensor_is_reduced}
Let \(e\in\mathbb{S}\) and let \(R\) be a commutative ring such that multiplication by integer divisors of \(e\) is injective on \(R\).
If \(R\) is reduced, then \(R\tensor_\Z\Z[\mu_e]\) is reduced.
\begin{proof}
Since \(R\tensor_\Z\Z[\mu_e]=\lim_{\rightarrow d\in\Z_{>0},\ d\divs e} R\tensor_\Z\Z[\mu_d]\) we may assume without loss of generality that \(e\in\Z_{>0}\). 
With \(S=\{e^m\,|\,m\in\Z_{\geq0}\}\) we consider the localization \(S^{-1}R\), which is reduced, and note that \(R\) injects into \(S^{-1}R\) since multiplication by \(e\) is injective. 
Since \((S^{-1}R)\tensor_\Z\Z[\mu_e]\cong S^{-1}(R\tensor_\Z\Z[\mu_e])\), we may also assume without loss of generality that \(R=S^{-1}R\).
Let \(\p\in\minspec R\), let \(K_\p\) be a field containing \(R/\p\) and note that \(\text{char}(K_\p)\nmid e\) as \(e\) is invertible in \(R/\p\).
We have that \(X^e-1\) and thus \(\Phi_e\), the \(e\)-th cyclotomic polynomial, is separable over \(K_\p\) since its derivative \(eX^{e-1}\) only has \(0\) as root.
Hence \(K_\p\tensor_\Z\Z[\mu_e]\cong K_\p[X]/\Phi_e\) is reduced by Theorem~A1.3 in \cite{Eisenbud}.
By Lemma~\ref{basic_comalg_wreath} we have an injection 
\[R\tensor_\Z\Z[\mu_e]\to \bigg( \prod_{\p\in\minspec R} K_\p \bigg)\tensor_\Z\Z[\mu_e] \cong \prod_{\p\in\minspec R}(K_\p\tensor_\Z\Z[\mu_e]), \]
where the isomorphism follows from the fact that \(\Z[\mu_e]\) is a free \(\Z\)-module of finite rank.
It follows that \(R\tensor_\Z\Z[\mu_e]\) is reduced as well.
\end{proof}
\end{lemma}

\begin{lemma}\label{lem:Noetherian_implies_finite_group}
Let \(k\) be a commutative Noetherian ring and let \(R\) be a reduced commutative \(k\)-algebra that is finitely generated as a module.
Then \(\Aut_{k\catt{-Alg}}(R)\) is a finite group.
\begin{proof}
Since \(R\) is reduced we have an injection \(R\to\prod_{\p\in\minspec R} R/\p\) by Lemma~\ref{basic_comalg_wreath} and thus an injection 
\[\Aut_{k\catt{-Alg}}(R)\hookrightarrow \bigg(\prod_{\p\in\minspec R} \Aut_{k\catt{-Alg}}(R/\p)\bigg)\rtimes\Aut_\catt{Set}(\minspec R).\]
We have by Corollary~1.3 in \cite{Eisenbud} that \(R\) is a Noetherian ring and thus \(\#\minspec R<\infty\) by Exercise~1.2 in \cite{Eisenbud}, hence it suffices to show that \(\Aut_{k\catt{-Alg}}(R/\p)\) is finite for all \(\p\in\minspec R\). As \(R/\p\) is a finitely generated \(k\)-algebra, we may then assume without loss of generality that \(R\) is a domain.

It suffices to show that each element \(x\in R\) has a finite orbit under \(\Aut_{k\catt{-Alg}}(R)\), because we may then apply this fact to the finite number of generators of \(R\).
Note that we may replace \(k\) by its image in \(R\) and thus assume that \(k\) is a domain.
Because \(R\) is finitely generated, \(x\) is the root of some non-zero polynomial \(f\in k[X]\).
As \(\Aut_{k\catt{-Alg}}(R)\) must permute the roots of this polynomial and \(f\) has only finitely many roots because \(R\) is a domain, we may conclude that the orbit of \(x\) is finite, as was to be shown. 
\end{proof}
\end{lemma}

\begin{lemma}\label{lem:decomp_bound}
Let \(k\) be a commutative Noetherian ring and let \(M\) be a finitely generated \(k\)-module.
There exists some constant \(c\in\Z_{\geq0}\) such that if \(\mathcal{M}\) is a decomposition of \(M\), then \(\mathcal{M}\) contains at most \(c\) non-zero elements.
\begin{proof}
Let \(n\in\Z_{\geq 0}\) be such that \(M\) can be generated by \(n\) elements.
First assume \(k\) is a local ring with maximal ideal \(\m\).
If \(N/\m N=0\) for some submodule \(N\subseteq M\), then \(N=0\) by Nakayama's lemma.
Hence if \(\mathcal{M}=\{M_i\}_{i\in I}\) is a decomposition of \(M\), then 
\[ n \geq \dim_{k/\m\catt{-Vec}}( M/\m M ) = \sum_{i\in I} \dim_{k/\m\catt{-Vec}}(M_i/\m M_i) \]
implies \(M_i/\m M_i=0\) and thus \(M_i=0\) for all but at most \(n\) elements \(i\in I\).

Now consider the general case. 
Let \(A=\text{Ass}_k(M)\) be the set of associated primes of \(M\), which is finite by Theorem 3.1 in \cite{Eisenbud} since \(M\) is finitely generated over a Noetherian ring.
We have an injection \(M\to\prod_{\p\in A} M_\p\) by Corollary 3.5.a in \cite{Eisenbud} and for any submodule \(N\subseteq M\) we have \(\text{Ass}_k(N)\subseteq\text{Ass}_k(M)\).
Let \(\mathcal{M}=\{M_i\}_{i\in I}\) be a decomposition of \(M\).
For each \(\p\in A\) we have a decomposition \(\{ M_{i,\p} \}_{i\in I}\) of the \(k_\p\)-module \(M_\p\) and thus \(M_{i,\p}=0\) for all but at most \(n\) elements \(i\in I\) as shown before.
For each \(i\in I\) such that \(M_i\neq 0\) there must exist some \(\p\in A\) such that \(M_{i,\p}\neq 0\), hence there are at most \(n\cdot\# A\) non-zero elements in \(\mathcal{M}\) 
\end{proof}
\end{lemma}

\begin{theorem}\label{thm:pgood_theorem}
Let \(e\in\mathbb{S}\), let \(k\) be a commutative Noetherian ring and let \(R\) be a reduced commutative \(k\)-algebra that is finitely generated as \(k\)-module, such that for each prime \(p\divs e\) we have \(\bigcap_{n\geq 0}p^n R = \{0\}\) and that multiplication by \(p\) is injective on \(R\).
Then \(R\) has a universal abelian \(e\)-torsion grading.
\begin{proof}
By Lemma~\ref{lem:decomp_bound} there exists some \(c\in\Z_{\geq 0}\) such that each grading of \(R\) has at most \(c\) non-zero homogenous components.
Hence any efficient abelian group-grading of \(R\) has exponent dividing \(c!\).
By replacing \(e\) with \(\gcd(e,c!)\) we may assume that \(e\in\Z_{>0}\).
Then \(k'=k\tensor_\Z\Z[\mu_e]\) is a Noetherian ring because \(\Z[\mu_e]\) is a finitely generated \(\Z\)-algebra, and thus \(R'=R\tensor_\Z\Z[\mu_e]\) is a Noetherian ring because it is finitely generated over \(k'\).
Since \(\Z[\mu_e]\) is a free \(\Z\)-module, we also have that \(\bigcap_{n\geq 0}p^n R' = \{0\}\) and that multiplication by \(p\) is injective on \(R'\) for all primes \(p\divs e\).
Let \(s\in\End_{k'\catt{-Mod}}(R')\) and assume \(p(1-ps)(x)=0\) for some \(x\in R'\) and prime \(p\divs e\).
Then \(x=ps(x)\) by injectivity of \(p\), so \(x=p^ns^n(x)\) for all \(n\in\Z_{\geq0}\).
Thus \(x\in\bigcap_{n\geq 0}p^nR'=\{0\}\) so \(x=0\) and thus \(p(1-ps)\) is injective.
It follows that \(\End_{k'\catt{-Mod}}(R')\) is \(p\)-good for all primes \(p\divs e\).
By Lemma~\ref{lem:tensor_is_reduced} have that \(R'\) is reduced and thus by Lemma~\ref{lem:Noetherian_implies_finite_group} we have that \(\Aut_{k'\catt{-Alg}}(R')\) is a finite group.
In particular, all \(\mu_e\)-diagonalizable elements of \(\Aut_{k'\catt{-Alg}}(R')\) commute by Proposition~\ref{prop:pgood_implies_commutative}.
It follows that \(X_e(R)\) as definined in Definition~\ref{def:X} is a finite abelian group.
Hence by Theorem~\ref{thm:commutative_implies_universal} there exists a universal abelian \(e\)-torsion grading of \(R\), as was to be shown.
\end{proof}
\end{theorem}

Note that reduced orders satisfy the conditions to Theorem~\ref{thm:pgood_theorem} for all \(e\in\mathbb{S}\) and thus have a universal abelian-group grading.
The condition that \(k\) be a Noetherian ring is rather strong, but Example~\ref{ex:non-compact_counter} shows that theorem becomes false when we drop it.

\clearpage
\section{Universal abelian group grading}\label{sec:abelian_algorithm}

In this section we present a proof of Theorem~\ref{AlgAbGrpThm} and derive Theorem~\ref{AlgGridThm} for abelian groups.

\subsection{Combinatorial algorithm}

In this section we consider a finite groupoid \(\mathcal{C}\) and a Steinitz number \(e\in\mathbb{S}\) such that the group \((\hat{\Z}/e\hat{\Z})^*\) acts on \(\mathcal{C}\).
Our goal is to study the following set, specifically its size and how to list its elements.

\begin{definition}\label{def_Xe}
Let \(\mathcal{C}\) be a finite groupoid, \(e\in\mathbb{S}\) and \(G\subseteq (\hat{\Z}/e\hat{\Z})^*\) a subgroup that has a left action on \(\mathcal{C}\).
Then \(G\) has an induced left-action on \(\wreath(\mathcal{C})\) by Remark~\ref{rem:induced_action_on_wreath} and a right action on the \(e\)-torsion of \(\wreath(\mathcal{C})\) by exponentiation.
With respect to these actions we define
\[ \Xe(\mathcal{C},G) = \{ \sigma\in\wreath(\mathcal{C}) \,|\, \sigma^e = 1,\, (\forall a\in G)\ {}^a\sigma = \sigma^a \}. \]
Beware that this notation hides the left action of \(G\) on \(\mathcal{C}\), which will remain unchanged throughout this section.
\end{definition}

In this section we will derive the following result.

\begin{theorem}\label{main_theorem_cat}
There is a deterministic algorithm that takes a quadruple \((\mathcal{C},e,G,*)\), where \(\mathcal{C}\) is a finite groupoid encoded as in Section~\ref{sec:wreath_conjugacy_computation}, where \(e>1\) is a prime power encoded in base \(1\), where \(G\subseteq(\Z/e\Z)^*\) is a subgroup and \(*\) is an action of \(G\) on \(\mathcal{C}\) which we encode by associating to each element of \(G\) a permutation on the morphisms of \(\mathcal{C}\), and computes \(X_e(\mathcal{C},G)\) in \(n^{O(m)}\) time, where \(n\) is the length of the input and \(m=\#(\obj(\mathcal{C})/G)\).  
With \(a=\#\aut(\mathcal{C})\) as defined in Definition~\ref{def:cat_morphisms} and \(c=2+\mathbbm{1}( 4 \divs e )\) we have the bound \(\# X_e(\mathcal{C},G) \leq (2ma^c)^m\). 
\end{theorem}

We will prove this using various reductions.
First we reduce to the case where all objects of \(\mathcal{C}\) are \((\mathcal{C}\rtimes G)\)-isomorphic.

\begin{lemma}\label{category_isomorphism_reduction}
With the same setting as in Definition~\ref{def_Xe}, consider the group \(\wreath(\mathcal{C})\rtimes G\) and let \(\mathcal{P}\) be the set of orbits of \(\wreath(\mathcal{C})\rtimes G\) acting on \(\obj(\mathcal{C})\). Then with \(\mathcal{C}|_P\) the full subcategory of \(\mathcal{C}\) with objects \(P\subseteq\obj(\mathcal{C})\) we have
\[ \Xe(\mathcal{C},G) = \prod_{P\in\mathcal{P}} \Xe(\mathcal{C}|_P, G). \tag{\ref{category_isomorphism_reduction}.1} \]
\begin{proof}
All \(P\in\mathcal{P}\) are \(G\)-invariant by definition, hence the action of \(G\) does in fact restrict to \(\mathcal{C}|_P\), making the right hand side of (\ref{category_isomorphism_reduction}.1) well-defined.
It follows directly from the definition that the right hand side of (\ref{category_isomorphism_reduction}.1) is contained in the left hand side.
For the converse, it suffices to note that each element \(\sigma\in\Xe(\mathcal{C},G)\) has orbits refining \(\mathcal{P}\) and that the relations \(\sigma^e=1\) and \({}^a\sigma=\sigma^a\) also hold for the restrictions \(\sigma|_P\in \wreath(\mathcal{C}|_P)\) for \(P\in\mathcal{P}\).
Hence \(\sigma|_P\in \Xe(\mathcal{C}|_P, G)\) and \(\sigma\in\prod_{P\in\mathcal{P}} \Xe(\mathcal{C}|_P,G)\).
\end{proof}
\end{lemma}

We may thus assume the groups we consider are `sufficiently transitive', i.e. \(\wreath(\mathcal{C})\rtimes G\) acts transitively on \(\obj(\mathcal{C})\).
In the next reduction we will strengthen this property to the individual elements of \(\wreath(\mathcal{C})\).
Namely, we note that for \(\sigma\in\Xe(\mathcal{C},G)\) the group \(G\) acts on \(\langle\sigma\rangle\) by exponentiation, hence we may define the following.

\begin{definition}
With the same setting as in Definition~\ref{def_Xe}, define
\[\Xew(\mathcal{C},G) = \{ \sigma\in\Xe(\mathcal{C},G) \,|\, \langle\sigma\rangle\rtimes G \text{ acts transitively on \(\obj(\mathcal{C})\)} \}.\]
\end{definition}

For partitions \(\mathcal{P}\) and \(\mathcal{O}\) of a set \(\Omega\) we write \(\mathcal{P}\succeq\mathcal{O}\) if and only if \(P\cap O\in\{\emptyset,O\}\) for all \(O\in\mathcal{O}\) and \(P\in\mathcal{P}\), i.e. when \(\mathcal{O}\) refines \(\mathcal{P}\).

\begin{lemma}\label{category_wide_reduction}
Consider the same setting as in Definition~\ref{def_Xe}. 
Assume \(\wreath(\mathcal{C})\rtimes G\) acts transitively on \(\obj(\mathcal{C})\) and let \(\mathcal{O}=\obj(\mathcal{C})/G\). Then
\[ \Xe(\mathcal{C},G) = \bigsqcup_{\mathcal{P}\succeq\mathcal{O}} \prod_{P\in\mathcal{P}} \Xew(\mathcal{C}|_P,G). \tag{\ref{category_wide_reduction}.1}\]
\begin{proof}
For all \(P\in\mathcal{P}\succeq\mathcal{O}\) the set \(P\) is \(G\)-stable, hence the action of \(G\) restricts to \(\mathcal{C}|_P\), making \(\Xew(\mathcal{C}|_P,G)\) well-defined.
Write \(X(\mathcal{P})=\prod_{P\in\mathcal{P}} \Xew(\mathcal{C}|_P,G) \subseteq \wreath(\mathcal{C})\) for all partitions \(\mathcal{P}\succeq\mathcal{O}\).
For \(\sigma\in X(\mathcal{P})\) one clearly has that the set of orbits of \(\langle\sigma\rangle\rtimes G\) is precisely \(\mathcal{P}\).
Hence for distinct partitions \(\mathcal{P}\) and \(\mathcal{Q}\) one has that \(X(\mathcal{P})\cap X(\mathcal{Q})=\emptyset\) and the union on the right side of (\ref{category_wide_reduction}.1) is indeed disjoint.
It follows directly from the definition that \( X(\mathcal{P}) \subseteq \Xe(\mathcal{C},G)\), so the right hand side of (\ref{category_wide_reduction}.1) is contained in the left.
Conversely, for any \(\sigma\in \Xe(\mathcal{C},G)\) the partition \(\mathcal{P}=\obj(\mathcal{C})/(\langle\sigma\rangle\rtimes G)\) must have \(\mathcal{O}\) as refinement.
As for all \(P\in\mathcal{P}\) the group \(\langle \sigma|_P\rangle \rtimes G\) acts transitively on \(P\), we have that \(\sigma|_P\in \Xew(\mathcal{C}|_P,G)\) and \(\sigma\in X(\mathcal{P})\), so the reverse inclusion also holds. 
\end{proof}
\end{lemma}

A final reduction step will restrict to the case where \(\langle\sigma\rangle\) itself is transitive.
Here we will use that \(G\) is abelian.

\begin{definition}
With the same setting as in Definition~\ref{def_Xe}, define
\[\Xet(\mathcal{C},G) = \{ \sigma\in \Xe(\mathcal{C},G) \,|\, \langle\sigma\rangle\text{ acts transitively on }\obj(\mathcal{C})\}.\]
\end{definition}

We will now show for \(\sigma\in\Xe(\mathcal{C},G)\) that the restriction of \(\sigma\) to a single orbit of \(\langle\sigma\rangle\) defines it in its entirety.

\begin{proposition}\label{transitive_reduction_prop}
Consider the same setting as in Definition~\ref{def_Xe}.
Let \(\mathcal{O}=\obj(\mathcal{C})/G\), let \(I_O\trianglelefteq G\) be the pointwise stabiliser of \(O\in\mathcal{O}\) and let \(I=\langle I_O \,|\, O\in \mathcal{O}\rangle\).
Write \(\mathcal{H}=\{H\,|\, I \trianglelefteq H \trianglelefteq G \}\) and fix some \(x\in\obj(\mathcal{C})\). 
For any \(H\in\mathcal{H}\) let
\[ \mathcal{P}(H) = \{ P \subseteq \obj(\mathcal{C}) \,|\, x\in P,\, (\forall O\in\mathcal{O})\, P\cap O \in \obj(\mathcal{C})/H \}. \]
Then for all \(\sigma\in\Xew(\mathcal{C},G)\) we have that \(P=\langle\sigma\rangle x\in \mathcal{P}(H_P)\), where \(H_P\trianglerighteq I\) is the setwise stabiliser of \(P\).
Moreover, the map
\begin{align*} 
\varphi:\Xew(\mathcal{C},G) &\to \bigsqcup_{H\in\mathcal{H}} \bigsqcup_{ P \in \mathcal{P}(H) } \Xet( \mathcal{C}|_{P}, H ) \\
\sigma &\mapsto \sigma|_{\langle \sigma \rangle x} \in \Xet(\mathcal{C}|_{\langle\sigma\rangle x},H_{\langle\sigma\rangle x})
\end{align*}
is well-defined and bijective.
\begin{proof}
We break up the proof into several parts.

{\em Well-definedness:} 
Let \(\sigma\in \Xew(\mathcal{C},G)\) be given, let \(P=\langle \sigma \rangle x\) be the orbit of \(x\) under \(\sigma\) and let \(H\trianglelefteq G\) be the setwise stabilizer of \(P\). 
Note that \(H\) acts on \(P\) and thus \(\Xet(\mathcal{C}|_P,H)\) is well-defined and \(\sigma|_P\in \Xet(\mathcal{C}|_P,H)\).
First we show that \(H\in\mathcal{H}\) and \(P\in\mathcal{P}(H)\), from which it follows that \(\sigma|_P\) is in fact an element of the codomain of \(\varphi\).
Since \(G\) normalizes \(\langle\sigma \rangle\) it permutes the orbits of \(\sigma\).
As \(P\cap O \neq \emptyset\) for all \(O\in\mathcal{O}\) by definition of \(\Xew(\mathcal{C},G)\), we have that since \(I_O\) stabilises a point \(y\in P\cap O\), it must stabilise the entire orbit \(P\) setwise.
Hence \(I_O \trianglelefteq H\) for all \(O\in\mathcal{O}\), so \(H\in\mathcal{H}\).
Let \(O\in\mathcal{O}\) be given.
Note that \(HP = P\) by definition and \(HO=O\) as \(H\subseteq G\), and thus \(H(P\cap O)=P\cap O\).
Furthermore, for all \(y,z\in P\cap O\) there exists some \(t \in G\) such that \(t y = z\) by transitivity.
By the same argument as before we have \(t P = P\), and thus \(t\in H\) and \(Hy=Hz\).
Thus \(P\cap O\) is an orbit of \(H\) acting on \(\obj(\mathcal{C})\) for all \(O\in\mathcal{O}\), so \(P\in\mathcal{P}(H)\).

For the well-definedness of \(\varphi\) it remains to show that the union defining the codomain of \(\varphi\) is indeed disjoint.
If \(P,P'\subseteq\obj(\mathcal{C})\) are distinct, then clearly \(\Xet(\mathcal{C}|_P,H)\cap \Xet(\mathcal{C}|_{P'},H')=\emptyset\) for any \(H,H'\trianglelefteq G\) for which this is defined. Thus we need to show that for distinct \(H,H'\in\mathcal{H}\) we have \(\mathcal{P}(H)\cap\mathcal{P}(H')=\emptyset\).
We have \(I_O\trianglelefteq H,H'\) for \(O=Gx\in\mathcal{O}\), so \(H/I_O\neq H'/I_O\).
As \(G/I_O\) acts regularly on \(O\), we have that \((H/I_O) x\neq (H'/I_O) x\), hence for any \(P\in\mathcal{P}(H)\) and \(P'\in\mathcal{P}(H')\) we have \(P\cap O = Hx \neq H' x = P'\cap O\).
So \(P\neq P'\), from which the claim follows.
We conclude that \(\varphi\) is well-defined.

{\em Injectivity:}
Let \(\sigma\in\Xew(\mathcal{C},G)\) and \(\omega=\varphi(\sigma) \in \Xet(\mathcal{C}|_P,H)\) for \(P=\langle\sigma\rangle x\) and \(H=H_P\).
For all \(y\in \obj(\mathcal{C})\) there exists some \(z\in P\) and \(\tau=\tau_a\in G\) such that \(t y= z\).
We use the convention of Section~\ref{sec:wreaths} that \(\sigma=((\sigma_y)_{y\in\obj(\mathcal{C})},s)\) and similarly for \(\tau\) and \(\omega\) is implicit.
We have 
\begin{align*} 
s(y) &= s t^{-1} (z) = t^{-1} s^a(z) = t^{-1} o^a (z) \quad\text{and} \tag{\ref{transitive_reduction_prop}.1}\\
\sigma_y &=  ( \tau^{-1} \sigma^a \tau )_y = (\tau_{s(y)})^{-1} \circ (\omega^a)_{z} \circ \tau_y. \tag{\ref{transitive_reduction_prop}.2}
\end{align*}
Hence \(\sigma\) is uniquely determined by \(\varphi(\sigma)\) and \(\varphi\) is injective.

Note that equations (\ref{transitive_reduction_prop}.1) and (\ref{transitive_reduction_prop}.2) in fact give us a way to construct, for any \(\omega\in\Xet(\mathcal{C}|_P,H)\), an element \(\sigma\in\wreath(\mathcal{C})\).
This construction a priori depends on the choice of \(\tau\) and \(z\).
If we choose for each class \(gH\in G/H\) a representative \(\tau_{gH}\) and require that \(\tau\) in (\ref{transitive_reduction_prop}.1) and (\ref{transitive_reduction_prop}.2) be taken in \(R=\{\tau_{gH}\,|\, gH\in G/H\}\), then for each \(y\in\obj(\mathcal{C})\) both \(\tau\) and \(z\) exist and are unique.
Now we obtain a map \(\psi_{H,P}:\Xet(\mathcal{C}|_P,H_P)\to \wreath(\mathcal{C})\), and a combined map \(\psi:\bigsqcup_{H\in\mathcal{H}}\bigsqcup_{P\in\mathcal{P}(H)} \Xet(\mathcal{C}|_P,H_P)\to \wreath(\mathcal{C})\) which is a left inverse of \(\varphi\).

{\em Surjectivity:}
For surjectivity of \(\varphi\), it remains to show that \(\psi\) is also a right inverse of \(\varphi\).
This is, however, clearly the case when \(\im(\psi)\subseteq \Xew(\mathcal{C},G)\).
Let \(H\in\mathcal{H}\), \(P\in\mathcal{P}(H)\) and \(\omega\in\Xet(\mathcal{C}|_P,H)\).
By definition, \(\sigma:=\psi(\omega)\) has \(P\) as an orbit.
For all \(y\in \obj(\mathcal{C})\) with corresponding \(\tau_a\) and \(z\) we have
\begin{align*} 
s^e(y) &= t^{-1} o^{ae}(z) = t^{-1} z = y  \quad\text{and} \\
(\sigma^e)_y &= \tau_{t^{-1}s(y)} \circ (\omega^{ae})_z \circ (\tau^{-1})_y = \tau_{t^{-1}(y)} \circ (\tau^{-1})_y = \id_y,
\end{align*}
so \(\sigma^e=1\).
Now let additionally \(\tau_b\in G\) be given. Then there exists some \(\tau_c\in R\) and \(\tau_d\in H\) such that \(\tau_b\tau_{a}^{-1} = \tau_c\tau_d\).
Hence
\begin{align*}
t_bs(y) &= t_b t_a^{-1} o^a(z) = t_c t_d o^a(z) = t_c o^{ad} t_d(z) = s^{adc} t_c t_d(z) = s^b t_b(y) \quad\text{and} \\
(\tau_b \sigma)_y &= (\tau_b\tau_a \omega^a \tau_a^{-1})_y = (\tau_c\tau_d \omega^a \tau_a^{-1})_y = (\tau_c \omega^{ad} \tau_d \tau_{a}^{-1})_y = (\sigma^{adc})_{t_d(y)} \circ (\tau_d)_y = (\sigma^b \tau_d)_y.
\end{align*}
Thus \(\sigma\in \Xew(\mathcal{C},G)\), as was to be shown. Therefore \(\psi\) is a right inverse of \(\varphi\).
It follows that \(\varphi\) is surjective. 
\end{proof}
\end{proposition}


\begin{lemma}\label{category_leaf_lemma}
Consider the same setting as in Definition~\ref{def_Xe}, let \(a=\#\aut(\mathcal{C})\) as defined in Definition~\ref{def:cat_morphisms} and let \(m\) be the number of orbits of \(G\) acting on \(\obj(\mathcal{C})\). 
If \(G\) is cyclic, then \(\#\Xet(\mathcal{C},G)\leq a^{m+1}\).
In this case, we may compute \(\Xet(\mathcal{C},G)\) in \(n^{O(m)}\) time with \(n\) the length of the input \((\mathcal{C},e,G,*)\), where \(e\in\Z_{>0}\) and \(*\) is the action of \(G\) on \(\mathcal{C}\).
\begin{proof}
Write \(G=\langle b\rangle\) and take any \(x\in\obj(\mathcal{C})\). Consider \(\sim_\mathcal{C}\) as in Definition~\ref{def:conjugacy} and \(\lambda\) as in Definition~\ref{def:lambda_sigma}.
Let 
\[S=\{\sigma\in\wreath(\mathcal{C}) \,|\, \sigma \text{ acts transitively on \(\obj(\mathcal{C})\), \(\sigma^e=1\) and \({}^b\lambda_\sigma(x) \sim_\mathcal{C} \lambda_\sigma(x)^b\)}\}\]
and note that it is closed under conjugation. Moreover, we have that \(\Xet(\mathcal{C},G)\subseteq S\).
Let \(\Phi\) be the set of \(\wreath(\mathcal{C})\)-conjugacy classes of \(S\), which corresponds to some subset of the set of conjugacy classes of \(\Aut_\mathcal{C}(x)\) by Theorem~\ref{conjugacy_theorem}.
Fix some representative \(\sigma_\varphi \in\varphi\) for each \(\varphi\in\Phi\).
Let \(\Gamma=\wreath(\mathcal{C}) \rtimes G\), let \(\tau_b\) be the image of \(b\) in \(\Gamma\) and let \(U_{\varphi}=\{\upsilon \in \wreath(\mathcal{C})\circ \tau_b \subseteq \Gamma \,|\, \upsilon \sigma_\varphi \upsilon^{-1} = \sigma_\varphi^b \}\).
If \(\upsilon_\varphi\in U_\varphi\), then all permutations \(\alpha\in\wreath(\mathcal{C})\) such that \(\alpha \upsilon_\varphi \alpha^{-1} = \tau_b \) give rise to an element \(\alpha \sigma_\varphi \alpha^{-1}\in \Xet(\mathcal{C},G)\), and we claim that we obtain all elements of \(\Xet(\mathcal{C},G)\) in this way.
Each \(\sigma\in\Xet(\mathcal{C},G)\) is conjugate to some \(\sigma_\varphi\), so there exists some \(\alpha\in\wreath(\mathcal{C})\) such that \(\alpha\sigma_\varphi \alpha^{-1}=\sigma\), while \(\alpha^{-1}\tau_b\alpha=\alpha^{-1} ({}^{b}\alpha) \tau_b\in U_\varphi\). 
Hence these are in fact all elements of \(\Xet(\mathcal{C},G)\).

It remains to count these elements.
Given \(\varphi\), we have that 
\[U_\varphi = \{\beta\in \wreath(\mathcal{C}) \,|\, \beta ({}^b\sigma_\varphi) \beta^{-1} = \sigma_\varphi^b \} \circ \tau_b.\]
We may compute \(U_\varphi\) using Theorem~\ref{conjugacy_theorem} in time \(n^{O(m)}\) and it contains \(\#\obj(\mathcal{C})\cdot \#C(\lambda_{\sigma_\varphi}(x))\) elements.
Furthermore \(\sum_{\varphi\in\Phi} \# C(\lambda_{\sigma_\varphi}(x))\leq\#\Aut_\mathcal{C}(x)\), so \(\sum_{\varphi\in\Phi} \# U_\varphi \leq \#\obj(\mathcal{C}) \cdot \#\Aut_\mathcal{C}(x) = a\).
Let then \(V(\upsilon) = \{ \alpha\in\wreath(\mathcal{C}) \,|\, \alpha \upsilon \alpha^{-1} = \tau_b \}\), which for any \(\upsilon\in\Gamma\) may be computed using Corollary~\ref{cor:semidirect_conjugacy_theorem} in time \(n^{O(m)}\) and contains at most \(a^m\) elements.
Hence, counting the number of elements we get
\[ \#\Xet(\mathcal{C},G) \leq \sum_{\varphi\in\Phi} \sum_{\upsilon\in U_\varphi} \#V(\upsilon) \leq a^m \cdot \sum_{\varphi\in\Phi} \# U_\varphi \leq a^{m+1},\] 
as was to be shown.
\end{proof}
\end{lemma}

\noindent\textbf{Proof of Theorem~\ref{main_theorem_cat}.}\\
\noindent
Write \(e=p^k\) for some prime \(p\) and let \(d=1+\mathbbm{1}(4\divs e)\).
By Lemma~\ref{lem:ugly_cite} the group \(G\subseteq(\Z/p^k\Z)^*\) has a cyclic subgroup \(G'\trianglelefteq G\) such that \([G:G']\leq d\). 
Note that \(\Xet(\mathcal{C},G)\subseteq\Xet(\mathcal{C},G')\) and \(\#\obj(\mathcal{C})/G'\leq dm\).
Then by Lemma~\ref{category_leaf_lemma} we have
\[\#\Xet(\mathcal{C},G)\leq\#\Xet(\mathcal{C},G')\leq a^{dm+1}. \tag{\ref{main_theorem_cat}.1}\]

We combine (\ref{main_theorem_cat}.1) with Proposition~\ref{transitive_reduction_prop}. 
For \(H\in\mathcal{H}\) we have \(\#\mathcal{P}(H)=[G:H]^{m-1}\) and assuming all objects of \(\mathcal{C}\) are \(\mathcal{C}\rtimes G\)-isomorphic we get \(\#\aut(\mathcal{C}|_P)=a/[G:H]\) for all \(P\in\mathcal{P}(H)\).
Thus
\[\#\Xew(\mathcal{C},G)\leq \sum_{H\in\mathcal{H}} \sum_{P\in\mathcal{P}(H)} \left(\frac{a}{[G:H]}\right)^{dm+1} = a^{dm+1} \sum_{H\in\mathcal{H}} [G:H]^{-2-(d-1)m} \leq 2a^{dm+1}. \tag{\ref{main_theorem_cat}.2}\]
In the last inequality we use \(S=\sum_{H\trianglelefteq G}[G:H]^{-2}\leq 2\).
Namely, if \(p>2\) then \(G\) is cyclic and \(H\) is uniquely determined by its size, so we have \(S\leq \sum_{i=1}^\infty i^{-2}=\pi^2/6\leq 2\).
If \(p=2\), then for each index greater than \(1\) there are at most 3 subgroups of \(G\) of this index, so \(S\leq1+3\sum_{i=1}^\infty 2^{-2i}=2\).

Still assuming all objects of \(\mathcal{C}\) are \(\mathcal{C}\rtimes G\)-isomorphic we combine (\ref{main_theorem_cat}.2) with Lemma~\ref{category_wide_reduction}.
Then we get for fixed \(\mathcal{P}\succeq \obj(C)/G\) that
\[\#\prod_{P\in\mathcal{P}} \Xew(\mathcal{C}|_P,G) \leq \prod_{P\in\mathcal{P}} 2(\#\aut(\mathcal{C}|_P))^{d\cdot\#(P/G)+1} \leq 2^m a^{\sum_{P\in\mathcal{P}}(d\#(P/G)+1)} = (2 a^c)^{m}.\]
There are at most \(m^m\) such \(\mathcal{P}\succeq\obj(\mathcal{C})/G\), hence
\[\# X_e(\mathcal{C},G) \leq \sum_{\mathcal{P}\succeq\obj(\mathcal{C})/G} (2 a^c)^{m} \leq (2ma^c)^m.\tag{\ref{main_theorem_cat}.3}\]
If we drop the assumption that all objects of \(\mathcal{C}\) are \(\mathcal{C}\rtimes G\)-isomorphic we obtain the same inequality by applying Lemma~\ref{category_isomorphism_reduction}.

What remains is to show we can compute \(X_e(\mathcal{C},G)\) in \(n^{O(m)}\) time.
The counting argument we used simply consists of chaining together various explicit bijections from the lemmata in this section.
The reader should verifiy that for any such bijection \(\varphi\) computing \(\varphi(\sigma)\) takes \(n^{O(m)}\) time.
This is trivial for all steps except possibly Proposition~\ref{transitive_reduction_prop}, where we have to note it is easy to solve the equations (\ref{transitive_reduction_prop}.1) and (\ref{transitive_reduction_prop}.2).
Hence we get a \((2ma^c)^m\cdot n^{O(m)}\subseteq n^{O(m)}\) time algorithm, since both \(m\) and \(a\) are bounded by \(n\).
\qed

\subsection{Computing cyclic gradings of reduced rational algebras}\label{sec:AlgAbGrpThm}

Using Theorem~\ref{main_theorem_cat} we can finally prove Theorem~\ref{AlgAbGrpThm}.

\begin{remark}\label{rem:numberfield_isomorphsism}
Given two number fields \(K\) and \(L\) represented by structure constants we may compute \(\alpha\in K\) and \(\beta\in L\) together with minimal polynomials \(f_\alpha,f_\beta\in\Q[X]\) such that \(K=\Q(\alpha)\) and \(L=\Q(\beta)\), for example using Theorem~1.7 from \cite{algcomalg2018}.
Then \(K\cong L\) precisely when \(f_\alpha\) and \(f_\beta\) have the same degree and \(f_\alpha\) has a linear factor over \(L\).
Moreover, an isomorphism \(\varphi:K\to L\) is uniquely determined by the image of \(\alpha\), which must be a root of \(f_\alpha\) in \(L\).
Hence we may use Theorem~4.5 of \cite{AKfactoring} to compute the isomorphisms between two number fields in polynomial time with respect to the length of the input.
\end{remark}

\begin{Algorithm}\label{alg:cycl_grad}
Let \(E\) be an \(n\)-dimensional reduced \(\Q\)-algebra, let \(p\) be prime and let \(k\geq 0\).
\begin{enumerate}[topsep=0pt,itemsep=-1ex,partopsep=1ex,parsep=1ex]
\item Let replace \(k\) by \(\min\{k,\lfloor\log_p n\rfloor\}\). If \(k=0\) simply output the trivial grading of \(E\) and terminate.
Compute \(E'=E\tensor_\Q \Q(\mu_{p^k})\) as in Section~\ref{sec:X_to_grad_bij_computation}. 
Using Algorithm~7.2 in \cite{algcomalg2018}, we compute \(\spec E'\), \(E'/\m\) for \(\m\in\spec E'\) and the isomorphism \(\varphi:E'\to \prod_{\m\in\spec E'} E'/\m\). 
\item For each pair \(\m,\n\in\spec E'\) compute \(\Iso_{\Q(\mu_{p^k})\catt{-Alg}}(E'/\m,E'/\n)\) using Remark~\ref{rem:numberfield_isomorphsism}, such that we can construct the concrete category \(\mathcal{C}\) with \(\obj(\mathcal{C})=\spec E'\) and \(\Hom_\mathcal{C}(\m,\n)=\Iso_{\Q(\mu_{p^k})\catt{-Alg}}(E'/\m,E'/\n)\) as in Section~\ref{sec:wreath_conjugacy_computation}.
Let \(G=(\Z/p^k\Z)^*\) and compute its left action on \(\mathcal{C}\) via \(G\to\Aut(\Q(\mu_{p^k}))\to \Aut(E')\) as in Corollary~\ref{cor:cyclotomic_automorphisms}.
\item Using Theorem~\ref{main_theorem_cat} compute \(X_{p^k}(\mathcal{C},G)\). 
\item For each \(\sigma\in X_{p^k}(\mathcal{C},G)\) compute the corresponding \(f_\sigma \in \Aut(E')\) as in Proposition~\ref{aut_wreath_iso} using \(\varphi\). Then compute the grading \(\overline{E}_\sigma\) corresponding to \(f_\sigma\) by Theorem~\ref{cyclic_correspondence} using Proposition~\ref{prop:computing_eigenspaces} and output \(\overline{E}_\sigma\). 
\end{enumerate}
\end{Algorithm}

\noindent\textbf{Proof of Theorem~\ref{AlgAbGrpThm}.}
We need to show the correctness of Algorithm~\ref{alg:cycl_grad} and that its runtime is in \(n^{O(m)}\).
By Proposition~\ref{prop:basic_grad_props}, an efficient grading with a cyclic group of prime power order must have a group of size at most \(\dim_\Q(E)=n\). 
Hence we may replace \(k\) by \(k=\min\{k,\lfloor\log_p n\rfloor\}\) in step (1).
By the correctness of the various algorithms we use up to step (3), we correctly compute \(X_{p^k}(\mathcal{C},G)\).
By definition of \(X_{p^k}(\mathcal{C},G)\) each \(\sigma\in X_{p^k}(\mathcal{C},G)\) corresponds to an element \(f_\sigma\in\Aut(E')\) which is also an element of \(X_{p^k}(E)\) as defined in Definition~\ref{def:X}.
Then by Theorem~\ref{cyclic_correspondence} each such \(f_\sigma\) encodes a \(\mu_{p^k}\)-grading we obtain by computing the eigenspaces, which we all output.
Hence the algorithm is correct.

Step (1) takes time polynomial in \(n\) by Theorem~1.10 in \cite{algcomalg2018}.
For each pair \(\m,\n\in\spec A\) computation of \(\Iso(A/\m,A/\n)\) is polynomial in \(n\) as well by Remark~\ref{rem:numberfield_isomorphsism}, and thus so is step (2).
In step (3) computation of \(X_{p^k}(\mathcal{C},G)\) takes \(n^{O(m)}\) time by Theorem~\ref{main_theorem_cat}.
For each \(\sigma\) we require only \(n^{O(1)}\) time in step (4), meaning that Algorithm~\ref{alg:cycl_grad} terminates in \(n^{O(m)}\) time. \qed

\subsection{Computing universal abelian group gradings of reduced orders}\label{sec:AlgGridThm}

Let \(R\) be a reduced order and write \(E=R\tensor_\Z\Q\).
We have a natural inclusion \(R\to E\) which for any \(e\in\mathbb{S}\) induces an injection \(X_e(R)\to X_e(E)\), where \(X_e\) is as in Definition~\ref{def:X}.
Moreover, \(E\) is a reduced \(\Q\)-algebra with \(\#\spec E=\#\minspec E=\#\minspec R\) and \(\#\spec E'=\#\minspec R'\) with \(E'=E\tensor_\Q\Q(\mu_e)\) and \(R'=R\tensor_\Z\Z[\mu_e]\).
With a slight modification to Algorithm~\ref{alg:cycl_grad} we can compute the cyclic gradings of \(R\).

\begin{Algorithm}\label{alg:abel_ugrad_order}
Let \(R\) be a reduced order of rank \(r\).
\begin{enumerate}[topsep=0pt,itemsep=-1ex,partopsep=1ex,parsep=1ex]
\item Let \(\mathcal{P}=\{ p^{\lfloor\log_p r\rfloor} \,|\, p\text{ prime},\, p\leq r\}\). Let \(E=R\tensor_\Z\Q\) be the reduced \(\Q\)-algebra with the same structure constants as \(R\).
\item For each \(q\in\mathcal{P}\) compute \(X_q(E)\) using Algorithm~\ref{alg:cycl_grad} and then compute \(X_q(R)\) as the subset of \(X_q(E)\) of morphisms for which the corresponding matrices have integer coefficients. 
Let \(S=\bigcup_{q\in\mathcal{P}} X_q(R)\).
\item Using Proposition~\ref{prop:computing_eigenspaces} compute the set \(Z\) of all \(z=\{\zeta_s\}_{s\in S}\in \mu_\infty^S\) such that \(U_z=\bigcap_{s\in S} R(s,\zeta_s)\neq 0\) and for each \(z\in Z\) compute \(U_z\).
\item Compute the group \(\Upsilon=\langle Z\rangle\) and return \(\overline{U}=(R,\Upsilon,\{U_\upsilon\}_{\upsilon\in\Upsilon})\).
\end{enumerate}
\end{Algorithm}

\noindent\textbf{Proof of Theorem~\ref{AlgGridThm}:}
We show that Algorithm~\ref{alg:abel_ugrad_order} correctly computes a universal abelian group-grading of \(R\) in time \(n^{O(m)}\), where \(n\) is the length of the input \(R\) and \(m=\#\minspec R=\#\spec E\).
The inclusion \(R\to E\) induces an inclusion \(X_q(R)\to X_q(E)\) when we interpret \(X_q(R)\) as the automorphisms of \(X_q(E)\) that restrict to an automorphism of \(R'=R\tensor_\Z\Z[\mu_e]\).
These are precisely the matrices with integer coefficients, since \(R\) and \(E\) have the same structure constants.
By Theorem~\ref{UniGridGradThm} the reduced order \(R\) has a universal abelian group grading and by Remark~\ref{rmk:joint_uni_grad_exist} all joint gradings of gradings of \(R\) exist.
Thus \(\overline{U}\) is in fact a grading of \(R\).
Moreover, its computation can be done in time \(n^{O(m)}\) by Theorem~\ref{AlgAbGrpThm} and Proposition~\ref{prop:computing_eigenspaces}.
It remains to show that \(\overline{U}\) is in fact the universal abelian group-grading of \(R\).

By Theorem~\ref{thm:grad_is_morph} each efficient grading \(\overline{R}=(R,\Gamma,\{R_\gamma\}_{\gamma\in\Gamma})\) can be written as the joint grading of some gradings with cyclic groups of prime power order using the structure theorem for finite abelian groups applied to \(\Gamma\).
For each efficient abelian group grading the group \(\Gamma\) has exponent at most \(r\), so it has exponent dividing \(e=\prod_{q\in\mathcal{P}} q\).
Hence all cyclic gradings of prime power order that can occur are the gradings associated to the elements of \(S\).
Thus there exists a map \(f:\Upsilon\to\Gamma\) such that \(f_*\overline{U}=\overline{R}\).
Since \(\overline{U}\) is efficient by construction, this \(f\) is unique and \(\overline{U}\) is universal by Proposition~\ref{prop:universal_implies_efficient}.
Hence the algorithm is correct. \qed

\newpage
\section{Universal grid grading}

In this section we present an algorithm to compute the universal grid grading of a reduced order.

\subsection{Restricting decompositions of finitely generated abelian groups}

\begin{definition}
Let \(k\) be a ring, let \(B\subseteq A\) be \(k\)-modules and let \(\mathcal{A}=\{A_i\}_{i\in I}\) be a decomposition of \(A\).
We say \(\mathcal{A}\) {\em restricts} to \(B\) if \(\mathcal{B}=\{A_i\cap B\}_{i\in I}\) is a decomposition of \(B\), and we call \(\mathcal{B}\) the restriction of \(\mathcal{A}\) to \(B\).
\end{definition}

In this section we will write \(S/T=\{T\}\cup\{\{s\} \,|\, s\in S\setminus T \}\) for sets \(T\subseteq S\), which comes with a natural map \(S\to S/T\) for which the set of fibres is precisely \(S/T\).

\begin{proposition}\label{prop:universal_mod_restriction_exists}
Let \(k\) be a ring, let \(B\subseteq A\) be \(k\)-modules and let \(\mathcal{A}=\{A_i\}_{i\in I}\) be a decomposition of \(A\) with \(I\) finite.
Then there exists a map \(u:I\to U\) for some set \(U\) that such that \(u_*\mathcal{A}\) restricts to \(B\) and such that it is universal with respect to this property, i.e. each map \(f:I\to J\) such that \(f_*\mathcal{A}\) restricts to \(B\) factors uniquely through \(u\).
\begin{proof}
We apply induction to \(n=\#I\).
If \(n=0\) we clearly have \(u=\id_\emptyset\) as universal map.

Now assume that \(n>0\) and that universal maps exist for all decompositions of \(A\) with less than \(n\) components.
For \(b\in B\) we may uniquely write \(b=\sum_{i\in I}a_i\) with \(a_i\in A_i\) and we define the support \(\supp_\mathcal{A}(b)=\{i\in I\,|\, a_i\not\in B\}\) and weight \(w_\mathcal{A}(b)=\#\supp_\mathcal{A}(b)\) of \(b\).
If \(\mathcal{A}\) restricts to \(B\), or equivalently each element of \(B\) has zero weight, we take \(u=\id_I\) and are done.
Otherwise, there exists some \(b\in B\) with \(2 \leq w_\mathcal{A}(b)\), because having \(w_\mathcal{A}(b)=1\) is clearly impossible.
We may choose \(b\in B\) such that \(\supp_\mathcal{A}(c)\subsetneq\supp_\mathcal{A}(b)\) for some \(c\in B\) implies \(w_\mathcal{A}(c)=0\).
Consider \(J=I/\supp_\mathcal{A}(b)\) together with the natural map \(q:I\to J\) and let \(\mathcal{A}'=q_*\mathcal{A}\).
By the induction hypothesis there exists a universal map \(v:J\to U\) such that \(v_*\mathcal{A}'\) restricts to \(B\) and we claim that \(u=v\circ q\) is a universal map for \(\mathcal{A}\).
Let \(f:I\to K\) be such that \(f_*\mathcal{A}\) restricts to \(B\).
To show that \(f\) factors uniquely through \(u\) it suffices to show it factors uniquely through \(q\) by universality of \(v\), which is precisely the case when \(\#f(\supp_\mathcal{A}(b))=1\).
Since \(f_*\mathcal{A}=\mathcal{C}=\{C_k\}_{k\in K}\) restricts to \(B\) we may uniquely write \(b=\sum_{k\in K} c_k\) with \(c_k\in B\cap C_k\).
We may also uniquely write \(b=\sum_{i\in I}a_i\) with \(a_i\in A_i\) and note that \(c_k=\sum_{i\in f^{-1}k} a_i\).
Hence \(\{i\}\subseteq\supp_\mathcal{A}(c_{f(i)})\subseteq\supp_\mathcal{A}(b)\) for all \(i\in\supp_\mathcal{A}(b)\), so \(\supp_\mathcal{A}(c_{f(i)})=\supp_\mathcal{A}(b)\) by choice of \(b\).
As all \(c_{f(i)}\) have the same non-empty support they must be equal and thus \(f(i)=f(j)\) for all \(i,j\in\supp_\mathcal{A}(b)\), as was to be shown. 
\end{proof}
\end{proposition}

We now set out to turn Proposition~\ref{prop:universal_mod_restriction_exists} into an algorithm.
Note that by the structure theorem of finitely generated groups, we may represent a finitely generated abelian group \(A\) by a matrix \(M_A\in\Z^{n\times n}\) for some \(n\in\Z_{\geq0}\) such that \(A\cong \Z^n/M_A\Z^n\), which we call a {\em matrix representation} of \(A\).
A subgroup \(B\) of \(A\) is represented by a {\em generator matrix} \(G_B\in\Z^{n\times k}\) for some \(k\in\Z_{\geq 0}\) such that \(B=G_B\Z^k/M_A\Z^n\).
We will prove the following theorem.

\begin{theorem}\label{thm:abelian_decomposition_lift}
There is a deterministic algorithm that takes a quadruple \((n,A,\mathcal{A},B)\) with \(n\in\Z_{\geq0}\), with \(A\) a finitely generated group given in matrix representation \(M_A\), with \(\mathcal{A}=\{A_i\}_{i\in I}\) a decomposition of \(A\) with each \(A_i\) given by a generator matrix \(G_i\) and with \(B\) a subgroup of \(A\) given by a generator matrix \(G_B\), and computes a universal map \(u:I\to U\) such that \(f_*\mathcal{A}\) lifts to \(B\) in polynomial time with respect to the input.
\end{theorem}

\begin{lemma}\label{lem:intersect_check}
For matrices \(A\in\Z^{n\times a}\) and \(B\in\Z^{n\times b}\) for \(a,b,n\in\Z_{\geq 0}\) we may compute a basis for \(A\Z^a+B\Z^b\) and \(A\Z^a\cap B\Z^b\) and we can verify whether \(A\Z^a\subseteq B\Z^b\), in polynomial time with respect to the length of the input \((A,B)\). 
\begin{proof}
We may compute the image and kernel of integer matrices in polynomial time using using the image and kernel algorithm from \cite{KernelAlgorithm} and compute determinants using Theorem 6.6 from \cite{determinant}.
Then a basis for \(A\Z^a+B\Z^b\) can be computed as a basis for the image of the matrix \((A|B)\).
We may compute a basis matrix \(K\in\Z^{(a+b)\times k}\) of the kernel of \((A|B)\) for some \(k\in\Z_{\geq0}\).
A basis for \(A\Z^a\cap B\Z^b\) we then obtain as a basis for the image of \((A|0_{n\times b})\cdot K\), where \(0_{n\times b}\in\Z^{n\times b}\) is the all-zero matrix.
We may replace \(A\) by a basis for the image of \(A\) such that \(A\) becomes injective.
For the inclusion \(A\Z^a\subseteq B\Z^b\) we compute a basis matrix \(M\) for the image of \((\mathbb{I}_a|0_{a\times b})\cdot K\).
We have that \(A\Z^a\subseteq B\Z^b\) precisely when \(M\) is invertible, which we may check by computing its determinant.
\end{proof}
\end{lemma}

\begin{Algorithm}\label{alg:lifting_algorithm}
Let \(n\in\Z_{\geq0}\), let \(M_A\in\Z^{n\times n}\), let \(\mathcal{A}=\{A_i\}_{i\in I}\) a sequence of matrices with \(A_i\in\Z^{n\times k_i}\) for some \(k_i\in\Z_{\geq 0}\) and let \(G_B\in\Z^{n\times m}\) for some \(m\in\Z_{\geq 0}\).
\begin{enumerate}[topsep=0pt,itemsep=-1ex,partopsep=1ex,parsep=1ex]
\item Let \(J=I\) and replace \(G_B\) by a basis matrix for \(G_B\Z^{m}+M_A\Z^n\).
\item For \(S\subseteq J\) write \(A_S\) for the matrix with as columns the columns of the matrices \(A_j\) for \(j\in S\) and let \(k_S=\sum_{j\in S} k_j\). Let \(S=J\).
\item For each \(i\in J\), let \(T=S\setminus\{i\}\) and if \(G_B\Z^m\cap\sum_{j\in T} A_j\Z^{k_j} \not\subseteq \sum_{j\in T}(G_B\Z^m\cap A_j\Z^{k_j})\) then update \(S=T\).
\item If \(S=J\), return the natural map \(I\to J\) and terminate. Otherwise update \(J=J/S\) and \(\mathcal{A}=\{A_j\}_{j\in J}\) and go to step 2. 
\end{enumerate}
\end{Algorithm}

\noindent\textbf{Proof of Theorem~\ref{thm:abelian_decomposition_lift}:}
We show that Algorithm~\ref{alg:lifting_algorithm} satisfies the requirements. 
Correctness follows from the proof of Proposition~\ref{prop:universal_mod_restriction_exists}, where we note that in step 4 of the algorithm \(S\) is an non-empty set that is minimal with respect to inclusion for which there exists an element of \(B\) with support \(S\).
With regards to runtime, the only non-trivial step in the algorithm is step 3, which can be executed in polynomial time by Lemma~\ref{lem:intersect_check}.
Then, since the size of \(J\) strictly decreases in step 4 we enter this step at most \(\#I\) times.
Thus the algorithm runs in polynomial time. \qed

\subsection{Restricting gradings of orders}

\begin{definition}
Let \(k\) be a commutative ring, let \(R\) be a \(k\)-algebra and let \(S\subseteq R\) be a subalgebra.
We say a grading \(\overline{R}=(R,G,\mathcal{R})\) {\em restricts} to \(S\) when \(\mathcal{R}\) restricts to \(S\) and hence \(\overline{S}=(S,G,\mathcal{S})\) is a grading, where \(\mathcal{S}\) is the restriction of \(\mathcal{R}\) to \(S\).
We then call \(\overline{S}\) the restriction of \(\overline{R}\) to \(S\).
\end{definition}

\begin{lemma}\label{lem:universal_grading_restriction}
Let \(k\) be a commutative ring, let \(R\) be a \(k\)-algebra and let \(S\subseteq R\) be a subalgebra.
For each loose grid-grading \(\overline{R}=(R,G,\mathcal{R})\) there exists a universal morphism of grids \(u:G\to H\) such that \(u_*\overline{R}\) restricts to \(S\).
\begin{proof}
Let \(u_1:G\to I\) be the universal map such that \(u_{1*}\mathcal{R}\) restricts to \(S\), which exists by Proposition~\ref{prop:universal_mod_restriction_exists}.
Consider the functor \(F:\catt{Grd}\to\catt{Set}\) that sends \(\Gamma\) to the set of \(\Gamma\)-gradings of \(R\) of which \(u_{1*}\mathcal{R}\) is a refinement.
By Proposition~\ref{prop:refinement_universal} there then exists a universal grid \(H\) and map \(u_2:I\to H\) such that \(\overline{U}=u_{2*}u_{1*}\overline{R}\) is a grid grading.
Note that \(\overline{U}\) restricts to \(S\) by construction, and since \(\overline{R}\) is loose we have that \(u=u_2\circ u_1\) is a morphism of grids.
It follows from universality that \(u\) satisfies the additional requirements.
\end{proof}
\end{lemma}

\begin{proposition}\label{prop:grading_lift_algorithm}
There is an algorithm that takes a pair \((R,\overline{E})\), where \(R\) is a reduced order represented by structure constants and  \(\overline{E}\) is a loose grid-grading of \(E=R\tensor_\Z\Q\), the \(\Q\)-algebra with the same structure constants as \(R\), and which computes both a universal morphism of grids \(u\) such that \(u_*\overline{E}\) restricts to \(R\) and the restriction of \(u_*\overline{R}\) to \(R\), in polynomial time with respect to the length of the input. 
\begin{proof}
Let \(n\) be the rank of \(R\) and write \(\overline{E}=(E,G,\{E_g\}_{g\in G})\).
Compute integer bases for \(S_g=E_g\cap \Z^n\) and let \(S=\sum_{g\in G}S_g\).
Then the determinant \(d\) of the basis matrix of \(S\) satisfies \(d\Z^n\subseteq S\).
Applying a change of basis to the basis matrix of \(dR\) we may assume that \(S=\Z^n\).
Now apply Theorem~\ref{thm:abelian_decomposition_lift} to the decomposition \(\mathcal{S}=\{S_g\}_{g\in G}\) of the finitely generated abelian group \(\Z^n\) and the subgroup \(R\). 
This gives a universal map \(u_1:G\to I\) such that \(\{E_g\}_{g\in G}\) restricts to \(R\) and let \(\mathcal{T}=\{T_i\}_{i\in I}\) be this restriction.
Write \(\pi_i:R\to T_i\) be the natural projection.
Now construct a graph \(\mathcal{G}\) on the vertices \(I\) where \(\{i,j\}\in\binom{I}{2}\) is an edge if and only if there are \(a,b\in I\) such that \(\pi_i( T_a \cdot T_b ) \neq 0\) and \(\pi_j(T_a \cdot T_b) \neq 0\), which we may compute in polynomial time.
Then with \(C\) the set of connected components of \(\mathcal{G}\) we obtain a map \(u_2:I\to C\) such \(\mathcal{U}=u_{2*}\mathcal{T}\) is a pre-grading.
We can now easily define a grid structure on \(C\) such that \(\overline{U}=(R,C,\mathcal{U})\) is a loose grid grading and \(u=u_2\circ u_1\) satisfies the requirements.
\end{proof}
\end{proposition}

\subsection{Computing universal grid gradings of reduced orders}

We can now present our final algorithm.

\begin{Algorithm}\label{GridUgradNonrecursive}
Let \(R\) be a reduced order.
\begin{enumerate}[topsep=0pt,itemsep=-1ex,partopsep=1ex,parsep=1ex]
\item Let \(E=R\tensor_\Z\Q\) and compute \(\spec E\), \(E_P=\prod_{\m\in P}E/\m\) and maps \(\pi_P:E\to E_P\) for each \(P\subseteq\spec E\), as well as \(R_P=\pi_P(R)\) using Theorem~1.10 from \cite{algcomalg2018}. 
\item For each \(P\subseteq\spec E\) compute a universal abelian group grading \(\overline{A}_P\) of \(R_P\) using Algorithm~\ref{alg:abel_ugrad_order} and let \(\overline{E}_P\) be the corresponding grading of \(E_P\).
\item For each partition \(\mathcal{P}\) of \(\spec E\) such that \(\overline{E}_P\) is loose for each \(P\in\mathcal{P}\),
compute \(\overline{E}_\mathcal{P}=\coprod_{P\in\mathcal{P}}\overline{E}_P\) and the corresponding universal \(u_\mathcal{P}\) such that \(u_{\mathcal{P}*}\overline{E}_\mathcal{P}\) restricts to \(R\) using Proposition~\ref{prop:grading_lift_algorithm}, and let \(\overline{R}_\mathcal{P}\) be that restriction.
\item Return any grading \(\overline{R}_\mathcal{P}\) with the most homogeneous components among all partitions \(\mathcal{P}\).
\end{enumerate}
\end{Algorithm}

\noindent\textbf{Proof of Theorem~\ref{AlgGridThm2} for grids:}
We show that Algorithm~\ref{GridUgradNonrecursive} computes the universal grid grading of \(R\) in \(n^{O(m)}\) time, where \(m=\#\minspec R=\#\spec E\).
First we note that there are \(2^m\) subsets \(P\subseteq\spec E\) and at most \(m^m\) partitions \(\mathcal{P}\) of \(\spec E\).
For each \(P\subseteq\spec E\) we take at most \(n^{O(m)}\) time and for each \(\mathcal{P}\) only \(n^{O(1)}\) time.
Thus the total complexity becomes \(n^{O(m)}\), as was to be shown.

Let \(\overline{U}=(R,\Upsilon,\{U_{\upsilon}\}_{\upsilon\in \Upsilon})\) be a universal grid grading of \(R\) and let \(\overline{V}\) be the corresponding grading of \(E\).
Consider \(\prid(V_1)\), which induces a partition of \(\prid(E)\), which in turn induces a partition \(\mathcal{P}\) of \(\spec E\) by Corollary~\ref{cor:spectrum_and_idempotents}.
For \(P\in\mathcal{P}\) let \(\overline{V}_P=(\pi_P(E),\Upsilon,\{\pi_P(V_\upsilon)\}_{\upsilon\in\Upsilon})\), which form the factors of \(\overline{V}\) as in Lemma~\ref{lem:coproduct_idempotent_intermediate}.
Note that each \(\overline{V}_P\) restricts to \(R_P\) and let \(\overline{U}_P\) be this restriction.
Because \(\overline{A}_P=(R_P,G_P,\{A_{P,g}\}_{g\in G_P})\) is a universal abelian group grading of \(R_P\), there exists a unique morphism \(f_P:G_P\to \Upsilon\) such that \(f_{P*} \overline{A}_P=\overline{U}_P\) and thus \(f_{P*}\overline{E}_P=\overline{V}_P\).
Then the combined map \(f:\coprod_{P\in\mathcal{P}} G_P \to \Upsilon\) satisfies \(f_* \overline{E}_\mathcal{P}=\overline{V}\).
Since \(\overline{V}\) restricts to \(R\), there exists a  unique \(g\) such that \(f=g\circ u_\mathcal{P}\) and \(g_*\overline{R}_\mathcal{P}=\overline{U}\) by universality of \(u_\mathcal{P}\).
As \(\overline{U}\) is a universal grid grading of \(R\), so must be \(\overline{R}_\mathcal{P}\) as it is loose.
Hence a universal grid grading of \(R\) is among the possible gradings returned by the algorithm in step (4).

Assume there is some partition \(\mathcal{Q}\) of \(\spec E\) such \(\overline{E}_Q\) is loose for all \(Q\in\mathcal{Q}\) and such that \(\overline{R}_\mathcal{Q}\) has at least as many homogeneous components as \(\overline{R}_\mathcal{P}\). 
Then there exists a unique morphism \(g_*:\overline{R}_\mathcal{P}\to\overline{R}_\mathcal{Q}\) which is surjective on the underlying grid by universality and must be bijective as the grid of \(\overline{R}_\mathcal{Q}\) is at least as large as the grid of \(\overline{R}_\mathcal{P}\) and both are finite.
Since \(\overline{R}_\mathcal{Q}\) is loose it must be universal as well.
Thus the algorithm is correct.
\qed\\

We get the final part of Theorem~\ref{AlgGridThm2} as a corollary. \\

\noindent\textbf{Proof of Theorem~\ref{AlgGridThm2} for groups:}
We compute the universal grid grading \(\overline{U}=(R,\Upsilon,\{U_\upsilon\}_{\upsilon\in \Upsilon})\) using Algorithm~\ref{GridUgradNonrecursive}.
We may take \(G=\Upsilon^\catt{grp}\) by Proposition~\ref{prop:hierarchy} and a finite presentation of \(G\) and a map \(f:\Upsilon\to G\) can be obtained as in Lemma~\ref{lem:groupification_adjunction} in clearly polynomial time. \qed

\bibliographystyle{unsrt}
\bibliography{biblio}

\clearpage
\appendix
\section{Category theory}

We provide a brief summary of the category theoretical terminology we employ.
Most of this can be found in \cite{categories} albeit with some possible notational differences.

\begin{definition}
A {\em category} \(\mathcal{C}\) consists of a class \(\obj(\mathcal{C})\) of {\em objects} and a class \(\hom(\mathcal{C})\) of {\em morphisms}. Each \(f\in \hom(\mathcal{C})\) has a unique associated {\em source} \(s(f)\) and {\em target} \(t(f)\) in \(\obj(\mathcal{C})\), written \(f:s(f)\to t(f)\), and we write \(\Hom_{\mathcal{C}}(A,B)\) for the subclass of morphisms \(f\in\hom(\mathcal{C})\) such that \(s(f)=A\) and \(t(f)=B\).
For \(f\in \Hom_\mathcal{C}(A,B)\) and \(g\in\Hom_{\mathcal{C}}(B,C)\) we have a {\em composition} \(g\circ f\in \Hom_{\mathcal{C}}(A,C)\), and the operation \(\circ\) is associative.
Furthermore, \(\Hom_{\mathcal{C}}(A,A)\) contains a (unique) {\em identity} element \(\id_A\) such that for all \(f\in\hom(\mathcal{C})\) we have \(\id_A\circ f = f\) if \(t(f)=A\) and \(f\circ \id_A = f\) if \(s(f)=A\).
We say \(f\in\Hom_\mathcal{C}(A,B)\) is {\em invertible} if there exists some \(g\in\Hom_\mathcal{C}(B,A)\) such that \(g\circ f=\id_A\) and \(f\circ g = \id_B\).
A category \(\mathcal{D}\) is a {\em subcategory} of \(\mathcal{C}\) if \(\obj(\mathcal{D})\subseteq \obj(\mathcal{C})\), \(\hom(\mathcal{D})\subseteq\hom(\mathcal{C})\) and \(\mathcal{D}\) has the same composition law and identity elements as \(\mathcal{C}\).
Any \(A\in\obj(\mathcal{C})\) is called {\em initial} if for all \(B\in\obj(\mathcal{C})\) there exists a unique morphism \(f:A\to B\).
\end{definition}

We write \catt{Set} for the category of sets, \catt{Grp} for the category of groups and \catt{Ab} for the category of abelian groups.
For any commutative ring \(k\) we write \catt{\(k\)-Mod} for the category of \(k\)-modules and \catt{\(k\)-Alg} for the category of \(k\)-algebras.

\begin{definition}\label{def:cat_morphisms}
For a category \(\mathcal{C}\) we define the classes of {\em endomorphisms} and {\em isomorphisms} of \(\mathcal{C}\) as
\[ \text{end}(\mathcal{C}) = \{f\in\hom(\mathcal{C}) \,|\, s(f)=t(f)\} \quad \text{and}\quad \iso(\mathcal{C})=\{ f\in\hom(\mathcal{C}) \,|\, f \text{ is invertible} \} \]
respectively. We also write \(\aut(\mathcal{C})=\text{end}(\mathcal{C})\cap \iso(\mathcal{C})\) for the class of {\em automorphisms} of \(\mathcal{C}\).
For given \(A,B\in\obj(\mathcal{C})\) we additionally define the classes 
\[\text{End}_{\mathcal{C}}(A)=\Hom_{\mathcal{C}}(A,A),\quad \text{Iso}_{\mathcal{C}}(A,B)=\Hom_{\mathcal{C}}(A,B)\cap\iso(\mathcal{C} ), \quad \text{Aut}_{\mathcal{C}}(A)=\text{Iso}_{\mathcal{C}}(A,A).\]
\end{definition}

\begin{definition}
Let \(\mathcal{C}\) be a subcategory of \(\mathcal{D}\).
Then 
\begin{enumerate}\itemsep0em 
\item \(\mathcal{C}\) is {\em locally small} if \(\Hom_{\mathcal{C}}(A,B)\) is a set for all \(A,B\in\obj(\mathcal{C})\).
\item \(\mathcal{C}\) is {\em small} if it is locally small and \(\obj(\mathcal{C})\) is a set.
\item \(\mathcal{C}\) is {\em finite} if \(\obj(\mathcal{C})\) and \(\Hom_{\mathcal{C}}(A,B)\) are finite sets for all \(A,B\in\obj(\mathcal{C})\).
\item \(\mathcal{C}\) is a {\em wide subcategory of \(\mathcal{D}\)} if \(\obj(\mathcal{C})=\obj(\mathcal{D})\). 
\item \(\mathcal{C}\) is a {\em full subcategory of \(\mathcal{D}\)} if \(\Hom_\mathcal{C}(A,B)=\Hom_{\mathcal{D}}(A,B)\) for all \(A,B\in\obj(\mathcal{D})\). 
\item \(\mathcal{C}\) is a {\em monoid} if \(\#\obj(\mathcal{C})=1\).
\item \(\mathcal{C}\) is a {\em groupoid} if \(\hom(\mathcal{C})=\iso(\mathcal{C})\), i.e. all morphisms of \(\mathcal{C}\) are invertible.
\item subclasses \(X,Y\subseteq\obj(\mathcal{C})\) are {\em disconnected} if for all \(K\in X\) and \(L\in Y\) there are no morphisms \(K\to L\) or \(L\to K\).
\item a subclass \(Z\subseteq\obj(\mathcal{C})\) is {\em connected} if \(Z\neq\emptyset\) and for all disconnected subclasses \(X,Y\subseteq\obj(\mathcal{C})\) such that \(X\cup Y=Z\) we have \(X=Z\) or \(Y=Z\).
\end{enumerate}
\end{definition}

\begin{definition}
Let \(\mathcal{C},\mathcal{D}\) be categories.
A {\em functor} \(F:\mathcal{C}\to\mathcal{D}\) is a mapping that assigns to each \(A\in\obj(\mathcal{C})\) an object \(F(A)\in\obj(\mathcal{D})\) and to each morphism \(f\in\Hom_{\mathcal{C}}(A,B)\) a morphism \(F(f)\in\Hom_{\mathcal{D}}(F(A),F(B))\) such that
\[F(\id_A) = \id_{F(A)} \quad \text{and}\quad F(g\circ f) = F(g) \circ F(f)\]
for all \(A,B,C\in\obj(\mathcal{C})\), \(f\in\Hom_{\mathcal{C}}(A,B)\) and \(g\in\Hom_{\mathcal{C}}(B,C)\).
If \(\mathcal{C}\) and \(\mathcal{D}\) are small categories and \(A,B\in\obj(\mathcal{C})\) we write \(F_{A,B}:\Hom_\mathcal{C}(A,B)\to\Hom_\mathcal{D}(F(A),F(B))\) for the map given by \(f\mapsto F(f)\).
The class of all small categories forms a category \texttt{Cat}, in which the functors are the morphisms.
\end{definition}

\end{document}